\documentclass{amsart}
\title{Perturbations of  Planar Algebras}

\author[\large P\lowercase{aramita} D\lowercase{as},
  S\lowercase{hamindra} G\lowercase{hosh and} V\lowercase{ed}
  G\lowercase{upta}]{\bf \large P\lowercase{aramita} D\lowercase{as},
  S\lowercase{hamindra} K\lowercase{umar} G\lowercase{hosh and}
  V\lowercase{ed} P\lowercase{rakash} G\lowercase{upta}}

\thanks{Das and Ghosh were supported by K.~U.~Leuven BOF Research
  Grant OT/08/032 and Gupta was supported by ERC Starting Grant
  VNALG-200749}

\date{}

\address{Department of Mathematics, Katholieke Universiteit Leuven,
  Belgium}
\email{paramita.das@wis.kuleuven.be, shami.ghosh@wis.kulueven.be,
  ved.gupta@wis.kuleuven.be}

\usepackage{amsmath}
\usepackage{amsfonts}
\usepackage{latexsym}
\usepackage{amssymb}
\usepackage{mathrsfs}
\usepackage{amscd}
\usepackage{hyperref}
\usepackage{psfrag}
\usepackage{graphicx}
\usepackage{fullpage}
\usepackage[all]{xy}
\usepackage{verbatim}
\usepackage{enumerate}
\newtheorem{thm}{Theorem}[section]
\newtheorem{lem}[thm]{Lemma}
\newtheorem{cor}[thm]{Corollary}
\newtheorem{prop}[thm]{Proposition}
\newtheorem{defn}[thm]{Definition}
\newtheorem{rem}[thm]{Remark}

\newenvironment{pf}{\noindent{\em Proof:}}{\hfill $\Box$\vspace*{1mm}}
\numberwithin{equation}{section}
\numberwithin{figure}{section}

\newcommand{\comments}[1]{}
\newcommand{\ra}{\rightarrow}
\newcommand{\lra}{\longrightarrow}

\newcommand{\mbb}{\mathbb}
\newcommand{\mcal}{\mathcal}

\newcommand{\N}{\mathbb N}

\newcommand{\C}{\mathbb{C}}
\newcommand{\R}{\mathbb{R}}
\newcommand{\mscr}{\mathscr}
\newcommand{\vlon}{\varepsilon}

\newcommand{\uset}{\underset}
\newcommand{\oset}{\overset}
\newcommand{\oline}{\overline}
\newcommand{\uline}{\underline}
\newcommand{\vphi}{\varphi}
\newcommand{\ot}{\otimes}


\begin{document}
\global\long\def\vlon{\varepsilon}\global\long\def\ul#1{\underline{#1}}
\global\long\def\os#1#2{\overset{#1}{#2}}\global
\long\def\us#1#2{\underset{#1}{#2}}\global\long\def\ous#1#2#3{
\overset{#1}{\underset{#3}{#2}}}

\maketitle

\begin{abstract}

  We analyze the effect of pivotal structures (on a $2$-category) on
  the planar algebra associated to a $1$-cell as in \cite{Gho08} and
  come up with the notion of {\em perturbations of planar algebras by
    weights} (a concept that appeared earlier in Michael Burns' thesis
  \cite{Bur03}); we establish a one-to-one correspondence between
  weights and pivotal structures. Using the construction of
  \cite{Gho08}, to each bifinite bimodule over $II_1$-factors, we
  associate a {\em bimodule planar algebra} in such a way that
  extremality of the bimodule corresponds to sphericality of the
  planar algebra. As a consequence of this, we reproduce an extension
  of Jones' theorem (\cite{Jon}) (of associating `subfactor planar
  algebras' to extremal subfactors). Conversely, given a bimodule
  planar algebra, we construct a bifinite bimodule whose associated
  bimodule planar algebra is the one which we start with, using
  perturbations and Jones-Walker-Shlyakhtenko-Kodiyalam-Sunder method
  of reconstructing an extremal subfactor from a subfactor planar
  algebra.  The perturbation technique helps us to construct an
  example of a family of non-spherical planar algebras starting from a
  particular spherical one; we also show that this family is
  associated to a known family of subfactors constructed by Jones.

\end{abstract}


\section{Introduction}
In the pioneering and celebrated work \cite{Jon83} of Jones, the
theory of subfactors saw a new opening with, among other ideas and
results, the introduction of the concepts of index and the tower of
basic construction for subfactors, which, over the years, had various
applications in the understanding of $II_1$-factors, knot theory,
quantum groups, TQFTs and other fields.  There have been a lot of
pathbreaking works in this theory since its initiation - see, for
instance, \cite{Jon83, PP86, PP88, Pop90a, Pop94, Pop95}.

Further, an important aspect of this paper of Jones was the evolution
of an invariant called the standard invariant of the subfactor, which
basically consists of a grid of finite dimensional $C^*$-algebras with
some rich structure. It has turned out that, among other invariants,
the standard invariant of a finite index subfactor is its most
important invariant. An instance to justify this claim is that, for
certain `good' family of subfactors (namely, the amenable ones, see
\cite{Pop94}), their standard invariants turn out to be complete
invariants. As such, it was motivating enough for people to work on
obtaining a better understanding of this invariant. Sorin Popa (in
\cite{Pop95}) gave an algebraic axiomatization of the standard
invariant as a grid of finite dimensional $C^*$-algebras, which he
called a {\em standard $\lambda$-lattice}. Conversely, given such a
$\lambda$-lattice, he constructed an {\em extremal} subfactor whose
standard invariant is the $\lambda$-lattice which he started
with. Ocneanu also came up with a `group-like' structure on the
standard invariant and called them {\em paragroups}. Subsequently,
Vaughan Jones (in \cite{Jon}) developed a very effective pictorial
reformulation of the standard invariant which he called {\em planar
  algebra}, and associated a planar algebra satisfying certain natural
conditions (referred as {\em subfactor planar algebra}) to any
extremal subfactor.  In the converse direction, starting from a
subfactor planar algebra Jones reconstructed a subfactor whose
associated planar algebra is isomorphic to the given one, using Popa's
characterization of $\lambda$-lattices.  Later, in \cite{Pop02}, Popa
extended his correspondence to {\em generalized $\lambda$-lattices} on
one hand and finite index subfactors (not necessarily extremal) on the
other. In 2003, Michael Burns in his thesis (see \cite{Bur03}),
established a similar correspondence replacing generalized
$\lambda$-lattices with planar algebras satisfying appropriate
conditions (without any assumption of sphericality). Very recently,
reconstruction of extremal subfactor from a subfactor planar algebra
has also been performed using random matrix and free probability
techniques in \cite{GJS} followed by simpler treatments using planar
algebra machinery in \cite{JSW08} and \cite{KS08}. Planar algebra
techniques have recently found applications in developing new methods
of constructing certain class of subfactors as well.

We now present a brief outline of the motivation and the results that
brought this article into existence.
\vspace*{2mm}

 \noindent  (I) Our motivation stemmed
solely from the investigation of the following:
\begin{enumerate}[(a)]
\item A construction of a planar algebra starting from a $1$-cell in a
  pivotal $2$-category was given by the second named author in his
  thesis (see \cite{Gho08}). This construction was purely algebraic,
  with the description of the action of tangles being given in terms
  of {\em graphical calculus of morphisms}, analogous to the ones used
  in \cite{Kas}.  However, the actual manifestation of the pivotal
  structure in the planar picture remained unclear and required
  further analysis.
\item In the operator algebra context,
a nice prototype is the $2$-category of bifinite bimodules over
$II_1$-factors. So, one would like to investigate the planar algebras
obtained using the method in \cite{Gho08} from a bifinite bimodule
${_A}{\mcal H}_B$ where $A$ and $B$ are $II_1$ factors.

\item  Another question in this context is whether the
Jones' planar algebra associated to an extremal finite index subfactor
$N \subset M$, is isomorphic to the planar algebra coming from the
bimodule ${_N}{L^2(M)}_M$.
\end{enumerate}
 In this paper, we make an attempt to answer
these and other natural questions, the answers of which we list below:
\vspace*{2mm}

\noindent (II) Summary of the main results:
\begin{enumerate}[(a)]
\item We answer the question (I)(c), that is, we show that one does
  not always get back Jones' subfactor planar algebraâ from the
  construction in \cite{Gho08} unless the right pivotal structure is
  chosen.
\item In order to find out the exact dependence of the planar algebra
  from the construction in \cite{Gho08} on the pivotal structure, we
  come up with the concept of weights of a planar algebra and
  perturbations of planar algebras by weights; later, we realized that
  such objects also appeared in Michael Burns' thesis to prove Jones'
  theorem for non-extremal subfactors. Finally, we establish a
  one-to-one correspondence between weights on a planar algebra and
  pivotal structures on the associated $2$-category
\item To each bifinite bimodule over $II_1$ -factors, using the
  construction in \cite{Gho08}, we associate a {\em bimodule planar
    algebra} in such a way that extremality of the bimodule
  corresponds to sphericality of the planar algebra. Moreover, this
  also shows that bimodules with different left and right dimensions,
  gives the right platform to investigate planar algebras with
  different modulii coming from shaded and unshaded loops.

\item Conversely, given a bimodule planar algebra, we construct a
  bifinite bimodule whose associated bimodule planar algebra is the
  one that we start with, using perturbations and
  Jones-Walker-Shlyakhtenko-Kodiyalam-Sunder method of reconstructing
  an extremal subfactor from a subfactor planar algebra.

\item We give explicit construction of examples of non-spherical
  planar algebras; more precisely we show that the perturbation class
  of the diagonal planar algebra with respect to the free group
  $\mathbb{F}_2$, generated by two free generators and trivial cocycle
  contains a continuum of non-spherical unimodular bimodule planar
  algebras with index greater than $4$; we also prove that this family
  is associated to a known family of subfactors constructed by Jones.
\end{enumerate}

\noindent (III) Some nice consequences:
\begin{enumerate}[(a)]
\item As a consequence of (II)(c), we reproduce an extension of Jones'
  theorem (of associating `subfactor planar algebras' to extremal
  subfactors). This was proved earlier by Michael Burns in his thesis
  \cite{Bur03}; the reconstruction of a non-extremal subfactor had
  appeared in \cite{Pop02}.

\item We show that the perturbation class of a bimodule planar algebra
  contains a unique spherical unimodular bimodule planar algebra which
  can also be characterized by the minimality of its index. Subfactor
  version of such results had appeared in the works of Hiai and Popa.
\end{enumerate}
All results in this article are derived using standard facts on
bimodules, subfactors and planar algebras, which can be found, for
instance, in \cite{Bis97,Jon,JS97,PP86,PP88,Pop90a,Pop94,Pop95}.
\vspace{0.2cm}

We now briefly describe the organization of this paper. Section
\ref{prelim} serves as a quick recollection of various definitions,
standard facts and basic aspects of planar algebras, pivotal
bicategories and the bicategory of bifinite bimodules.

In Section \ref{pa}, we define weight of a planar algebra $P$ and
perturbation of $P$ by the weight. A planar algebra with modulus
$(\delta_-, \delta_+)$ can be {\em normalized}, that is, perturbed
with an appropriate scalar weight to get a {\em unimodular} planar
algebra (that is, having $\delta_- = \delta_+$) although the {\em
  index} (:= the product of the $\delta$'s) remains unchanged. If the
actions of the $0$-tangles in the normalization are invariant under
spherical isotopy, then the planar algebra is called spherical; this
is a slight modification of Jones' definition of sphericality in order
to accommodate non-unimodular planar algebras. At the very end of this
section, we make few immediate observations involving perturbations,
$*$-structures and positivity in planar algebras.

In the first part of Section \ref{weights}, we associate a strict
$2$-category to a planar algebra and show that weights of the planar
algebra are in one-to-one correspondence with pivotal structures on
the associated $2$-category. Conversely, if we start with a bicategory
with two pivotal structures, then the planar algebras obtained from
any $1$-cell using \cite{Gho08} method, are perturbations of each
other.  This section ends with a quick recollection (from
\cite{Gho08}) of the method of associating a planar algebra to an
$1$-cell in a strict $2$-category.

Section \ref{bimod-pa} is an omnibus section and is the crux of this
paper. In this section, we first formalize what we mean by a bimodule
planar algebra; then, following the above-mentioned method of
\cite{Gho08}, we associate a bimodule planar algebra to each bifinite
bimodule such that the extremality of the bimodule exactly corresponds
to the sphericality of its associated planar algebra. This, in turn,
provides an extension to Jones' theorem \cite[Theorem $4.2.1$]{Jon},
that is, we associate a unimodular bimodule planar algebra to an
arbitrary finite index subfactor. Such extension was also obtained by
Burns; however, our techniques are completely independent and rely on
simple graphical calculus of morphisms in the pivotal $2$ category.
In the converse direction, given any bimodule planar algebra, we
obtain a bifinite bimodule whose associated bimodule planar algebra is
isomorphic to the one that we started with, through the application of
perturbations and following the strategy of \cite{JSW08,KS08}.

In the first part of Section \ref{eg}, we show that the perturbation
class of every bimodule planar algebra contains a unique spherical
unimodular bimodule planar algebra which can also be characterized by
the property of having the minimal index in the perturbation
class. Minimizing indices of conditional expectations onto a subfactor
already appeared in the literature in the work of Hiai (in
\cite{Hia88}) and then Popa (in \cite{Pop94}); we are now able to
connect this circle of ideas with our notion of perturbation of planar
algebra. In the second part, we construct concrete examples of
nonspherical planar algebras purely algebraically; more precisely we
show that the perturbation class of the diagonal planar algebra with
respect to the free group $\mbb{F}_2$, generated by two free
generators and trivial cocycle contains a continuum of non-spherical
unimodular bimodule planar algebras with index greater than $4$. As
suggested by Jones, we prove that these planar algebras are isomorphic
to the ones associated to the (non-extremal) subfactors that he
constructed in \cite{Jon83} in order to prove that every index greater
than $4$ is realized.

In the final section, we discuss some questions pertaining to
perturbations and weights of a planar algebra.


\section{Preliminaries}\label{prelim}
This section is mainly a recollection of various definitions, standard
facts and setting up of notations which will be used in the subsequent
sections.
\subsection{Planar algebras}
\begin{figure}[h]
 \psfrag{k}{$k$}
 \psfrag{2k}{$2k$}
 \psfrag{e}{$\vlon$}
 \psfrag{-e}{$-\vlon$}
 \psfrag{M}{$ M_{\vlon k} =$}
 \psfrag{m}{$(\vlon k,\vlon k)\rightarrow \vlon k$}
 \psfrag{1ek}{$1_{\vlon k} =$}
\psfrag{1}{$\emptyset \rightarrow \vlon k$}
 \psfrag{RI}{$RI_{\vlon k} =$}
  \psfrag{ri}{$\vlon k \rightarrow \vlon (k+1)$}
 \psfrag{LI}{$LI_{\vlon k}=$}
  \psfrag{li}{$\vlon k \rightarrow -\vlon (k+1)$}
 \psfrag{Id}{$I_{\vlon k}=$}
\psfrag{id}{$\vlon k \rightarrow \vlon k$}
\psfrag{E}{$E_{\vlon (k+1)}=$}
\psfrag{jp}{$\emptyset \rightarrow \vlon (k+2)$}
 \psfrag{RE}{$RE_{\vlon (k+1)} =$}
 \psfrag{re}{$\vlon(k+1) \rightarrow \vlon k$}
 \psfrag{LE}{$LE_{\vlon (k+1)} =$}
 \psfrag{le}{$\vlon(k+1) \rightarrow -\vlon k$}
     \includegraphics[scale=0.2]{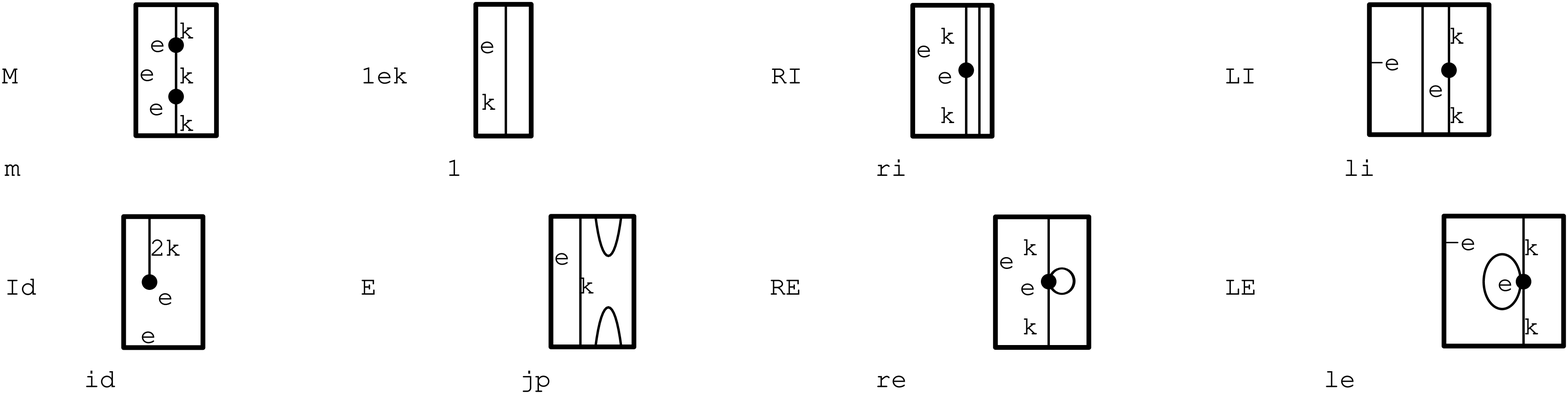}\vspace*{-3mm}
 \caption{Generating tangles.}\label{tangles}
 \end{figure}

 Since its inception in \cite{Jon}, the formalism of planar algebras
 has undergone gradual modifications - see, for instance, \cite{Jon,
   Jon00, KS04, Gho08}. The starting ingredient for defining a planar
 algebra is the operad of tangles. A {\em tangular diagram} $T$
 consists of a subset $D_0$ (referred as the {\em external disc}) of
 $\mbb{R}^2$, homeomorphic to the unit disc along with: (a) finitely
 many (possibly none) non-intersecting subsets $D_1,\, \ldots,\, D_b$
 (referred as {\em internal discs}) in the interior of $D_0$, each of
 which is also homeomorphic to the unit disc, (b) the boundary of each
 disc (internal or external) having even number of marked points
 numbered clockwise, (c) non-intersecting paths (called {\em
   strings}) in $D_0\setminus \left[\uset{i=1}{\oset{b}{\cup}} Int(
   D_i)\right]$, which are either loops or meet the boundaries of the
 discs exactly at two distinct marked points in such a way that every
 marked point is an endpoint of a string and (d) a checker-board
 shading on the connected components of $ Int(D_0)
 \setminus\left[\left( \uset{i=1}{\oset{b}{\cup}} D_i \right) \cup
   \{\text{strings}\} \right]$. We will usually indicate the
 checker-board shading and the numbering of the marked points simply
 by putting a $-$ (resp., $+$) sign in the shaded (resp., unshaded)
 connected component near the boundary segment between the last and
 first marked points. Such a sign along with half the number of marked
 points on a disc is called its {\em color}. A {\em tangle} is the
 class of a tangular diagram under the equivalence of {\em planar
   isotopy} (preserving the shading and the distinguished boundary
 components). If a tangle $T$ has $b > 0$ (resp., no) internal disc(s)
 and the color of the disc $D_i$ is $\vlon_i k_i$ for $0 \leq i \leq
 b$, then the tangle is usually expressed as $T: (\vlon_1 k_1, \ldots,
 \vlon_b k_b) \ra \vlon_0 k_0$ (resp., $T: \emptyset \ra \vlon_0
 k_0$). See Figure \ref{tangles} for illustrations. Before we proceed
 further, we fix some notations.
 \begin{enumerate}
\item We will consider the natural binary operation on $\{-,+\}$
given by $++ := +$, $+-:=-$, $-+:=-$ and $--:=+$.
\item In a tangle, we will replace (isotopically) parallel strings by
a single strand labelled by the number of strings, and an internal
disc with color $\vlon k$ will be replaced by a bold dot with the
sign $\vlon$ placed at the angle corresponding to the distinguished
boundary component of the disc. For example,\\
\psfrag{replace}{will be replaced by}
\psfrag{2}{$2$}
\psfrag{4}{$4$}
\psfrag{e}{$\vlon$}
\includegraphics[scale=0.2]{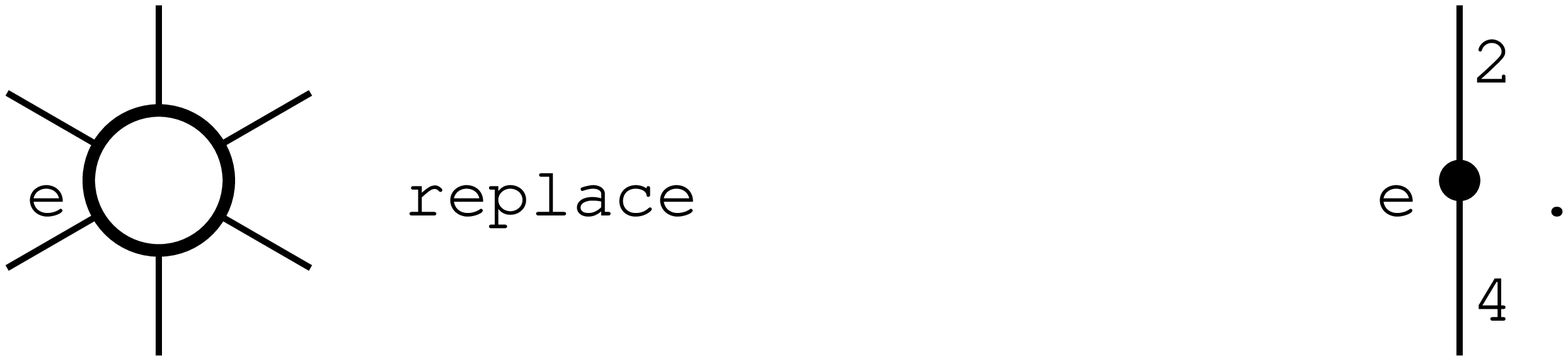}
\item We will denote the set of all possible colors of discs in
  tangles by $Col:= \left\{ \vlon k: \vlon \in \{+,-\},\, k \in
  \mbb{N}_0\right\}$ where $\N_0 := \N \cup \{ 0 \}$.
 \end{enumerate}
 There is a natural notion of {\em composition of tangles} under which
 the tangles listed in Figure \ref{tangles} generate the whole operad
 of tangles. (See \cite{KS04}.)

A {\em planar algebra} $P$ is a `representation' of the operad of
tangles, that is, it consists of complex vector spaces $\{P_{\vlon k}:
\vlon k \in Col \}$ and for every tangle $T : (\vlon_1 k_1,\vlon_2
k_2, \, \ldots, \vlon_b k_b) \ra \vlon_0 k_0$ (resp., $T : \emptyset
\ra \vlon_0 k_0 $), there exists an action of $T$ given by a
multi-linear map $P_T : \times_{i = 1}^b P_{\vlon_i k_i} \ra
P_{\vlon_0 k_0}$ (resp., a vector $P_T \in P_{\vlon_0 k_0}$) such that
the action preserves (i) composition and (ii) identity (that is,
$P_{I_{\vlon k}} = id_{P_{\vlon k}}$). Note that $\{ P_{\vlon k} \}_{k
  \in \N_0}$ has a unital filtered algebra structure with
multiplication, unit and inclusion given by the actions of $M_{\vlon
  k}$, $1_{\vlon k}$ and $RI_{\vlon k}$, respectively.
\begin{defn}A planar algebra $P$ is said to
\begin{enumerate}
\item be connected (resp., finite dimensional) if $dim ( P_{\pm 0} ) =
  1$ (resp., $dim ( P_{\vlon k} ) < \infty$ for all $\vlon k \in
  Col$).
\item have modulus $(\delta_-, \delta_+) \in \C^2$ if
  \psfrag{e}{$\pm$} \psfrag{P}{$P$} \psfrag{rhs}{$= \delta_\pm
    P_{1_{\pm 0}} = \delta_\pm 1_{P_{\pm 0}}$ and in this case, the
    scalar $\delta_+ \delta_-$ is call-} $\!\!\!\!\!\!\!
  \vcenter{\includegraphics[scale=0.2]{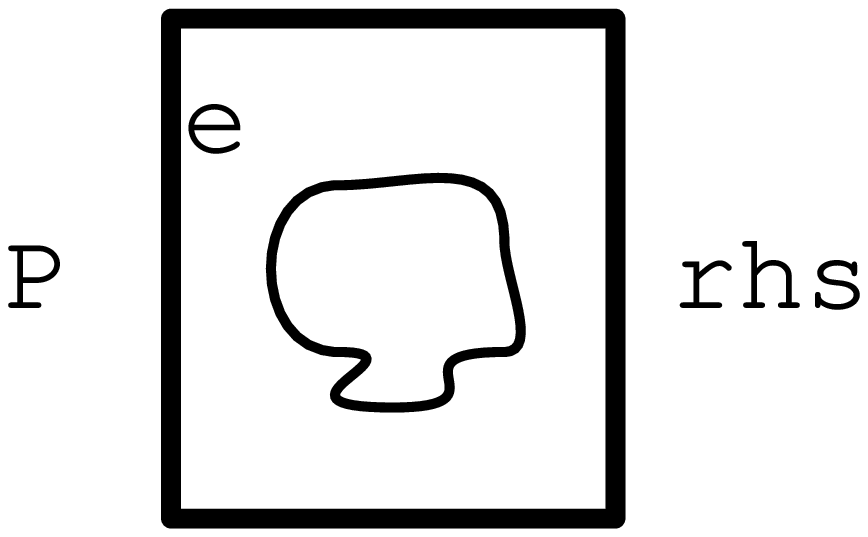}}$
  ed the index of $P$.
\item be unimodular if it has modulus $(\delta_-, \delta_+)$ such that
  $\delta_- = \delta_+$.
\item be a $\ast$-planar algebra if there exists a conjugate-linear
  involutions $\ast:P_{\vlon k} \rightarrow P_{\vlon k}$ for all
  $\vlon k \in Col$ satisfying the $\ast$-condition:
\begin{equation}\label{*-condition}
P_T(x_1,\, \ldots,\, x_b)^* = P_{T^*}(x_1^*,\, \ldots,\, x_b^*) \quad
(resp.,\ (P_T)^* = P_{T^*} )
\end{equation}
for each tangle $T: (\vlon_1 k_1\times \cdots \times \vlon_b k_b) \ra
\vlon_0 k_0$ (resp., $T : \emptyset \ra \vlon_0 k_0$) and $x_i \in
P_{\vlon_i k_i}$, $1 \leq i \leq b$, where the adjoint $T^*$ of a
tangle $T$ is obtained by reflecting it about a horizontal line
keeping the shading and distinguished boundary components intact.
\item be $C^\ast$-planar algebra if each $P_{\vlon k}$ is a
  $C^*$-algebra such that its multiplication is the same as that
  induced by $M_{\vlon k}$ and $P$ forms a $\ast$-planar algebra with
  respect to the $\ast$.
\item be irreducible if $dim ( P_{+1} ) = 1$ (equivalently, $dim ( P_{-1} ) = 1$)
\end{enumerate}
\end{defn}
\noindent Note that a connected planar algebra $P$ has modulus, and
two canonical {\em picture traces} $ P_{TR^{r}_{\vlon k}} : P_{\vlon
  k} \ra P_{\vlon 0}\cong \mbb{C} $ and $P_{TR^{l}_{\vlon k}}:
P_{\vlon k} \ra P_{(-)^k \vlon 0} \cong \mbb{C}$ induced by the trace tangles $TR^r_{\vlon k} := $ 
\psfrag{k}{$k$}
\psfrag{e}{$\vlon$}
\psfrag{ek}{$(-)^k \vlon$}
\includegraphics[scale=0.20]{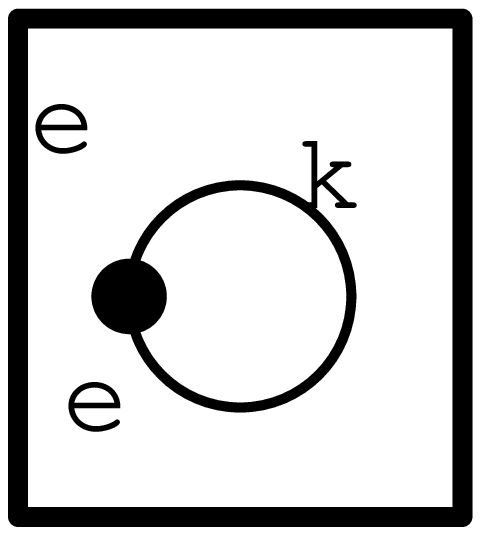}
and $TR^l_{\vlon k} := $
\includegraphics[scale=0.20]{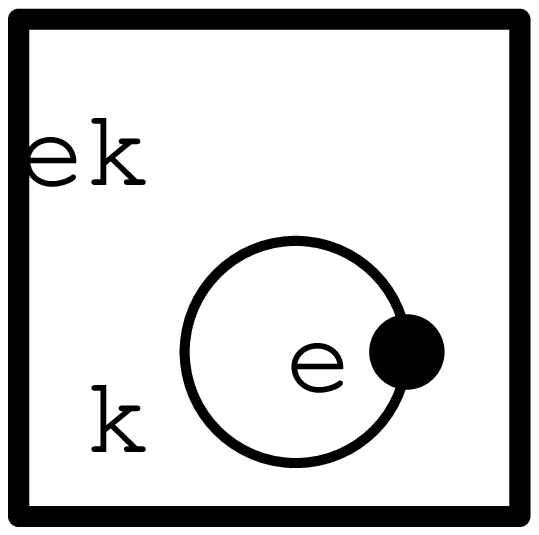} .
A connected $C^*$-planar algebra $P$ is said to be {\em positive} if
the canonical picture traces are positive definite.

We now recall the definition of the {\em $n$-th dual} of a planar algebra $P$ with modulus $(\delta_- , \delta_+)$, denoted by $\lambda_n (P)$. For a tangle $T$, let $\lambda_n(T)$ be the tangle obtained from $T$ by attaching $n$ parallel string on the disc $D$ (internal or external) enclosing closely the sign $\vlon_D$ near the distinguished boundary component of the disc.

\noindent {\em Vectors spaces:} For all colors $\vlon k$, $\lambda_n (P)_{\vlon k} := Range( P_{\lambda_n ( I_{\vlon k} ) } )$.

\noindent {\em Action of tangles:} For all tangle $T$, $\lambda_n (P)_T := \left[ \us{D \in \{\text{internal disc(s) of } T \} }{\prod} \left( \ous{n}{\prod}{l=1} \delta_{(-)^l \vlon_D} \right)^{-1} \right] P_{\lambda_n (T)}$.

Note that $\lambda_0 (P) = P$, $\lambda_n ( \lambda_m (P) ) = \lambda_{m+n} (P)$ and $\lambda_n(P)$ has modulus $(\delta_{(-)^{n+1}} , \delta_{(-)^n})$ for all $m , n \in \N_0$.

For any mathematical category, studying the morphisms in that category
is very crucial. A {\em morphism} $\varphi:P \rightarrow Q$ for two
planar algebras $P$ and $Q$, is a collection of linear maps
$\varphi_{\vlon k}: P_{\vlon k} \ra Q_{\vlon k}$, $\vlon k \in Col$
which are equivariant with the action of tangles (that is, $\vphi
\circ P_T = Q_T \circ \vphi$ or $\vphi \circ P_T = Q_T$ according as
$T$ has at least one internal disc or none). Given such a double
sequence of linear maps, it is not always necessary to verify its
equivariance with every tangle. For instance, it suffices to do the
same only for the actions of a generating set of tangles mentioned in
Figure \ref{tangles}.
\begin{rem}\label{planar-morphism}
Let $P$ and $Q$ be connected planar algebras with non-zero
moduli. Then, a linear map $ \varphi : P \ra Q$ is a planar algebra
morphism if it is equivariant with the actions of any of the following
sets of tangles:
\begin{enumerate}
\item $\us{k \in \N_0}{\cup} \{M_{\pm k}, RI_{\pm k}, LI_{\pm k},
  E_{\pm 1}\}$\comments{used in Theorem 5.2 (i)}
\item $\us{k \in \N_0}{\cup} \{M_{+k}, RI_{+k}, LI_{-k}, E_{+(k+1)}, LE_{+ (k+1)} \}$
\comments{\item $\{M_{+k}, RI_{+k}, LE_{+(k+1)}, RE_{+ (k+1)}, E_{+(k+1)} : k \in
  \N_0 \}$}
\end{enumerate}
\end{rem}
Sometimes, in order to obtain morphisms between planar algebras, it is
even enough to obtain a morphism on the `positive parts'. More
precisely, we have the following lemma.
\begin{lem}\label{pa-positive-isom}
Let $P$ and $Q$ be planar algebras having the same pair of non-zero
moduli, say $(\delta_-, \delta_+)$, and suppose there exist linear
maps $\varphi_{+k}: P_{+k} \ra Q_{+k}$, $k \geq 0$,
equivariant with the action of tangles with discs only of positive
colors. Then there exists a unique planar algebra morphism
$\tilde{\varphi}: P \oset{\sim}{\ra}Q$ such that $ \varphi_{+k} =
\tilde{\varphi}_{+k}$ for all $k \geq 0$. Moreover, if $\varphi$ is an
isomorphism, so is $\tilde \varphi$.
\end{lem}
 \begin{pf}
Consider $\tilde{\varphi} : P \ra Q$ given by $\tilde{\varphi}_{+k} =
\varphi_{+k}$ and
\[
P_{-k} \ni x \oset{\tilde{\varphi}_{-k}}{\longmapsto} \delta_-^{-1}
Q_{LE_{+(k+1)}} \circ \varphi_{+(k+1)} \circ P_{LI_{-k}} (x) \in Q_{-k}
\]
for all $k \geq 0$. It is straight forward to check that $\tilde
\varphi$ is equivariant with the action of tangles of all types (not
necessarily having positive colors for all discs). The uniqueness
follows from the definition.

For the second part, since $\varphi$ is equivariant with the action of
$E'_{+(k+1)} := LI_{-k} \circ LE_{+(k+1)}$, we have $\varphi_{+(k+1)}
(Ran\, P_{E'_{+(k+1)}}) \subseteq Ran\, Q_{E'_{+(k+1)}} $, with
equality if $\varphi$ is surjective; also, $\varphi_{-k}$ is injective
if so is $\varphi$ because $Q_{LE_{+(k+1)}}$ restricted to $ Ran\, Q_{
  E'_{+(k+1)} } $ is injective. Thus, $\tilde{\varphi}$ is an
isomorphism if so is $\varphi$.
\end{pf}
Let $P$ be a planar algebra. Then, by a $P$-{\em labelled} (resp.,
{\em semi-labelled}) tangle, we mean a tangle whose all (resp., some)
internal discs are labelled by elements of $P$ such that an internal
disc of color $\vlon k$ is labelled by an element of $P_{\vlon k}$.
For simplicity, we will replace a $P$-labelled internal disc by a bold dot as before with the label being placed at the angle corresponding to the distinguished boundary component of the disc.
\subsection{Bicategories}
In this subsection, for the sake of completeness, we recall the notion
of bicategories and structures of rigidity and pivotality on
them which will be used later. Most of the
material in this subsection can be found in any standard textbook on
bicategories.
\begin{defn}
A bicategory $\mathcal{B}$ consists of:
\begin{itemize}
\item a class $\mathcal{B}_{0}$ whose elements are called $0$-cells;
\item for $\alpha, \beta \in \mathcal{B}_{0}$, a category
  $\mathcal{B}(\alpha,\beta)$ whose object $X$ shall be called
  $1$-cells and denoted by $\alpha\overset{X}{\ra}\beta$, and whose
  morphism $X_{1}\overset{f }{ \ra }X_{2}$ from
  $\alpha\overset{X_{1}}{\ra}\beta$ to
  $\alpha\overset{X_{2}}{\ra}\beta$ shall be called $2$-cells;
\item for $\alpha, \beta, \gamma \in \mathcal{B}_{0}$, there exists a
  functor $ \otimes :\mathcal{B}(\beta,\gamma)\times
  \mathcal{B}(\alpha,\beta) \rightarrow \mathcal{B}(\alpha,\gamma)$;
\item Associativity constraint: For each triple
$\alpha\overset{X}{ \ra }\beta$,
$\beta\overset{Y}{\ra }\gamma$,
$\gamma\overset{Z}{ \ra }\delta$ of $1$-cells, there
exists an isomorphism $(Z\otimes Y)\otimes X\overset{\alpha
_{Z,Y,X}}{\longrightarrow } Z \otimes (Y\otimes X)$ in
$Mor(\mathcal{B}(\alpha,\delta))$;
\item Identity object: for each $0$-cell $\alpha$ , there
  exists a $1$-cell $\alpha \stackrel{1_{\alpha}}{\ra} \alpha$ (called
  the identity on $\alpha$);
\item Unit Constraint: for each $1$-cell
  $\alpha\overset{X}{\ra }\beta$, there exist isomorphisms
  $1_{\beta}\otimes X\overset{\lambda _{X}}{ \ra }X$ and
  $X\otimes 1_{\alpha}\overset{\rho _{X}}{\ra} X $ in
  $Mor(\mathcal{B}(\alpha,\beta))$
\end{itemize}
 such that $\alpha _{Z,Y,X}$, $\lambda _{X}$ and $\rho _{X}$ are
  natural in $Z$, $Y$ and $X$, and satisfy the pentagon and the
  triangle axioms (which are exactly similar to the ones in the
  definition of a tensor category).
\end{defn}

On a bicategory $\mathcal{B}$, one can perform the operation $op$
(resp., $co$) and obtain a new bicategory $\mathcal{B}^{op}$ (resp.,
$\mathcal{B}^{co}$) by setting (i)
$\mathcal{B}^{op}_0=\mathcal{B}_0=\mathcal{B}^{co}_0$, (ii)
$\mathcal{B}^{op}(\beta,\alpha)= \mathcal{B}(\alpha,\beta) =
\left(\mathcal{B}^{co}(\alpha,\beta) \right)^{op}$ as categories
(where $op$ of a category is basically reversing the directions of the
morphisms).

A bicategory will be called a \textit{strict} $2$\textit{-category} if
the associativity and the unit constraints are identities. A
\textit{$\mbb{C}$-linear} bicategory $\mathcal{B}$ is a bicategory
such that $\mathcal{B}(\alpha,\beta)$ is a $\C$-linear category for
every $ \alpha$, $\beta\in \mathcal{B}_{0}$ and the functor $\otimes $
is additive.
\begin{defn}
Let $\mathcal{B}$, $\mathcal{B}^{\prime }$ be bicategories. A {\it
  weak functor} $ F=(F,\varphi ):\mathcal{B\rightarrow B}^{\prime }$
consists of:
\begin{itemize}
\item a function $F:\mathcal{B}_{0}\mathcal{\rightarrow B}_{0}^{\prime
}$,
\item for all $\alpha$,$\beta\in \mathcal{B}_{0}$, there exists a
  functor $F^{\alpha,\beta}: \mathcal{B}(\alpha,\beta)\rightarrow
  \mathcal{B}^{\prime }(F(\alpha),F(\beta))$ written simply as $F$,
\item for all $\alpha$, $\beta$, $\gamma\in \mathcal{B}_{0}$, there
  exists a natural isomorphism $\varphi ^{\alpha,\beta,\gamma}:\otimes
  ^{\prime }\circ (F^{\beta,\gamma}\times F^{\alpha,\beta})\rightarrow
  F^{\alpha,\gamma}\circ \otimes $ written simply as $\varphi $ (where
  $\otimes $ and $\otimes ^{\prime }$ are the tensor functors of $
  \mathcal{B}$ and $\mathcal{B}^{\prime }$ respectively),
\item for all $\alpha\in \mathcal{B}_{0}$, there exists an invertible
  (with respect to composition) $2$-cell $\varphi
  _{\alpha}:1_{F(\alpha)}\rightarrow F(1_{\alpha})$,
\end{itemize}
satisfying commutativity of certain diagrams (consisting of $2$-cells)
which are analogous to the hexagonal and rectangular diagrams
appearing in the definition of a tensor functor.
\end{defn}
\begin{defn}
Let $F=(F,\varphi )$, $G=(G,\psi ):\mathcal{B\rightarrow B}^{\prime }$
be weak functors. A weak transformation $\sigma :F\rightarrow G$
consists of:
\begin{itemize}
\item for all $\alpha\in \mathcal{B}_{0}$, there
exists a $1$-cell $\sigma_{\alpha} \in
ob\left(\mathcal{B}^\prime (F(\alpha),G(\alpha)) \right)$,
\item for all $\alpha$, $\beta\in \mathcal{B}_{0}$, there exists a
  natural transformation $\sigma ^{\alpha,\beta}:(\sigma
  _{\beta}\otimes ^{\prime }F^{\alpha,\beta})\rightarrow
  G^{\alpha,\beta}\otimes ^{\prime }\sigma _{\alpha}$ written simply
  as $\sigma $ (where $(\sigma _{\beta}\otimes ^{\prime
  }F^{\alpha,\beta})$, $G^{\alpha,\beta}\otimes ^{\prime }\sigma
  _{\alpha}:\mathcal{B}(\alpha,\beta)\rightarrow
  \mathcal{B(}F(\alpha),G(\beta))$ are functors defined in the obvious
  way), satisfying the following:

\noindent for all $X\in ob(\mathcal{B}(\beta,\gamma))$, $Y\in
ob(\mathcal{B}(\alpha,\beta))$ where $\alpha$, $\beta$, $\gamma\in
\mathcal{B}_{0}$, the following two diagrams commute:
\[
\xymatrix @-0.8pc { \sigma _{\gamma}\otimes ^{\prime }F(X)\otimes
  ^{\prime }F(Y) \ar[rrr]^{ \sigma _{X}\otimes ^{\prime }id_{F(Y)}}
  \ar[d]_{id_{\sigma _{\gamma}}\otimes ^{\prime }\varphi _{X,Y}} & & &
  G(X)\otimes ^{\prime }\sigma _{\beta}\otimes ^{\prime }F(Y)
  \ar[rrr]^{ id_{G(X)}\otimes ^{\prime }\sigma _{Y} } & & &
  G(X)\otimes ^{\prime }G(Y)\otimes ^{\prime }\sigma _{\alpha}
  \ar[d]^{\psi _{X,Y}\otimes ^{\prime }id_{\sigma _{\alpha}}}
  \\ \sigma _{\gamma}\otimes ^{\prime }F(X\otimes Y)
  \ar[rrrrrr]_{\sigma _{X\otimes Y}} & & & & & & G(X\otimes
  Y)\otimes^{\prime} \sigma _{\alpha}}
\]
\vspace*{-4mm}
\[
\xymatrix @-0.5pc { \sigma _{\alpha}\otimes ^{\prime }id_{F(\alpha)}
  \ar[rr]^{\rho _{\sigma _{\alpha}}^{\prime }} \ar[d]_{id_{\sigma
      _{\alpha}}\otimes ^{\prime }\varphi _{\alpha} } & & \sigma
  _{\alpha} & & id_{G(\alpha)}\otimes ^{\prime } \sigma _{\alpha}
  \ar[ll]_{\lambda _{\sigma _{\alpha}}^{\prime }} \ar[d]^{ \psi
    _{\alpha}\otimes ^{\prime }id_{\sigma _{\alpha}}\ , } \\ \sigma
  _{\alpha}\otimes ^{\prime }F(1_{\alpha}) \ar[rrrr]_{\ \sigma
    _{1_{\alpha}} } & & & & G(1_{\alpha})\otimes ^{\prime }\sigma
  _{\alpha} }
\]
\noindent where $\lambda ^{\prime }$ and $\rho ^{\prime }$ are the left and right
unit constraints of $\mathcal{B}^{\prime }$ respectively.
\end{itemize}
\end{defn}

When such a weak transformation exists, we say that $F$ and $G$ are
{\em weakly isomorphic}. We have the following useful {\em Coherence
  Theorem} for bicategories. See \cite{Lei} for a proof.

\begin{thm}\label{coherence} 
Every bicategory $\mcal B$ is biequivalent to some strict $2$-category
$\mathcal{B}^{\prime }$, i.e., there exist weak functors $F:
\mathcal{B\rightarrow B}^{\prime }$ and $G:\mathcal{B}^{\prime
}\rightarrow \mathcal{B}$ such that $id_{\mathcal{B}}$
(resp., $id_{\mathcal{B}^{\prime }}$) is weakly isomorphic to $G\circ
F$ (resp., $F\circ G$).
\end{thm}

\noindent In view of this, time and again we will supress (and will
not mention about it) the associativity and unit constraints to give a
simpler look to expressions involving these constraints.

\noindent {\bf Rigid Structure on a bicategory.} Let $\alpha
\overset{X}{\ra}\beta$ be a $1$-cell in a bicategory $ \mathcal{B}$. A
\textit{right dual of} $X$ is a $1$-cell $\beta\overset{X^{\#
}}{\longrightarrow }\alpha$ such that there exist $2$-cells $X^{\#
}\otimes X\overset{e_{X}}{\ra }1_{\alpha}$ and
$1_{\beta}\overset{c_{X}}{ \ra }X\otimes X^{\# }$  satisfying
\begin{equation*}
\begin{tabular}{rcl}
$(id_{X}\otimes e_{X})\circ \left( c_{X}\otimes id_{X}\right) =id_{X}$
  & and & $(e_{X}\otimes id_{X^{\# }})\circ \left( id_{X^{\# }}\otimes
  c_{X}\right) =id_{X^{\# }}.$
\end{tabular}
\end{equation*}
 A bicategory is said to be \textit{(right) rigid} if right dual
 exists for every $1$-cell. Further, in a rigid bicategory
 $\mathcal{B}$, one can consider right dual as an invertible weak
 functor $\# =(\# ,s ):\mathcal{B}\rightarrow \mathcal{B}^{op\, co}$ in
 the following way:
\begin{itemize}
\item for each $1$-cell $X$, we fix a triplet $(X^{\# },e_{X},c_{X})$
  so that when $X=1_{\alpha}$ for a $0$-cell $\alpha$, then
  $X^{\# }=1_{\alpha}$, $ e_{X}=\lambda_{1_{\alpha}}$ ($=\rho
  _{1_{\alpha}}$, see \cite{Kas} for a proof) and $ c_{X}=\lambda
  _{1_{\alpha}}^{-1}=\rho _{1_{\alpha}}^{-1}$;
\item $\# $ induces identity map on $\mathcal{B}_{0}$;
\item for each pair of $0$-cells $\alpha$ and $\beta$,  define the 
contravariant functor $ \#
  :\mathcal{B}(\alpha,\beta)\rightarrow \mathcal{B}(\beta,\alpha)$ as
  follows: for each $X$, $Y\in ob(\mathcal{B}(\alpha,\beta))$ and
  $2$-cell $f :X\rightarrow Y$, set $\# (X)=X^{\# }$ and $\# (f )$,
  denoted by $f ^{\# }$, be given by the following composition\\
$Y^{\# }\overset{\rho _{Y^{\# }}^{-1}}{\longrightarrow }Y^{\# }\otimes
1_{\alpha}\overset{id_{Y^{\# }}\otimes c_{X}}{\longrightarrow }Y^{\#
}\otimes X\otimes X^{\# }\overset{id_{Y^{\# }}\otimes f \otimes
  id_{X^{\# }}} {\longrightarrow }Y^{\# }\otimes Y\otimes X^{\#
}\overset{e_{Y}\otimes id_{X^{\# }}}{\longrightarrow }1_{\beta}\otimes
X^{\# }\overset{\lambda _{X^{\# }}}{\longrightarrow }X^{\# }\text{;}$
\item for all $\alpha$, $\beta$, $\gamma\in \mathcal{B}_{0}$, the
  natural isomorphism $ s :\otimes \circ (flip) \circ (\#
  ^{\beta,\gamma}\times \# ^{\alpha,\beta} )\rightarrow \#
  ^{\alpha,\gamma}\circ \otimes $ is defined by:\\ for $X\in
  ob(\mathcal{B}(\alpha,\beta))$, $Y\in
  ob(\mathcal{B}(\beta,\gamma))$, the invertible $2$-cell $s _{X,Y}$
  is given by the composition \\ $ X^{\#
  }\otimes Y^{\#} \overset{id_{(X^{\# }\otimes Y^{\#})}\otimes
    c_{(Y\otimes X)}}{\longrightarrow} X^{\# }\otimes Y^{\# }\otimes
  (Y\otimes X)\otimes (Y\otimes X)^{\# } \overset{id_{X^{\# }}\otimes
    e_{Y}\otimes id_X \otimes id_{(Y\otimes X)^{\#
  }}}{\longrightarrow}(X^{\# }\otimes X)\otimes (Y\otimes X)^{\#
  }\overset{e_{X}\otimes id_{(Y\otimes X)^{\# }}}{\longrightarrow
  }(Y\otimes X)^{\#}$;
\item for all $\alpha\in \mathcal{B}_{0}$, the invertible $2$-cell $s
  _{\alpha}:1_{\alpha}\rightarrow 1_{\alpha}$ is given by identity
  morphism on $1_{\alpha}$.
\end{itemize}

 Note that the above prescription of the dual functor $(\#, s)$
 carries forward almost verbatim to another weak functor $(\tilde{\#},
 \tilde{s}): \mcal{B}^{op\, co} \ra \left( \mcal{B}^{op\,
   co}\right)^{op\, co} = \mcal{B}$. This allows us to consider the
 composition $ (\tilde{\#}, \tilde{s}) \circ (\#, s):
 \mcal{B}\ra\mcal{B}$. This is again a weak functor and we abuse
 notation to denote it by $(\#\#, t)$ and call it the {\em bi-dual
   functor}.
\begin{defn}
A bicategory $\mcal B$ is said to be pivotal if $\mcal B$ is (right)
rigid and there exists a weak transformation $a : id_{\mcal B} \ra \#
\#$ such that $a_{\vlon } = 1_{\vlon}$ for all $\vlon \in \mcal{B}_0$.
\end{defn}
We now recall some useful standard properties of a pivotal bicategory.

\begin{prop}\label{a-properties} 
Let $\mcal B$ be a pivotal bicategory with pivotality given by the
weak transformation $a : id_{\mcal B} \ra \#\#$. Then, for all
$1$-cells $\alpha \stackrel{Y}{ \ra}\beta $ and
$\beta\stackrel{X}{\ra} \gamma  $, we have
\begin{enumerate}
\item $a_{X \otimes Y} = t_{X,Y} \circ (a_X \otimes a_Y)$ and
\item $a_{X}^{\#} = a_{X^{\#}}^{-1}$. 
\end{enumerate}
\end{prop}

\noindent Let $\alpha_i\oset{X_i}{\ra}\alpha_{i+1}$, $1 \leq i \leq n$
be $1$-cells in $\mcal B$. Suppressing associativity, consider the
morphism $t_{X_1, \ldots, \, X_n} := (t_{X_1\otimes\, \cdots\,
  \otimes X_{n-1}, X_n}) \circ\, \cdots\, \circ ( t_{X_1\otimes X_2,
  X_3} \otimes id_{X_4^{\# \#} \otimes\, \cdots\, \otimes X_n^{\#
    \#}}) \circ (t_{X_1, X_2} \otimes id_{X_3^{\# \#} \otimes\,
  \cdots\, \otimes X_n^{\# \#}}) $ $\in Mor(X_1^{\# \#} \otimes\,
\cdots\, \otimes X_n^{\# \#}, (X_1 \otimes \cdots \otimes X_n)^{\# \#}
).  $ Then, a simple iterative application of Proposition
\ref{a-properties}(1) gives:
\begin{cor}\label{t-property}
$a_{X_1 \otimes\, \cdots\, \otimes X_n} = t_{X_1, \ldots,\, X_n} \circ
  (a_{X_1} \otimes\cdots \otimes a_{X_n})$ for all $1$-cells
  $\alpha_i\oset{X_i}{\ra}\alpha_{i+1}$, $1 \leq i \leq n$.
\end{cor}

\subsection{Bicategory of bifinite bimodules}\label{bimod-bicat}
In this subsection, for the sake of completeness, we first recall
certain standard facts about subfactors and modules over
$II_1$-factors (which can be found, for instance, in \cite{Bis97,
  EK98, GHJ89, Jon83, JS97, PP86, PP88,Pop94, Pop95, Sun92}). And
while doing so, we also illustrate how the collection of all bimodules
inherit a canonical structure of a pivotal bicategory. We make a
stand-in assumption that all Hilbert spaces are separable and their
inner products are linear in second and conjugate-linear in first
variable.

The following proposition gives a characterization of the basic
construction.
\begin{prop}\label{bc-fact} 
Let $A \subset B \subset C$ be unital inclusions of $II_1$-factors and
suppose there is a projection $ e \in \mscr{P}(C)$ satisfying
\begin{enumerate}
\item $exe = E_A(x) e$ for all $x \in B$ and
\item $Be = Ce$.
\end{enumerate}
Then the above tower is an instance of basic construction.
\end{prop}
\comments{\begin{pf} Recall each $II_1$-factor is algebraically
    simple, thus in view of Property $(2)$ in the hypothesis, $BeB$ is
    an ideal in $C$ and hence $BeB = C$. Also, we know that $B_1 = B
    e_1 B$, where $A \subset B \subset^{e_1} B_1$ is a basic
    construction of $A \subset B$. Thus, the map $C= BeB \ni \sum beb
    \mapsto \sum be_1 b \in Be_1B = B_1$ is an onto algebra
    homomorphism and the fact that $BeB = C$ implies (as in \cite[$\S$
      5.3]{JS97}) that it injective as well.
\end{pf}}

Given a finite index subfactor $N \subset M$, a {\em left} (resp.,
{\em right}) {\em basis of $M$ over $N$} is a finite subset $B$ of $M$
satisfying any of the following equivalent conditions:

(i) $x = \us{b \in B}{\sum} E_N (x b^\ast) b$ (resp., $x = \us{b \in
  B}{\sum} b E_N (b^\ast x)$) for all $x \in M$,

(ii) $x = \us{b \in B}{\sum} b^\ast E_N (b x)$ (resp., $x = \us{b \in
  B}{\sum} E_N (x b) b^\ast$) for all $x \in M$,

(iii) $\us{b \in B}{\sum} b^\ast e b = 1$ (resp., $\us{b \in B}{\sum}
b e b^\ast = 1$) where $e$ is a Jones projection in a basic
construction of $N \subset M$.\\

\noindent Proof of existence of such basis can be found in
\cite{PP86}.
\vspace*{2mm}

Let $A$ and $B$ be $II_1$-factors and $_A \mcal{H}$ (resp.,
$\mcal{K}_B$) be a left $A$-module (resp., right $B$-module) such that
$dim({_A}\mcal{H})  < \infty$ (resp., $dim(\mcal{K}_B) < \infty$),
equivalently, $A^\prime := {_A}\mcal{L} (\mcal{H})$ (resp., $B^\prime
:= \mcal{L}_B (\mcal{K})$) is a $II_1$-factor. Consider the set of
bounded vectors $\left( _A \mcal{H}\right)^o := \{\xi \in \mcal{H}:
\text{there exists } k>0 \text{ s.t. } \| a \xi \|^2 \leq k\; tr_A(a
a^\ast) \text{ for all } a \in A\}$ (resp., $\left( \mcal{K}_B
\right)^o := \{\eta \in \mcal{K}: \text{there exists } k>0 \text{
  s.t. } \| \eta b \|^2 \leq k \; tr_B(b b^\ast) \text{ for all } b \in
B\}$) which forms a dense subspace of $\mcal{H}$ (resp., $\mcal{K}$)
and is closed under the action of $A^\prime$ (resp.,
$B^\prime$). Using the Radon-Nikodym derivative with respect to the
faithful trace, one can obtain the $A$- (resp., $B$-) valued inner
product $ _A\langle \cdot , \cdot \rangle : (_A\mcal{H})^o \times
(_A\mcal{H})^o \ra A$ (resp., $ \langle \cdot , \cdot \rangle_B :
(\mcal{K}_B)^o \times (\mcal{K}_B)^o \ra B$) defined by the equation
$tr_A(a{_A}\langle \xi, \xi^\prime \rangle) = \langle \xi, a
\xi^\prime \rangle$ (resp., $tr_B( \langle \eta, \eta^\prime \rangle_B
b) = \langle \eta, \eta^\prime b \rangle$) for all $\xi, \xi^\prime
\in ({_A}\mcal{H})^o, a \in A$ (resp., $\eta, \eta^\prime \in
(\mcal{K}_B)^o, b \in B$). It is easy to check that the inner product
has the following properties:
\begin{enumerate}
\item $_A\langle \xi, \xi \rangle$ (resp., $\langle \eta, \eta
  \rangle_B \geq 0$),
\item $ _A\langle \xi, \xi^\prime \rangle^* =\, _A\langle \xi^\prime,
  \xi \rangle$ (resp., $\langle \eta, \eta^\prime \rangle_B ^* =
  \langle \eta^\prime, \eta \rangle_B$),
\item $ _A\langle \xi, a \xi^\prime \rangle = a\, _A\langle \xi,
  \xi^\prime \rangle$ and $ _A\langle a\xi, \xi^\prime \rangle = {_A}
  \langle \xi, \xi^\prime \rangle a^\ast$ (resp., $\langle \eta,
  \eta^\prime b \rangle_B = \langle \eta, \eta^\prime \rangle_B \, b$
  and $\langle \eta b, \eta^\prime \rangle_B = b^\ast \langle \eta,
  \eta^\prime \rangle_B$),
\item $ _A\langle \xi, x \xi^\prime \rangle = {_A}\langle x^* \xi,
  \xi^\prime \rangle$ (resp., $\langle \eta, y \eta^\prime \rangle =
  \langle y^* \eta, \eta^\prime \rangle_B$) \comments{
\item $tr_A(_A\langle \xi, \eta\rangle ) = \langle \xi, \eta\rangle$,
  $tr_B(\langle \xi, \eta\rangle_B ) = \langle \xi, \eta\rangle$}
\end{enumerate}
for all $a \in A$, $x \in A^\prime$, $\xi,\xi^\prime \in
({_A}\mcal{H})^o$ (resp., $b \in B$, $y \in B^\prime$,
$\eta,\eta^\prime \in (\mcal{K}_B)^o$).  Also, given a finite index
subfactor $C$ of $A$ (resp., $B$), one can use a basis for $C \subset
A$ (resp., $C \subset B$) to obtain $({_A}\mcal{H})^o =
({_C}\mcal{H})^o$ (resp., $(\mcal{K}_B)^o = (\mcal{K}_C)^o$) where the
$C$-valued inner product is given by the $A$- (resp., $B$-) valued
inner product composed with the trace preserving conditional
expectation onto $C$.  Further, there exists a finite subset
$\{\xi_i\}_i$ (resp., $\{\eta_j\}_j$) of $(_A\mcal{H})^o$ (resp.,
$(\mcal{K}_B)^o$) satisfying $id_{({_A}\mcal{H})^o} = \us{i}{\sum}
{_A}\langle \xi_i, \cdot \rangle \xi_i$ (resp., $id_{(\mcal{K}_B)^o} =
\us{j}{\sum} \eta_j \langle \eta_j, \cdot \rangle_B$); such a subset
is called {\em basis} for the module. Such a basis also satisfies the
following conditions which are completely straight forward to verify.
\begin{prop}\label{trace-A'}
(i) $tr_{A'} (x) = [dim ({_A} \mcal{H})]^{-1} \us{i}{\sum} \langle
  \xi_i, x \xi_i \rangle$ (resp., $tr_{B'} (y) = [dim(\mcal{K}_B)]^{-1}
  \us{j}{\sum} \langle \eta_j, \eta_j y \rangle$) for all $x \in A'$
  (resp., $y \in B'$) which implies $({_A}\mcal{H})^o =
  ({_{A'}}\mcal{H})^o$ (resp., $(\mcal{K}_B)^o =
  ({_{B'}}\mcal{K})^o$),

(ii) $\us{i}{\sum} {_{A'}}\langle \xi_i, \xi_i \rangle =
  dim({_A}\mcal{H}) 1_{A'}$ (resp., $\us{j}{\sum} {_{B'}}\langle
  \eta_j, \eta_j \rangle = dim( \mcal{K}_B) 1_{B'}$).
\end{prop}
The dual Hilbert space or the contragredient $\oline{\mcal{H}} := \{
\oline{\xi}:\xi \in \mcal{H}\}$ (resp., $\oline{\mcal{K}} := \{
\oline{\eta}:\eta \in \mcal{K}\}$) where bar being a conjugate linear
unitary, can be equipped with a right $A$- (resp., left $B$-) module
structure given by $\oline{\xi} a = \oline{a^\ast \xi}$ for $a\in A,
\xi \in \mcal{H}$ (resp., $b \oline{\eta} = \oline{\eta b^\ast }$ for
$b \in B, \eta \in \mcal{K}$). Note that (i) $(\oline{\mcal{H}}_A)^o =
\oline{({_A}\mcal{H})^o}$ (resp., $({_B}\oline{\mcal{K}})^o =
\oline{(\mcal{K}_B)^o}$), (ii) $\langle \oline{\xi} ,
\oline{{\xi}^\prime} \rangle_A = {_A}\langle \xi^\prime , \xi \rangle$
for $\xi, \xi^\prime \in ({_A}\mcal{H})^o$ (resp., ${_B}\langle
\oline{\eta} , \oline{{\eta}^\prime} \rangle = \langle \eta^\prime ,
\eta \rangle_B$ for $\eta, \eta^\prime \in (\mcal{K}_B)^o$),
(iii) $\{\oline{\xi_i}\}_i$ (resp., $\{\oline{\eta_j}\}_j$) is a basis
for $\oline{\mcal{H}}_A$ (resp., $_B\oline{\mcal{K}}$) and (iv)
$dim({_A}\mcal{H}) = dim(\oline{\mcal{H}}_A)$ (resp., $dim(\mcal{K}_B)
= dim({_B}\oline{\mcal{K}})$).


Next, we briefly recall few aspects of bimodules over $II_1$-factors
$A$ and $B$. Let $\mcal{H}$ be an $A$-$B$-bimodule. If $dim
(\mcal{H}_B) < \infty$, then the Jones index of the subfactor $A \subset
B'$ turns out to be $[B': A] = dim ({_A}\mcal{H})\, dim (\mcal{H}_B) =:
index ( {_A}{\mathcal H}_B )$. A bimodule $_A{\mathcal H}_B$ is called
{\em bifinite} if $index ( {_A}{\mathcal H}_B ) < \infty$, and a
bifinite bimodule $_A{\mathcal H}_B$ is called {\em extremal} if the
canonical traces of the $II_1$-factors $A'$ and $B'$ coincide on the
intertwiner space $_A{\mathcal L}_B ({\mathcal H})$ (which has finite
complex dimension due to the finiteness of the index). Note that if
$index (_A\mcal{H}_B) = 1$, then $_A\mcal{H}_B$ is an irreducible
$A$-$B$-bimodule. Also, if $index(_A\mcal{H}_B) < \infty$, then
$(_A\mcal{H})^o = (\mcal{H}_B)^o$; we will write $\mcal{H}^o$ for this
space of bounded vectors.
      
Given two bifinite bimodules $_A\mcal{H}_B$ and $_B\mcal{K}_C$, one
can consider the tensor product $\mcal{H} \us{B}{\otimes} \mcal{K}$
defined as the completion of the space $\mcal{H}^o \us{B}{\otimes}
\mcal{K}^o := \frac{\displaystyle \mcal{H}^o \us{alg}{\otimes}
  \mcal{K}^o }{\displaystyle span \left\{\xi a \otimes \eta - \xi
  \otimes a \eta \left| a \in A, \xi \in \mcal{H}^o, \eta \in
  \mcal{K}^o \right. \right\} }$ with respect to the inner product
$\langle \xi \us{B}{\otimes} \eta,\xi^\prime \us{B}{\otimes}
\eta^\prime\rangle = \langle \xi ,\xi^\prime {_B}\langle \eta ,
\eta^\prime \rangle \rangle = \langle \eta, \langle \xi , \xi^\prime
\rangle_B \eta^\prime\rangle$ for $\xi, \xi^\prime \in \mcal{H}^o,
\eta, \eta^\prime \in \mcal{K}^o$, which is equipped with the obvious
$A$-$C$ bimodule structure. The following is a list of very useful
properties of this tensor product.
\begin{enumerate}
\item $(\mcal{H} \us{B}{\otimes} \mcal{K})^o = \mcal{H}^o
  \us{B}{\otimes} \mcal{K}^o$.
\item The $A$- (resp., $C$-) valued inner product is given by ${_A}
  \langle \xi \us{B}{\otimes} \eta,\xi^\prime \us{B}{\otimes}
  \eta^\prime\rangle = {_A}\langle \xi, \xi^\prime {_B}\langle \eta,
  \eta^\prime \rangle \rangle$ (resp., $\langle \xi \us{B}{\otimes}
  \eta,\xi^\prime \us{B}{\otimes} \eta^\prime\rangle_C = \langle \eta,
  \langle \xi , \xi^\prime \rangle_B \eta^\prime\rangle_C$) for $\xi,
  \xi^\prime \in \mcal{H}^o, \eta, \eta^\prime \in \mcal{K}^o$.
\item If $\{ \xi_i \}_i$ and $\{ \eta_j \}_j$ are basis for
  ${_A}\mcal{H}$ and $_B\mcal{K}$ (resp., $\mcal{H}_B$ and
  $\mcal{K}_C$) respectively, then $\{\xi_i \us{B}{\otimes}
  \eta_j\}_{i,j}$ forms a basis for $_A\mcal{H} \us{B}{\otimes}
  \mcal{K}$ (resp., $\mcal{H} \us{B}{\otimes} \mcal{K}_C$).
\item The left dimension, the right dimension and the index of the
  bifinite bimodules are multiplicative with respect to this tensor
  product.
\item The map $\mcal{L}_B(\mcal{H}) \ni x \mapsto x \us{B}{\otimes}
  id_{\mcal{K}} \in \mcal{L}_C (\mcal{H} \us{B}{\otimes} \mcal{K})$
  (resp., ${_B}\mcal{L}(\mcal{K}) \ni y \mapsto id_{\mcal{H}}
  \us{B}{\otimes} y \in {_A}\mcal{L} (\mcal{H} \us{B}{\otimes}
  \mcal{K})$) is an inclusion of unital $*$-algebras.
\end{enumerate}
Let $\mcal{H}$ be an $A$-$B$-bimodule with $dim({_A}\mcal{H}) <
\infty$ (resp., $dim(\mcal{H}_B) < \infty$) and $\{ \xi_i\}_i$ (resp.,
$\{ \eta_j\}_j$) be a basis for $_A\mcal{H}$ (resp.,
$\mcal{H}_B$). Then, it is easy to see that the bounded vector
$\us{j}{\sum} {\eta}_j \uset{B}{\otimes} \oline{\eta_j}$ (resp.,
$\us{i}{\sum} \oline{\xi_i} \uset{A}{\otimes} \xi_i$) is independent
of the basis and is $A$-$A$-central, that is, $ a (\us{j}{\sum}
{\eta}_j \uset{B}{\otimes} \oline{\eta_j}) = (\us{j}{\sum} {\eta}_j
\uset{B}{\otimes} \oline{\eta_j}) a $ for all $a \in A$ (resp.,
$B$-$B$-central).
\begin{lem}\label{basic-con-intertwiners}
Let $\mcal{K}$ be a right $B$-module with $dim(\mcal{K}_B) < \infty$
and $\mcal{H}$ be a bifinite $A$-$B$ bimodule for $II_1$ factors $A$ and
$B$. Then, the inclusion of $II_1$ factors
\[
\begin{array}{c}
{\mcal{L}_B}(\mcal{K}) \hookrightarrow {\mcal{L}_A}(\mcal{K}
\uset{B}{\otimes} \oline{\mcal{H}}) \hookrightarrow
     {\mcal{L}_B}(\mcal{K} \uset{B}{\otimes} \oline{\mcal{H}}
     \uset{A}{\otimes} \mcal{H} )\\
\hspace*{-5mm} x \mapsto x \uset{B}{\otimes} id_{\oline{\mcal H}}, \ \,  y
\mapsto y \uset{A}{\otimes} id_{\mcal H}
\end{array}
\]
is an instance of basic construction with Jones projection $e$ given
(on bounded vectors) by
\[
\mcal{K} \uset{B}{\otimes} \oline{\mcal{H}} \uset{A}{\otimes} \mcal{H}
\ni \xi \uset{B}{\otimes} \bar{\eta} \uset{A}{\otimes} \zeta
\oset{e}{\longmapsto} \frac{1}{dim\, _A\mcal{H}} \sum_i (\xi
\uset{B}{\otimes} \bar{\xi}_i \uset{A}{\otimes} \xi_i) \langle \eta,
\zeta \rangle_B \, \in \mcal{K} \uset{B}{\otimes} \oline{\mcal{H}}
\uset{A}{\otimes} \mcal{H},
\]
where $\{\xi_i \}_i$ is a basis for $_A \mcal{H}$.
\end{lem}
\begin{pf}
Set $N = {\mcal{L}_B}(\mcal{K})$, $M = {\mcal{L}_A}(\mcal{K}
\uset{B}{\otimes}\oline{\mcal{H}} )$ and $M_1 = {\mcal{L}_B}(\mcal{K}
\uset{B}{\otimes} \oline{\mcal{H}} \uset{A}{\otimes} \mcal{H} )$. Let
$\{\eta_j\}_j$ (resp., $\{\sigma_k\}_k$) be a basis for $\mcal{H}_B$
(resp., $\mcal{K}_B$). Then, using Proposition \ref{trace-A'}, it is
completely routine to check that the map $E_N: M \ra N$ given by
$E_N(x) (\xi) = [dim ({_A}\mcal{H})]^{-1} \us{i, i', k}{\sum} \sigma_k
\langle \xi_{i'} , \langle \sigma_k \uset{B}{\otimes} \bar{\xi}_{i'} ,
x(\xi \uset{B}{\otimes} \bar{\xi}_i ) \rangle_A \xi_i \rangle_B$ for
all $\xi \in \mcal{K}^o, x \in M$, is the unique $tr_M$ preserving
conditional expectation from $M$ onto $N$. Further, it can also be
readily shown that $e \in \mscr{P} (N' \cap M_1)$ and $e (x
\uset{A}{\otimes} id_{\mcal H})e = e( E_N(x)\uset{B}{\otimes}
id_{\oline{\mcal{H}}\uset{A}{\otimes} \mcal{H} }) $ for all $x \in
M$. In view of Proposition \ref{bc-fact}, it just remains to show that
$M_1e = Me$. Let $x_1 \in M_1$ and consider the map $x : \mcal{K}
\uset{B}{\otimes}\oline{\mcal{H}} \ra \mcal{K}
\uset{B}{\otimes}\oline{\mcal{H}}$ given (on bounded vectors) by $x (
\xi \uset{B}{\otimes} \bar{\eta} ) = \sigma_k \uset{B}{\otimes}
\bar{\xi}_{i'} \, _A\langle \eta, \eta_j \langle \sigma_k
\uset{B}{\otimes} \bar{\xi}_{i'} \uset{A}{\otimes} \eta_j, x_1( \xi
\uset{B}{\otimes} \bar{\xi}_i \uset{A}{\otimes} \xi_i )\rangle_B $ for
all $ \xi \in (\mcal{K}_B)^o, \eta \in (\mcal{H}_B)^o$. Clearly, $x
\in M$, and it involves nothing more than straightforward verification
to show that $(x \uset{A}{\otimes} id_{\mcal H})e = x_1 e$.
\end{pf}

We denote the bicategory of bifinite bimodules with ${\mathcal B}$
whose $0$-cells are $II_1$-factors; for $II_1$-factors $A$ and $B$,
the objects of the category ${\mathcal B} (B,A)$ are bifinite $A$-$B$
bimodules and morphisms or $2$-cells are $A$-$B$ linear maps between
such bimodules (which are automatically bounded). The tensor functor
is given by the usual relative tensor product of bimodules and for
each $II_1$-factor $A$ the identity object in $\mcal{B}(A, A)$ is the
canonical $A$-$A$-bimodule $L^2(A)$. There is a natural associativity
constraint for relative tensor product of bimodules. Further, for an
$A$-$B$-bimodule $\mcal H$, the unit constraints are given by the
canonical isomorphisms $ L^2(A) \uset{A}{\otimes} \mcal{H}
\oset{A\text{-}B}{\cong} \mcal{H}$ and $ \mcal{H} \uset{B}{\otimes}
L^2(B) \oset{A\text{-}B}{\cong} \mcal{H}$. Thus, $\mcal B$ has a
natural bicategory structure. For the (right) rigid structure on
$\mcal B$, for each $A$-$B$-bimodule $\mcal H$, we set
$(_A\mcal{H}_B)^\# =\, _B\oline{\mcal H}_A$ and define the evaluation
and coevaluation maps $e_{\mcal H}\in\, _B \mcal{L}_B
(\oline{\mcal{H}} \uset{A}{\otimes} \mcal{H} , L^2(B))$ and $c_{\mcal
  H} \in \, _A \mcal{L}_A (L^2(A), \mcal{H}\uset{B} {\otimes}
\oline{\mcal{H}} ) $ respectively, (on bounded vectors) by
\vspace*{-2mm}
\[
e_{\mcal H}(\bar{\xi} \uset{A}{\otimes}\eta ) = \langle\xi , \eta
\rangle_B\text{ and } c_{\mcal H}(\hat{a}) = \sum_i a (\eta_j
\uset{B}{\otimes} \bar{\eta}_j)\ \text{ for all $\xi, \eta \in
  \mcal{H}^o$, $a \in A$, $b\in B$},
\]
\vspace*{-2mm}
where $\{\eta_j \}$ is a basis for the right $B$-module
$\mcal{H}_B$. Thus, $\mcal{B}$ indeed inherits a canonical rigid
structure. Finally, the canonical isomorphism $_A\mcal{H}_B
\oset{A\text{-}B}{\cong} \, _A\oline{\oline{\mcal H}}_B$ for any
bifinite $A$-$B$-bimodule $_A\mcal{H}_B$, equips $\mcal B$ with a
pivotal structure. Note that, for $\theta \in\, _A\mcal{L}_B (\mcal H,
\mcal K)$, it can be easily shown that $\theta^{\#}(\bar{\xi}) =
\oline{\theta^*(\xi)}$ and, hence, $
\theta^{\#\#}\left(\bar{\bar{\xi}}\right) =
\oline{\oline{\theta(\xi)}}$ for all $\xi \in \mcal{H}^o$, where
$\theta^*$ is the usual adjoint of the intertwiner $\theta$.

\section{Perturbations of planar algebras}\label{pa}

In this section, we define perturbation of a planar algebra to
obtain a new one which has the same filtered algebra structure but the
action of Jones projections and conditional expectation tangles
differ. As we will see, this will turn out to be crucial in the
following sections.  We must mention here that such objects already
appeared in the work of Michael Burns \cite{Bur03} while extending
Jones theorem (of associating a spherical planar algebra to an
extremal subfactor) in the non-extremal case; however, for us, these
arose in a purely different context, namely, while detecting the
effect of different pivotal structures on the planar algebra
associated to a $1$-cell in a strict $2$-category (as in
\cite{Gho08}).
\begin{defn}\label{weight}
Let $P$ be a planar algebra. An invertible element $z \in P_{+1} $ is
said to be a weight of $P$ if $z_{\vlon k} \in \mcal{Z}(P_{\vlon k})$
for all $\vlon k \in Col$, where
\vspace*{-2mm}
\[
\psfrag{z}{$z$}
\psfrag{zinv}{$z^{-1}$}
\psfrag{kstrings}{$\cdots k$\text{ strings}}
\psfrag{+}{$+$}
\psfrag{-}{$-$}
z_{+k} := P_{
\includegraphics[scale=0.3]{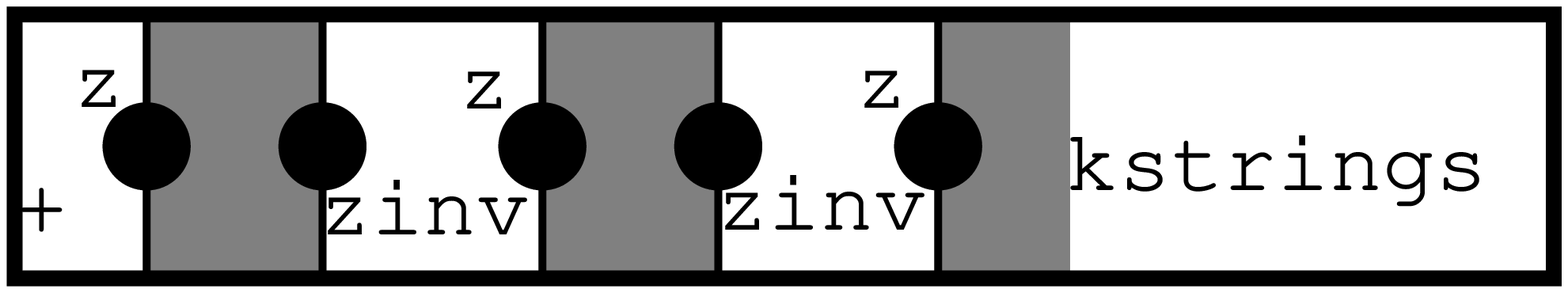}
}\text{ and } z_{-k} := P_{
\includegraphics[scale=0.3]{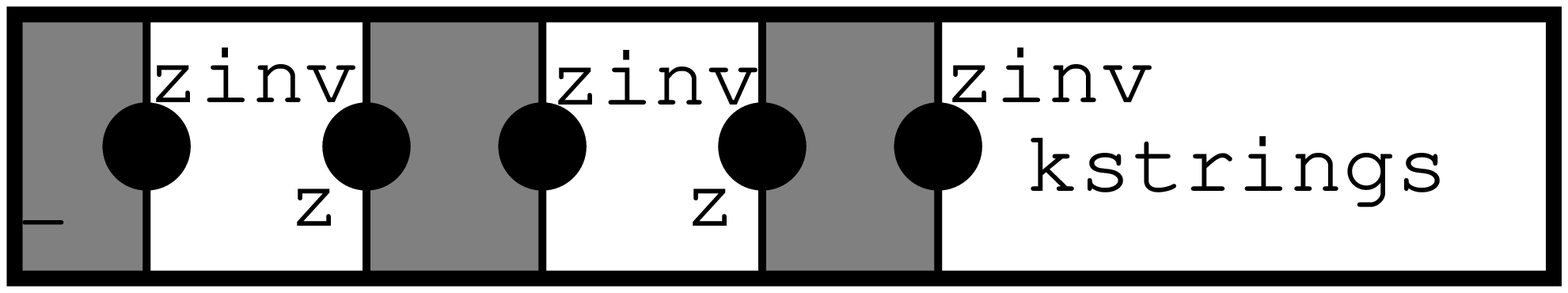}
}.
\]
\vspace*{-4mm}
\end{defn}
Given a weight $z$ of a planar algebra $P$ and an invertible
decomposition $z = ab$ for $a,b$ invertible in $P_{+1}$, we now
construct a new planar algebra $P^{(a,b)}$ as follows:

\noindent (i) {\em Vector spaces}: $P^{(a,b)}_{\vlon k} : = P_{\vlon
  k}$ for all $ \vlon k \in Col$.

\begin{figure}[h]
\psfrag{a}{$a$}
\psfrag{b}{$b$}
\psfrag{ainv}{$a^{-1}$}
\psfrag{binv}{$b^{-1}$}
\psfrag{by}{by}
\psfrag{;}{;}
\includegraphics[scale=0.15]{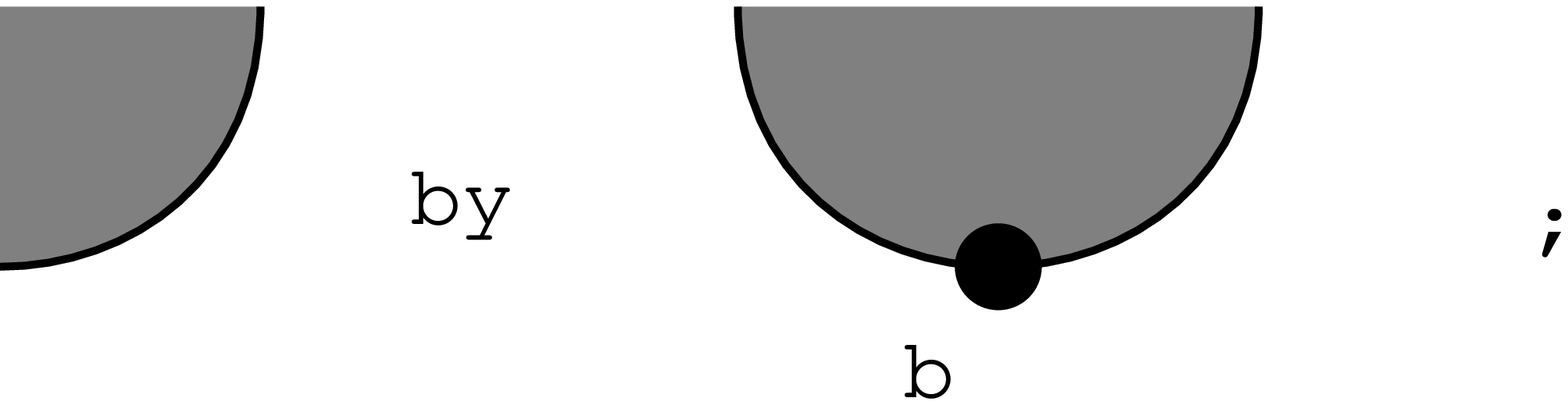}
\vspace*{-4mm}\caption{Perturbing a planar algebra}\label{max-min}
\end{figure}

\noindent (ii){\em Actions of tangles}: Let $T$ be a tangle and
$\hat{T}$ be a standard form representative (see \cite[$\S$4]{Gho08})
of the isotopy class of $T$. We replace each local maximum and minimum
appearing in $\hat T$ as in Figure \ref{max-min} and call the
resulting semi-labelled tangular diagram $\hat{T}^{(a,b)}$.
\vspace*{-4mm}
Define $P^{(a,b)}_T := P_{\hat{T}^{(a,b)}}$. The immediate
thing to check is the well-definedness of $P^{(a,b)}_T$. Note
that the above prescription is invariant under the sliding, wiggling
and $360^{\circ}$-rotation moves, applying a finite sequence of which
takes one standard form representative to another. Hence,
$P^{(a,b)}_T$ is well defined and $P^{(a,b)}$ is a planar algebra. As
mentioned before, note that $P^{(a,b)}$ has same filtered algebra
structure as $P$ whereas the action of Jones projection tangles and conditional
expectation tangles differ.

We will refer $P^{(a,b)}$ as the {\em perturbation of $P$
by the decomposition $z = ab$ of the weight $z$}.
\begin{rem}
For an invertible decomposition $z = a b$ of a weight $z$ of a planar
algebra $P$ and any $\lambda \in \C \setminus \{ 0\}$, $P^{(a,b)} =
P^{(\lambda a, \lambda^{-1}b)}$. Further, the planar algebras $P^{(a,
b)}, P^{(b, a)}, P^{(z, 1)},$ and $ P^{(1, z)}$ are all
isomorphic. Hence, up to isomorphism, the perturbation of $P$ only
depends upon the weight $z$.
\end{rem}
\noindent To see the first part, observe that the number of local 
maxima is the same as the number of local minima in a standard form
representative of a tangle, which results in cancellation of the
scalars appearing in tangle maps due to $\lambda$ in the latter
perturbation.  In the second part, for instance, the isomorphism
$P^{(a, b)} \cong P^{(z, 1)}$ is obtained by the maps
\[
P^{(a, b)}_{+ k }= P_{+k} \ni x \longmapsto P_{
\psfrag{x}{$x$}
\psfrag{+}{$+$}
\psfrag{b}{$b$}
\psfrag{binv}{$b^{-1}$}
\psfrag{kstrings}{$\cdots k \text{ strings}$} 
\includegraphics[scale=0.2]{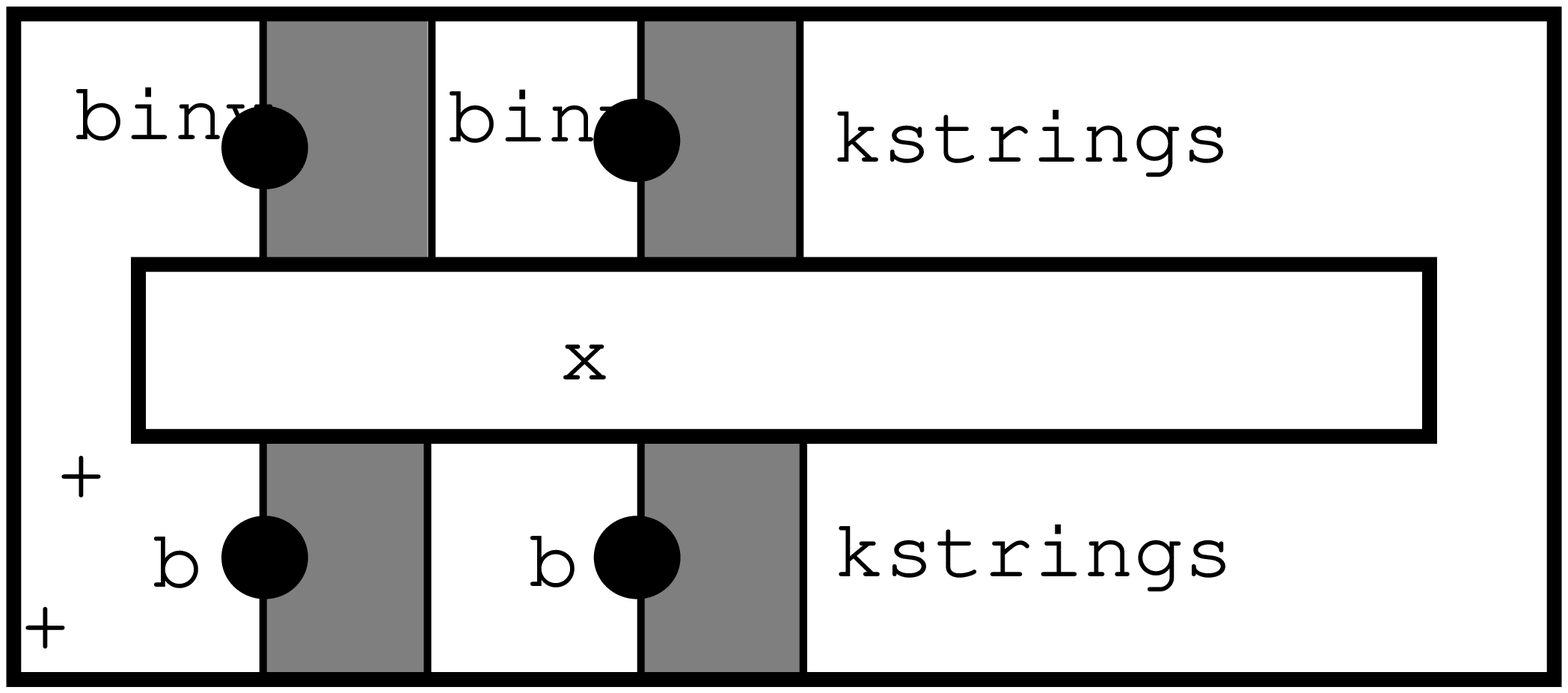}
}
\; \in P_{+ k} = P^{(z, 1)}_{+ k},
\]
\[
P^{(a, b)}_{- k } = P_{-k} \ni x \longmapsto P_{
\psfrag{x}{$x$}
\psfrag{+}{$+$}
\psfrag{-}{$-$}
\psfrag{b}{$b$}
\psfrag{binv}{$b^{-1}$}
\psfrag{kstrings}{$\cdots k \text{ strings}$} 
\includegraphics[scale=0.2]{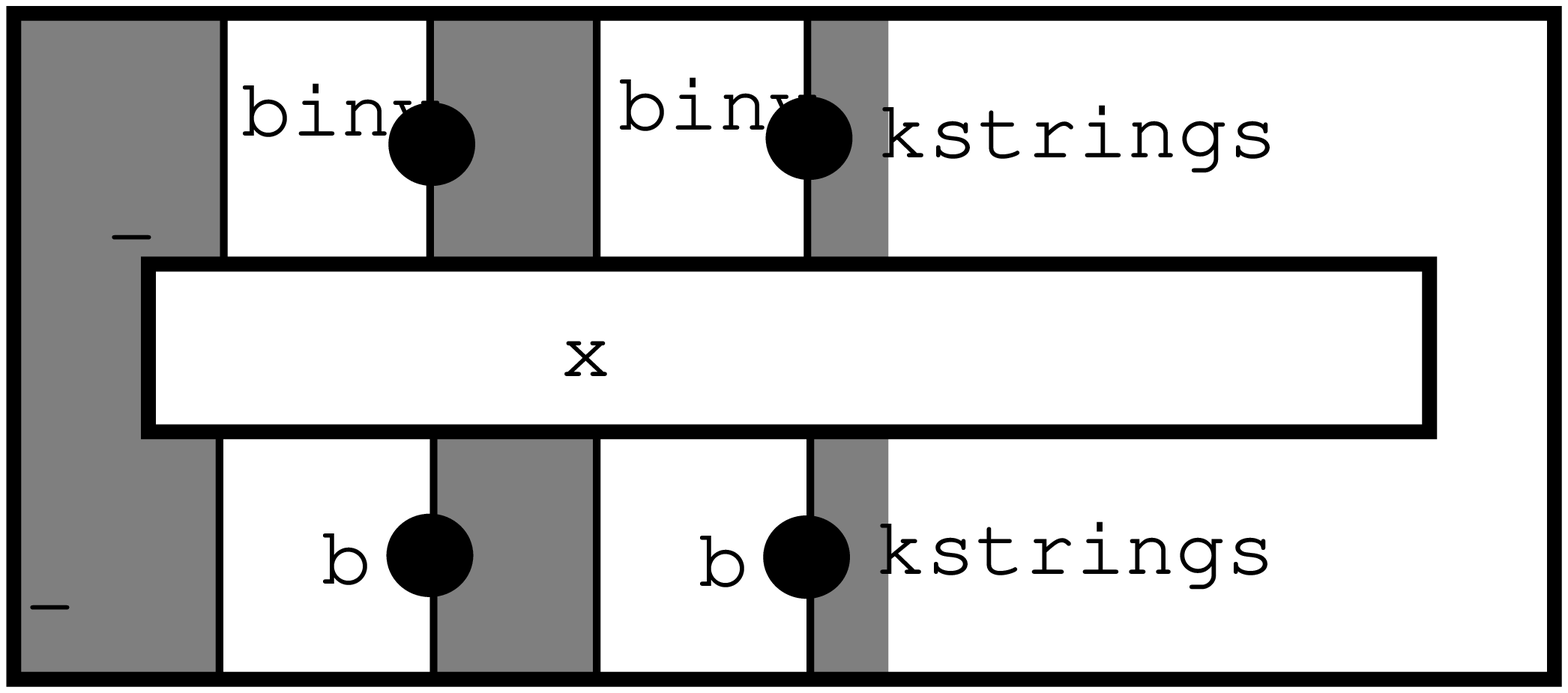}
}
\; \in P_{- k} = P^{(z, 1)}_{- k}.
\]
It is straight forward to verify the equivariance of this isomorphism
with the actions of the generating tangles in Figure \ref{tangles} and
hence all tangles.

A trivial example of a weight of a planar algebra $P$ is a non-zero
scalar $\lambda \in \C $. By the above remark, $P^{(\lambda, \mu)} =
P^{(\lambda \mu, 1)} = P^{(1, \lambda \mu)} = P^{(\mu, \lambda)}$ for
all non-zero scalars $\lambda$ and $\mu$. We will usually refer to
such perturbations as {\em scalar perturbations}.

A perturbation of a planar algebra with modulus need not have modulus (except
for perturbations by scalar weights). However, if the planar algebra
is connected then so are its perturbations, but the moduli of the
perturbations might vary. For instance, for a planar algebra $P$ with
modulus $(\delta_-, \delta_+)$, the scalar perturbation $P^{(\lambda,
1)}$ has modulus $( \lambda^{-1}
\delta_-,\, \lambda \delta_+ )$.

\begin{defn}
The normalized planar algebra associated to a planar algebra $P$ with
modulus $(\delta_-, \delta_+)$ is its scalar perturbation by the weight
$ \sqrt{\frac{\delta_-}{\delta_+}}$.
\end{defn}
\noindent Note that the normalization of  $P$ is a unimodular planar
algebra.  Although scalar perturbations change the modulus, the index
however remains the same; in the last section, we will come across an
example of a perturbation class whose normalized planar algebras
realize all indices greater than or equal to $4$.

\begin{defn}
A connected planar algebra $P$ is said to be spherical if the actions
of $0$-tangles in its normalization are invariant under spherical
isotopy.
\end{defn}
Note that this property is equivalent to demanding that the normalized
left and right picture traces on $P_{+1}$ are identical. The relevance
of the above definition will become clear in the section on bimodule
planar algebras, where we establish a correspondence between
sphericality and extremality. We must also point out that a
non-unimodular planar algebra could be spherical according to the
above definition which is not allowed in Jones' original definition
in \cite{Jon}.

In general, a perturbation of a $\ast$-planar algebra need not be a
$\ast$-planar algebra. However, for certain specific weights, the
perturbations also turn out to be $\ast$-planar algebras. For
instance, if $P$ is a $\ast$- (resp., positive) planar algebra, it is
routine to verify that a perturbation of the type $P^{(a, \lambda
  a^*)}$ for non-zero real (resp., positive) scalar $\lambda$ becomes
a $\ast$- (resp., positive) planar algebra with $\ast$-structure
coming from the original one.

\section{Weights and Pivotality}\label{weights}

In this section, we mention how one can canonically associate a
pivotal $\C$-linear strict $2$-category to a planar algebra;
conversely, from \cite{Gho08}, we recall how to associate a planar
algebra to a $1$-cell in a pivotal $\C$-linear strict $2$-category. We
then establish a relation between weights and perturbations of planar
algebras and pivotal structures on bicategories.

\subsection{Planar algebras to bicategories}
 Let $P$ be a planar algebra. From $P$, we first describe a
 $\C$-linear strict $2$-category $\mcal{B}$ and see that it inherits
 canonical rigid and pivotal structures. Set $\mcal{B}_0 = \{+,-\}$;
 for $\vlon, \eta \in \mcal{B}_0$, set $ob(\mcal{B}(\vlon, \eta)) =
 2\N_0 + \delta_{\vlon \neq \eta}$.
 To avoid confusion, we will write an object $k \in
 ob(\mcal{B}(\vlon, \eta))$ as $\vlon k $ (whence $\eta =
 (-)^k \vlon$). For two objects $\vlon k, \vlon l$, set $Mor(\vlon
 k, \vlon l) = P_{(-)^k \vlon \left( \frac{ k+l }{2}\right)} =
 P_{(-)^l \vlon \left( \frac{ k+l }{2}\right)}$. Composition of
 morphisms be given by the bilinear map
\[
Mor(\vlon l, \vlon m) \times Mor(\vlon k, \vlon l) \ni (x, y) \mapsto
x \circ y := P_{
\psfrag{x}{$x$}
\psfrag{y}{$y$}
\psfrag{k}{$k$}
\psfrag{l}{$l$}
\psfrag{m}{$m$}
\psfrag{e}{$(-)^k \vlon$}
\includegraphics[scale=0.2]{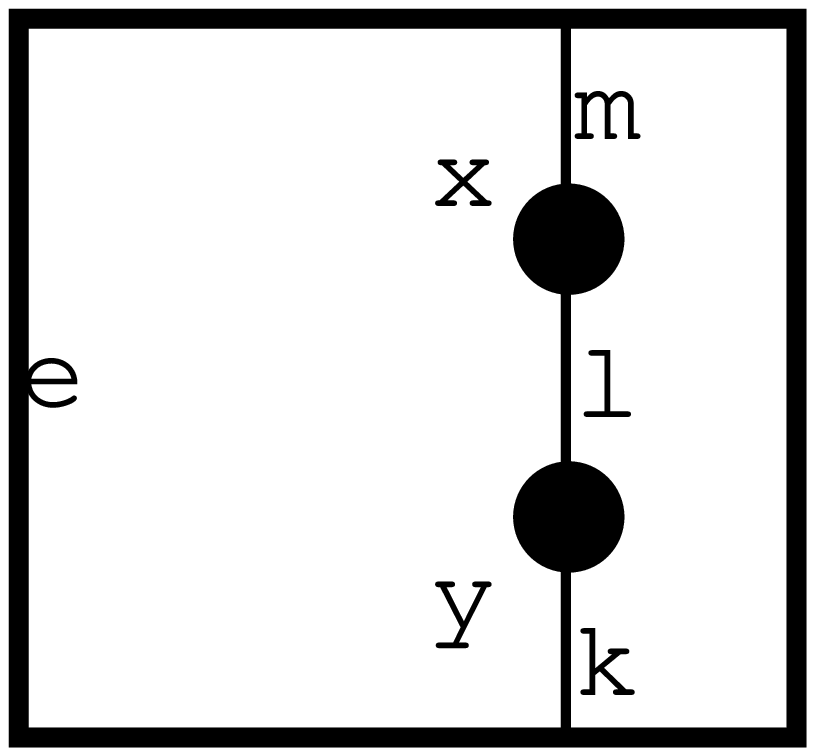}
} \; \in Mor(\vlon k, \vlon m).
\] 
For each $1$-cell $\vlon k$, the identity morphism is
 given by $P_{1_{\vlon k}}$. The tensor functor $\otimes
\, : \, \mcal{B}(\eta, \sigma) \times \mcal{B}(\vlon, \eta) \ra
 \mcal{B}(\vlon, \sigma)$ is defined by $\eta l \otimes \vlon k := \vlon
 (k+l)$ and 
\[ 
Mor(\eta m, \eta n) \times Mor(\vlon k, \vlon l) \ni (x, y) \mapsto
x \otimes y := P_{
\psfrag{x}{$x$}
\psfrag{y}{$y$}
\psfrag{k}{$k$}
\psfrag{l}{$l$}
\psfrag{m}{$m$}
\psfrag{n}{$n$}
\psfrag{s}{$\sigma$} 
\includegraphics[scale=0.2]{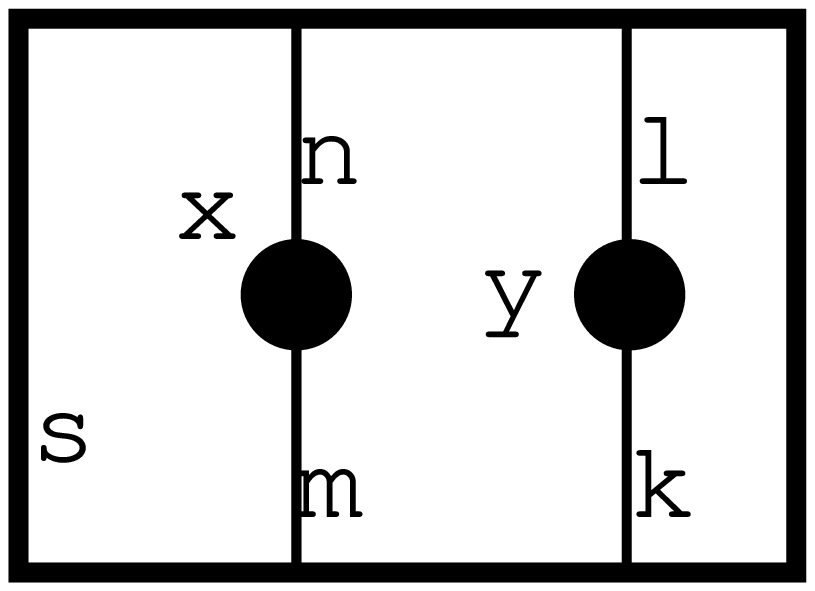}
} \; \in Mor(\vlon (k+m) , \vlon (l + n)).
\]
For each $0$-cell $\vlon$, set the identity object $1_{\vlon } = \vlon
0 \in ob(\mcal{B}(\vlon, \vlon))$. With the above structure, it is
easily seen that $\mcal{B}$ is a strict $2$-category.  We now describe
a (right) rigid structure $\#$ on $\mcal{B} $ as follows:

{\em $0$-cells:} $\vlon^{\#} := \vlon$.

{\em $1$-cells:} $(\vlon k)^{\#} = (-)^k \vlon k$.

The evaluation map $e_{\vlon k}: (\vlon k)^{\#} \otimes \vlon k = \vlon
2k \ra \vlon 0$ (resp., coevaluation map $ c_{\vlon k} : (-)^k \vlon
0 \ra \vlon k \otimes (\vlon k)^{\#} = (-)^k \vlon 2k$) will be given
by $ e_{\vlon k}:= P_{
\psfrag{e}{$\vlon$}
\psfrag{k}{$k$}
\includegraphics[scale=0.15]{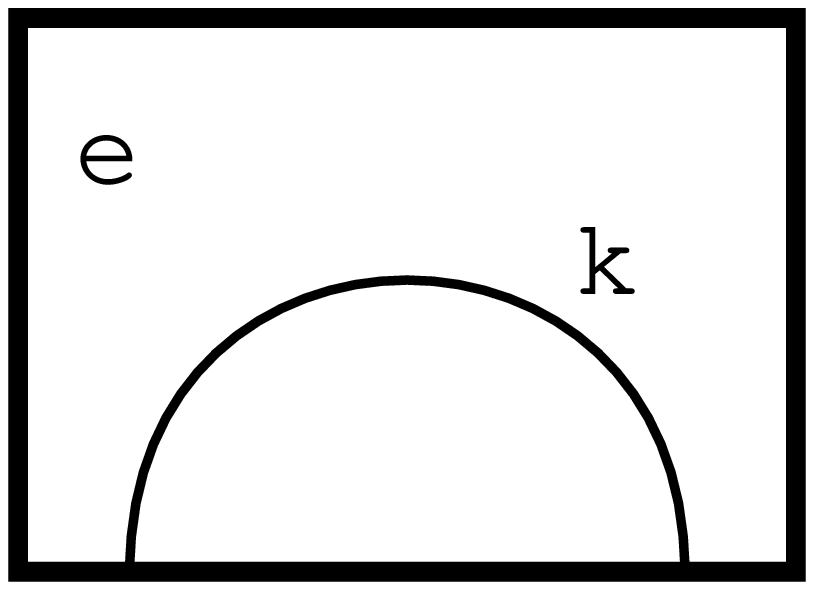}
}$
(resp., $c_{\vlon k}:=\ P_{
\psfrag{e}{$(-)^k \vlon$}
\psfrag{k}{$k$}
\includegraphics[scale=0.15]{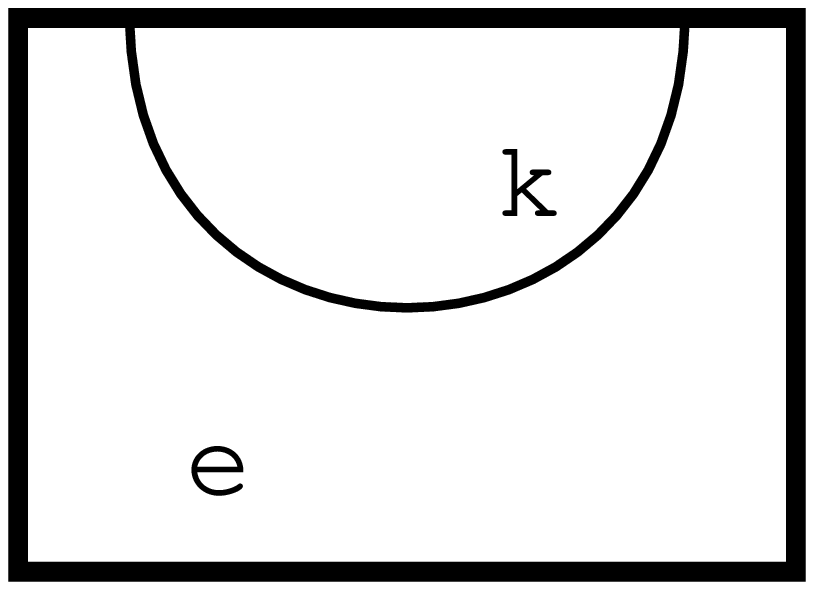}
}$).
Thus, we obtain a weak functor $\# : \mcal{B} \ra \mcal{B}^{op\, co}$
which yields $x^{\#} = P_{
\psfrag{x}{$x$}
\psfrag{k}{$k$}
\psfrag{l}{$l$}
\psfrag{e}{$\vlon$}
\includegraphics[scale=0.15]{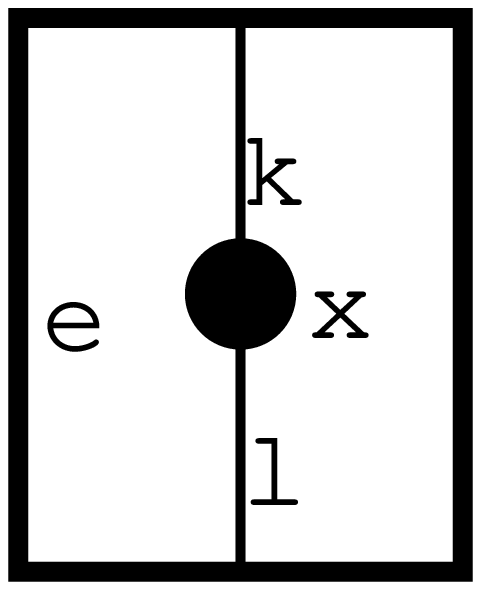}
}$ for all $x \in Mor(\vlon k, \vlon l)$, and
satisfies $\# \circ \# = id_{\mcal{B}}.$

This strict $2$-category $\mcal B$ also inherits a canonical pivotal
structure $a: id_{\mcal B} \ra \# \circ \# = id_{\mcal{B}}$, which is
identity on $\mcal{B}_0$, the objects $a_{\vlon} := {\vlon 0} \in
ob(\mcal{B}(\vlon, \vlon))$ and the morphisms $a_{\vlon k} := id_{\vlon
k} : \vlon k \ra \vlon k$. In fact, we have a correspondence between
the pivotal structures on $\mcal B$ and the weights of $P$.

\begin{prop}\label{weight-pivotal}
 There is a one-to-one correspondence between weights of a planar
algebra $P$ and pivotal structures on the strict $2$-category $\mcal
B$ associated to $P$ as above.
\end{prop}
\begin{pf}
Given a weight $z$ of $P$, define  $a_{\vlon k} := z_{(-)^k \vlon k} \in Mor(\vlon
k, \vlon k)$. For naturality of $a$, consider $f \in Mor(\vlon
k, \vlon l)$. Note that $a_{\vlon l}^{-1} \circ
f \circ a_{\vlon k} = z_{(-)^k \vlon \frac{k+l}{2}}^{-1} \circ f \circ
z_{(-)^k \vlon \frac{k+l}{2}} = f$ where the first equality readily follows
from pictures and the second one holds because $z$'s are central. Also,
the tensor condition $a_{(-)^k \vlon l} \otimes a_{\vlon k} = a_{\vlon
(k+l)}$ is easy to verify. Setting $a_{\vlon} = 1_{\vlon}$, we have a
weak transformation $ a: id_{\mcal B} \ra \# \circ \# =
id_{\mcal{B}}$.  This proves that $a$ gives a pivotal structure.

 Conversely, if $a: id_{\mcal B} \ra \# \circ \# = id_{\mcal B}$ is a
pivotal structure, set $z = a_{- 1} \in Mor(-1,-1) = P_{+1}$. Then,
the tensor property of $a$, namely $a_{(-)^k \vlon l} \otimes a_{\vlon
k} = a_{\vlon (k+l)}$, along with naturality of $a$ and
Proposition \ref{a-properties}$(2)$ implies that $z_{\vlon
k} \in \mcal{Z}(P_{\vlon k})$ and hence $z$ is a weight of $P$.
\end{pf}

\subsection{Bicategories to planar algebras}

We first briefly recall from \cite{Gho08} the planar algebra
 associated to any $1$-cell in a pivotal $\C$-linear strict
 $2$-category $\mcal B$.  Let $X \in ob(\mcal{B}(-,+) )$ for
 $\{+,-\} \subset \mcal{B}_0$, $\# : \mcal{B} \ra \mcal{B}^{op\,
 co} $ be a right rigid structure with respect to evaluation and
 coevaluation $e$ and $c$ respectively, and $a: id_{B} \ra \#\circ \#$
 be a pivotal structure. The ingredients of the planar algebra $P$
 associated to $(\mcal B, \#, a, X)$ are as follows:\\ {\em Vectors
 spaces:} $P_{\vlon k} := End(X_{\vlon k})$, where

\begin{equation*}
\begin{tabular}{ccc}
$X_{\vlon k}$ & $:=$ & $\left\{
\begin{tabular}{c}
$X\otimes X^{\# }\otimes X\otimes X^{\# }\otimes X\otimes \cdots k$ many
tensor factors if $\varepsilon =+$, \\
$X^{\# }\otimes X\otimes X^{\# }\otimes X\otimes X^{\# }\otimes \cdots
k$ many tensor factors if $\varepsilon =-$,
\end{tabular}
\right. $
\end{tabular}
\end{equation*}
if $k \geq 1$ and $X_{\vlon 0}: = 1_{\vlon} \in
ob(\mcal{B}(\vlon, \vlon))$.
\vspace*{2mm}

\begin{figure}[h]\vspace*{-3mm}
\hspace*{-105mm}
\psfrag{ar}{$\leadsto$}
\psfrag{k}{$k$}
\psfrag{x}{$x$}
\psfrag{e}{${\! \vlon}$}
\psfrag{x}{$x$}
\psfrag{m}{$_m$}
\psfrag{l}{$l$}
\psfrag{ev}{$id_{X_{(-)^k k}}\otimes \left( e_{X^{\#}} \circ (a_X \otimes id_{X^{\#}} )
\right) \otimes id_{X_{+ l}} 
\in Mor(X_{(-)^k (k+l+2)}, X_{(-)^k (k+l)})$}
\psfrag{e*}{$id_{X_{(-)^{k+1} k}}\otimes e_{X} \otimes id_{X_{-l}} 
\in Mor(X_{(-)^{k+1} (k+l+2)}, X_{(-)^{k+1} (k+l)} )$}
\psfrag{c}{$id_{X_{(-)^k k}}\otimes c_{X} \otimes id_{X_{+ l}} 
\in Mor(X_{(-)^k (k+l)}, X_{(-)^k (k+l+2)})$}
\psfrag{c*}{$id_{X_{(-)^{k+1} k}}\otimes\left( (id_{X^{\#}} \otimes a_{X}^{ -1} ) \circ c_{X^{\#}} \right) \otimes id_{X_{- l}} \in Mor(X_{(-)^{k+1} (k+l)}, X_{(-)^{k+1} (k+l+2)})$}
\psfrag{xin}{$id_{X_{(-)^k \vlon k}} \otimes x \otimes id_{X_{(-)^m \vlon l}} \in P_{(-)^k \vlon (k+l+m)} \text{ for all } x \in P_{\vlon m}$}
\includegraphics[scale=0.4]{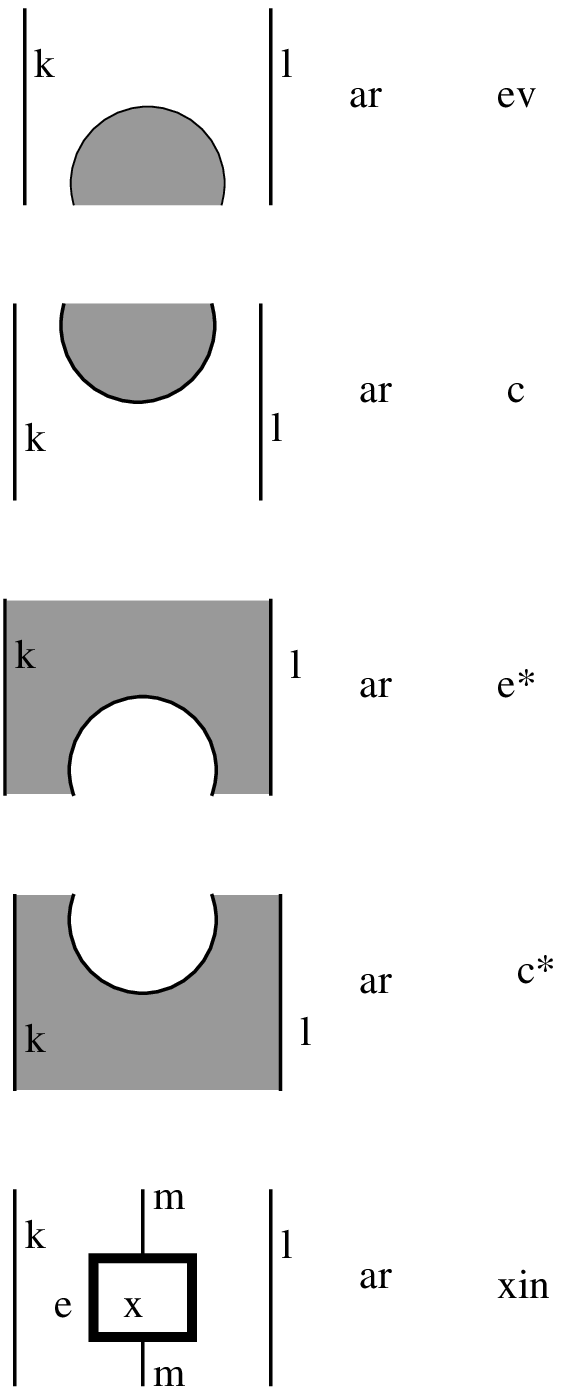}\vspace*{-3mm}
\caption{Prescription for bicategory planar algebra}\label{bicat-pa}
\end{figure}

\noindent {\em Actions of tangles:} Given a tangle $T: (\vlon_1 k_1, 
\ldots \vlon_b k_b) \ra \vlon_0 k_0$ and elements $x_i \in P_{\vlon_i k_i},\ 
1 \leq i \leq b$, the action $P_T(x_1,\ldots,x_b)$ is given by $(i)$
choosing a standard form representative $T_1$ in the isotopy class of
$T$ labelled with $x_i$ in the $i$-th internal box, $(ii)$ cutting
$T_1$ into horizontal stripes so that each stripe should have at most
one local maximum, minimum or (labelled) internal rectangle, $(iii)$
assigning a $2$-cell to each horizontal stripe as prescribed in
Figure \ref{bicat-pa} and $(iv)$ successively composing these
$2$-cells with the one coming from the bottom stripe being the rightmost
in the composition.

\begin{prop}
Let $a$ and $\tilde{a}$ be two pivotal structures on $\mcal B$ as
above, $P$ and $\tilde P$ be the planar algebras associated to a
$1$-cell $X$ with respect to $a$ and $\tilde a$ respectively. Then,
$P$ and $\tilde P$ are perturbations of each other.
\end{prop}
\begin{pf}
Set $z = {a}_X^{-1} \circ \tilde{a}_X \in P_{+1} = End(X)$. We assert
that $z$ is a weight of $P$ and  ${P}^{(z,1)} = \tilde P$.  It
follows from the definition of action of tangles and \cite[Lemma
$4.2$]{Gho08} that $z^{\#} = P_{
\psfrag{z}{$z$}
\psfrag{-}{$-$}
\includegraphics[scale=0.2]{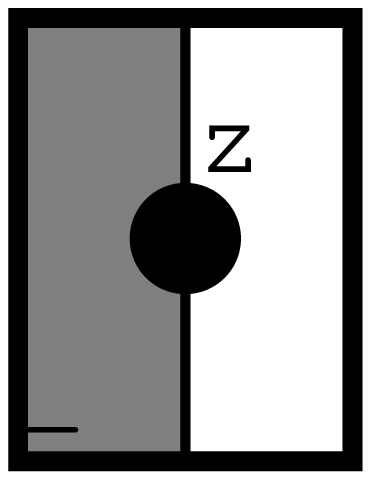}
}$. 
Note that $z_{+k} $ (as in Definition \ref{weight}) is given by
\begin{eqnarray*}
z_{+k} & = & z \otimes (z^{-1})^{\#} \otimes z \otimes \cdots
\ k\text{-tensors} \\ & = & (
a_{X} \otimes a_{X^{\#}} \otimes a_{X} \otimes \cdots )^{-1} \circ
(\tilde{a}_X \otimes
\tilde{a}_{X^{\#}} \otimes  \tilde{a}_X \otimes \cdots )\\
& = &
a_{X_{+k}}^{-1}\circ \tilde{a}_{X_{+k}}\ \in 
\mcal{Z}(P_{+k}),
\end{eqnarray*}
where in the second (resp., last) equality, we have used the pivotal
property of $a$ and $\tilde a$ as in Proposition \ref{a-properties}
(resp., Corollay \ref{t-property}). And on similar lines one
establishes that $z_{-k} \in \mcal{Z}(P_{-k})$. This shows that $z$ is
indeed a weight of $P$. Now, it is a matter of routine verification
that ${P}^{(z,1)}_T = \tilde P_T$ for a set of generating tangles as
in Figure \ref{tangles}. Thus, we conclude that ${P}^{(z,1)} = \tilde
P$.
\end{pf}

\section{Bimodule planar algebras}\label{bimod-pa}

\subsection{Planar algebra associated to a bimodule}
In this subsection, we associate a `bimodule planar algebra' to
a bifinite bimodule with a natural correspondence
between extremality and sphericality. This
will pave the way for us to associate a unimodular bimodule planar
algebra to a finite index subfactor, giving an extension of Jones'
Theorem \cite[Theorem 4.2.1]{Jon} to an arbitrary finite index
subfactor (not necessarily extremal).

\begin{defn}
A finite dimensional, connected, positive $C^*$-planar algebra is called
bimodule planar algebra.
\end{defn}

Let $_A {\mathcal H}_B$ be a finite index bimodule for $II_1$ factors
$A$ and $B$. Before stating the next theorem, we set up some notations
that will be used throughout this section:\\ ${\mathcal H}_{+ 0} :=
L^2 (A)$, ${\mathcal H}_{- 0} : = L^2 (B)$ and for $k\geq 1$,
${\mathcal H}_{\vlon k}$ is the tensor product (over $A$ or $B$) of
$k$-many modules $\mathcal H$ and $\overline{\mathcal H}$ alternately
with $\mathcal H$ (resp., $\overline{\mathcal H}$) being the left-most
module if $\vlon=+$ (resp., $-$).

\begin{thm}\label{bimod-pa-theorem}
Let $_A {\mathcal H}_B$ be a finite index bimodule for $II_1$-factors
$A$ and $B$. Then,\\
(i) $P$ defined by\vspace*{-2mm}
\[
P_{\vlon k} = \left\{
\begin{array}{ll}
{_A}{\mathcal L}_A ({\mathcal H}_{+ k}) \text{ or } {_A}{\mathcal L}_B
  ({\mathcal H}_{+ k}),  \text{ if } \vlon = +, \\  {_B}{\mathcal L}_B
  ({\mathcal H}_{- k})\text{ or }{_B}{\mathcal L}_A ({\mathcal
    H}_{- k}), \text{ if } \vlon = -,
\end{array}
\right.
\]
according as $k$ is even or odd, has a unique bimodule planar algebra
structure with $\ast$-structure coming from the usual adjoints of
intertwiners, satisfying:

(a) action of multiplication tangles matches with the composition
of operators in the intertwiner spaces,

(b) $P_{RI_{\vlon k}} = \left\{
\begin{tabular}{ll}
$P_{\vlon k}\ni T \mapsto T \underset{A}{\otimes} id_{\mathcal H} \in
  P_{\vlon (k+1)}$, & if either $\vlon = +$ and $k$ is even, or $\vlon
  = -$ and $k$ is odd,\\ $P_{\vlon k}\ni T \mapsto T
  \underset{B}{\otimes} id_{\overline{\mathcal H}} \in P_{\vlon
    (k+1)}$, & otherwise,
\end{tabular}
\right.$

(c) $P_{LI_{\vlon k}} = \left\{
\begin{tabular}{ll}
$P_{\vlon k}\ni T \mapsto id_{\overline{\mathcal H}}
  \underset{A}{\otimes} T \in P_{-\vlon (k+1)}$, & if $\vlon =
  +$,\\ $P_{\vlon k}\ni T \mapsto id_{\mathcal H} \underset{B}{\otimes}
  T \in P_{-\vlon (k+1)}$, & if $\vlon = -$ and
\end{tabular}
\right.$

(d) $P_{E_{+1}}$ (resp., $P_{E_{-1}}$) is given by 
\begin{align*}
 & {\mathcal H}^o \underset{B}{\otimes} {\overline{\mathcal H}}^o \ni
  \xi \underset{B}{\otimes} \overline{\eta} \overset{A-A}{\longmapsto}
  \underset{j}{\sum} {_A}\langle \eta , \xi \rangle \left( \eta_j
  \underset{B}{\otimes} \overline{\eta}_j \right) \in {\mathcal H}^o
  \underset{B}{\otimes} {\overline{\mathcal H}}^o\\ \text{(resp., } &
           {\overline{\mathcal H}}^o \underset{A}{\otimes} {\mathcal
             H}^o \ni \overline{\eta} \underset{A}{\otimes} \xi
           \overset{B-B}{\longmapsto} \underset{i}{\sum} \langle \eta
           , \xi \rangle_B \left( \overline{\xi}_i
           \underset{A}{\otimes} \xi_i \right) \in {\overline{\mathcal
               H}}^o \underset{A}{\otimes} {\mathcal H}^o \text{ )}
\end{align*}
where $\{\xi_i\}_{i}$ (resp., $\{\eta_j\}_{j}$) is any basis for
$_A{\mathcal H}$ (resp., ${\mathcal H}_B$). (We will refer $P$ (resp.,
normalized $P$) as the {\em bimodule} (resp., {\em normalized} or {\em
  unimodular}) {\em planar algebra associated to $_A{\mathcal H}_B$}).

(ii) $P$ (as in (i)) is spherical if and only if $_A {\mathcal H}_B$
is extremal (as in $\S$\ref{bimod-bicat}).

(iii) If $B' := {\mathcal L}_B ({\mathcal H})$, then the normalized
planar algebras associated to $_A {\mathcal H}_B$ and $_A L^2
(B')_{B'}$ are isomorphic as $\ast$-planar algebras.
\end{thm}

\begin{pf} 
$(i)$ The planar algebra (in \cite{Gho08}) associated to the $1$-cell
  $_A\mcal{H}_B$ in the pivotal bicategory of bifinite bimodules ($\S$
  \ref{bimod-bicat}) has same vector spaces as of $P$; we provide $P$
  with the same planar algebra structure, the prescription for the
  actions of tangles for which is given in Figure \ref{bicat-pa}. It
  is clear from this prescription that $P$ satisfies the conditions
  $(a)$-$(d)$. Since $P$ is connected, the fact that $P$ with the
  above mentioned $\ast$-structure is a $C^\ast$-planar algebra, can
  be verified readily by checking the $\ast$-condition
  \eqref{*-condition} for the tangles $M_{\vlon k}, RI_{\vlon k},
  LI_{\vlon k}, E_{\vlon k}$. It now remains to show that the
  $C^*$-planar algebra $P$ is positive. For this, it is enough to show
  that the tangle maps $P_{RE_{\vlon k}}$ and $P_{LE_{\vlon k}}$ are
  positive definite for all $\vlon k \in Col$. We prove this only for
  $P_{RE_{+2k}}$; the others can be verified using similar
  arguments. For each $x \in P_{+2k} = \!_A\mcal{L}_A(\mcal{H}_{+2k})$
  and $\zeta \in \mcal{H}_{+(2k-1)}^o$, we have
\begin{eqnarray*}
P_{RE_{+2k}}(x)(\zeta) & = &\sum_i\left((id_{\mcal{H}_{+(2k-1)}}
\uset{B}{\otimes} e_{\mcal{H}}) \circ (x \uset{A}{\otimes} id_{\mcal
  H} ) \right) ( \zeta \uset{B}{\otimes} \bar{\xi}_i \uset{A}{\otimes}
\xi_i ) \\ & = & \sum_{i, i', k} (id_{\mcal{H}_{+(2k-1)}}
\uset{B}{\otimes} e_{\mcal{H}}) \left( \gamma_k \uset{B}{\otimes}
\bar{\xi}_i \langle \gamma_k \uset{B}{\otimes} \bar{\xi}_{i'}, x(\zeta
\uset{B}{\otimes} \bar{\xi}_i ) \rangle_A \uset{A}{\otimes} \xi_i
\right) \\ & = & \sum_{i, i', k} \gamma_k \left\langle {\xi}_{i'}
, \langle \gamma_k \uset{B}{\otimes} \bar{\xi}_{i'}, x(\zeta
\uset{B}{\otimes} \bar{\xi}_i ) \rangle_A \xi_i \right \rangle_B,
\end{eqnarray*}
where $\{ \gamma_k\}_k$ is a basis for $(\mcal{H}_{+(2k-1)})_B$. This
gives
\[
\langle \zeta,\, P_{RE_{+2k}}(x)(\zeta) \rangle = \sum_{k, i, i'}\left
\langle \zeta \uset{B}{\otimes} \oline{ \langle \gamma_k
  \uset{B}{\otimes} \bar{\xi}_{i'}, x(\zeta \uset{B}{\otimes}
  \bar{\xi}_i ) \rangle_A \xi_i}, \gamma_k \uset{B}{\otimes}
\bar{\xi}_{i'} \right \rangle = \sum_i\langle \zeta \uset{B}{\otimes}
\oline{\xi}_i,\, x(\zeta \uset{B}{\otimes} \oline{\xi}_i ) \rangle.
\]
Thus, $P_{RE_{+2k}}(x)$ is positive semi-definite. In addition, the
above also gives
\[
\sum_k\langle \gamma_k,\, P_{RE_{+2k}}(x)(\gamma_k) \rangle =
\sum_{k,i}\langle \gamma_k \uset{B}{\otimes} \oline{\xi}_i,\, x(\gamma_k
\uset{B}{\otimes} \oline{\xi}_i ) \rangle = dim (\mcal{H}_{+2k})_A
\;  tr_C(x),
\]
where $C:= \mcal{L}_A(\mcal{H}_{+2k})$ and the last equality is a
consequence of Proposition \ref{trace-A'}. This proves that $
P_{RE_{+2k}}$ is also faithful. Thus, $P$ is a bimodule planar
algebra; and for uniqueness, we appeal again to
Remark \ref{planar-morphism}.
\vspace*{1mm}

$(ii)$ It is enough to prove that the normalized left and the right
picture traces on $P_{+1}$ are given by the unique traces on the
$II_1$-factors $A' := {_A}\mathcal{L}(\mcal{H})$ and $B':=
\mathcal{L}_B(\mcal{H})$ respectively. We will only exhibit a proof
for the left one. First, note that
\begin{align*}
 \psfrag{+}{$\!\! _+$}
  \psfrag{-}{$_-$} P_{\ \includegraphics[scale=0.35]{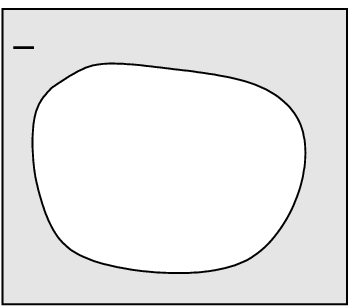}}
& =
e_{\mathcal H} \circ (id_{\overline{\mathcal H}} \otimes
a^{-1}_{\mathcal H}) \circ c_{\overline{\mathcal H}} =
\underset{i}{\sum} \langle \xi_i,\xi_i \rangle_B
\hat{}=dim(_A{\mathcal H}) 1_{P_{-0}}\text{ and } \\
 \psfrag{+}{$\!
  _+$} \psfrag{-}{$\! \! _-$}
  P_{\ \includegraphics[scale=0.35]{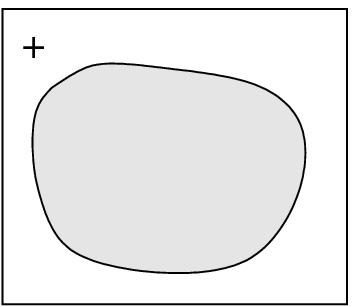}} & =
e_{\overline{\mathcal H}} \circ (a_{\mathcal H} \otimes
id_{\overline{\mathcal H}}) \circ c_{\mathcal H} = \underset{j}{\sum}
{_A} \langle \eta_j,\eta_j \rangle\hat{}=dim({\mathcal H}_B)
1_{P_{+0}}
\end{align*}
where $\{\xi_i\}_{i}$ (resp., $\{\eta_j\}_{j}$) is a basis for
$_A{\mathcal H}$ (resp. ${\mathcal H}_B$).  So, the modulus of the
planar algebra $P$ is $( dim {_A\mathcal H},\, dim {\mathcal H}_B
)$. Further, for each $x \in P_{+1}={_A}{\mathcal L}_B ({\mathcal
  H})$,
\[
{\mathcal H}_{-0} := L^2 (B) \ni \hat{1} \overset{P_{LE_{+1}}
  (x)}{\longmapsto} \underset{i}{\sum} \langle \xi_i , x \xi_i
\rangle\hat{}_B = \underset{i}{\sum} \langle \xi_i , x \xi_i \rangle\hat{1} =
dim({_A}{\mathcal H}) tr_{A'} (x)\hat{1} \in L^2(B) =: {\mathcal
H}_{-0}.
\]
The first equality follows from the fact that
$\underset{i}{\sum} \langle \xi_i , x \xi_i
\rangle_B$ commutes with $B$ and the second one comes again from 
 Proposition \ref{trace-A'}.

\vspace*{1mm}

$(iii)$ Note that the planar algebra associated to $_A{\mathcal H} _B$
may very well be non-unimodular but the index of $P$ is same as that
of $_A{\mathcal H}_B$.  Set $\delta_{-} = dim({_A}{\mathcal H} )$,
$\delta_{+} = dim({\mathcal H}_B )$ and $\delta = \sqrt{\delta_{+}
  \delta_{-}}$.  Let $Q$ denote the bimodule planar algebra associated
to $_A L^2 (B')_{B'}$. So, $\tilde{P}:=
P^{\left(\sqrt{\frac{\delta_{-}}{\delta_{+}}}, 1\right)}$ and
$\tilde{Q}:= Q^{(\delta,1)}$ are the normalizations of $P$ and $Q$
respectively. To establish an isomorphism between $\tilde{P}$ and
$\tilde{Q}$, we consider the following $B'$-$B'$ linear maps
$u:{\mathcal H} \underset{B}{\otimes} \overline{{\mathcal H}}
\rightarrow {\mathcal K} := L^2(B')$ and $w:{\mathcal H}
\underset{B}{\otimes} \overline{{\mathcal H}} \rightarrow
\overline{\mathcal K}$ given by
\begin{align*}
& {\mathcal H}^o \underset{B}{\otimes} \overline{{\mathcal H}}^o \ni \xi
\underset{B}{\otimes} \overline{\eta} \overset{u}{\mapsto}
         {_{B'}}\langle \eta , \xi \rangle \hat{1}\in {\mathcal K} \;
         \text{ implying } \; {\mathcal K} \ni \hat{1}
         \overset{u^*}{\mapsto} \delta^{-1}_{+} \underset{j}{\sum}
         \left( \eta_j \underset{B}{\otimes} \overline{\eta}_j \right)
         \in {\mathcal H}^o \underset{B}{\otimes} \overline{{\mathcal
             H}}^o \text{ and}\\
& {\mathcal H}^o \underset{B}{\otimes} \overline{{\mathcal H}}^o \ni \xi
\underset{B}{\otimes} \overline{\eta} \overset{w}{\mapsto} {_{B'}}
\langle \eta , \xi \rangle \overline{\hat{1}} \in \overline{\mathcal
  K} \; \text{ implying } \; \overline{\mathcal K} \ni \overline{\hat{1}}
\overset{w^*}{\mapsto} \delta^{-1}_{+} \underset{j}{\sum} \left(
\eta_j \underset{B}{\otimes} \overline{\eta}_j \right) \in {\mathcal
  H}^o \underset{B}{\otimes} \overline{{\mathcal H}}^o.
\end{align*}
Clearly, $u$ and $w$ are co-isometries; further, $index ( {_{B'}}
\mcal{H} \underset{B}{\otimes}\overline{ \mcal{H}} {_{B'}} ) = 1$
implies $u$ and $w$ are unitaries as well.  Let $v: {\mathcal H}
\underset{B}{\otimes} \overline{{\mathcal H}} \rightarrow {\mathcal K}
\underset{B'}{\otimes} \overline{{\mathcal K}}$ be the $B'$-$B'$
linear unitary obtained from $u$ after composing with the unitary
taking $\hat{1} \in {\mathcal K}$ to $( \hat{1} \underset{B'}{\otimes}
\bar{\hat{1}} ) \in {\mathcal K} \underset{B'}{\otimes}
\overline{{\mathcal K}}$.  Define $\varphi : \tilde{P} \rightarrow
\tilde{Q}$ in the following way:
\[
\tilde{P}_{+ 2l} = {_A}{\mathcal L}_A ({\mathcal H}_{+2l}) \ni x
\overset{\varphi_{+2l}}{\longmapsto} Ad_{v_l} (x) \in {_A}{\mathcal L}_A
({\mathcal K}_{+2l})=\tilde{Q}_{+2l} \text{ for $l \geq 1$,}
\]
where $v_l := \left( v \underset{A}{\otimes} v \underset{A}{\otimes}
\cdots l \text{ tensors} \right) \in {_{B'}}{\mathcal L}_{B'}
({\mathcal H}_{+2l},{\mathcal K}_{+2l})$, and for $l \geq 0$,
\begin{align*}
\tilde{P}_{+(2l+1)} = {_A}{\mathcal L}_B ({\mathcal H}_{+(2l+1)}) \ni
x & \overset{\varphi_{+(2l+1)}}{\longmapsto}   Ad_{(v_l \underset{A}{\otimes}
  u)} \left(x \underset{B}{\otimes} id_{\overline{\mathcal H}} \right)
\in {_A}{\mathcal L}_{B'} ({\mathcal
  K}_{+(2l+1)})=\tilde{Q}_{+(2l+1)},
\\
\tilde{P}_{-(2l+1)} = {_B}{\mathcal L}_A ({\mathcal H}_{-(2l+1)}) \ni
x & \overset{\varphi_{-(2l+1)}}{\longmapsto}  Ad_{(w \underset{A}{\otimes}
  v_l)} \left(id_{\mathcal H} \underset{B}{\otimes} x \right) \in
{_{B'}}{\mathcal L}{_A} ({\mathcal K}_{-(2l+1)})=\tilde{Q}_{-(2l+1)},
\\
\tilde{P}_{-(2l+2)} = {_B}{\mathcal L}_B ({\mathcal H}_{-(2l+2)}) \ni
x &\overset{\varphi_{-(2l+2)}}{\longmapsto} Ad_{(w \underset{A}{\otimes}
  v_l \underset{A}{\otimes} u)} \left(id_{\mathcal H}
\underset{B}{\otimes} x \underset{A}{\otimes} id_{\overline{\mathcal
    H}} \right) \in {_{B'}}{\mathcal L}_{B'} ({\mathcal
  K}_{-(2l+2)})=\tilde{Q}_{-(2l+2)}.
\end{align*}
Clearly, each $\varphi_{\vlon k}$ is an injective $\ast$-algebra
homomorphism and each $\varphi_{+ 2l}$ is surjective. To see that
$\varphi_{+(2l+1)}$ (resp., $\varphi_{-(2l+1)}$) is surjective, note
that
\begin{align*}
& Ad_{(v_l \underset{A}{\otimes} u)^*} \left( {_A}{\mathcal L}_{B'}
  ({\mathcal K}_{+(2l+1)}) \right) = {_A}{\mathcal L}_{B'} ({\mathcal
    H}_{+(2l+1)} \underset{B}{\otimes} \overline{\mathcal H})=\left\{x
  \underset{B}{\otimes} id_{\overline{\mathcal H}} : x \in
           {_A}{\mathcal L}_B ({\mathcal H}_{+(2l+1)})
           \right\}\\ \text{(resp., }& Ad_{(w \underset{A}{\otimes}
             v_l)^*} \left({_{B'}}{\mathcal L}{_A} ({\mathcal
             K}_{-(2l+1)}) \right) = {_{B'}}{\mathcal L}{_A}
           ({\mathcal H} \underset{B}{\otimes} {\mathcal H}_{-(2l+1)})
           = \left\{ id_{\mathcal H} \underset{B}{\otimes} x : x \in
           {_B}{\mathcal L}{_A} ({\mathcal H}_{-(2l+1)}) \right\}
           \text{)}
\end{align*} 
where the second equality follows from Lemma
\ref{basic-con-intertwiners} since ${_B}\overline{\mathcal H}_{B'}$
(resp., ${_{B'}}{\mathcal H}_B$) has index $1$. Same arguments also
imply $\varphi_{-(2l+2)}$ is surjective. Thus, $\varphi$ is an
isomorphism preserving actions of multiplication tangles.
\vspace{2mm}

\noindent{\em Actions of inclusion tangles:} We only show that $\varphi$
is equivariant with respect to the action of $RI_{+k}$; the same for
$RI_{-k}$ and $LI_{\vlon k}$ can be verified along similar lines. For
$x \in \tilde{P}_{+(2l+1)}$, note that
\begin{align*}
[Ad_{v_l \uset{A}{\otimes} u}(x \uset{B}{\otimes}
  id_{\overline{\mathcal H}})] \uset{B'}{\otimes}
id_{\overline{\mathcal K}} =Ad_{v_l \uset{A}{\otimes} u
  \uset{B'}{\otimes} id_{\overline{\mathcal K}}}(x \uset{B}{\otimes}
id_{\overline{\mathcal H}} \uset{B'}{\otimes} id_{\overline{\mathcal
    K}}) =Ad_{v_l \uset{A}{\otimes} v \circ \theta}(x
\uset{B}{\otimes} id_{\overline{\mathcal H}} \uset{B'}{\otimes}
id_{\overline{\mathcal K}}) =Ad_{v_l \uset{A}{\otimes} v}(x
\uset{B}{\otimes} id_{\overline{\mathcal H}}),
\end{align*}
where ${\mcal H} \uset{B}{\otimes} \overline{\mcal H}
\uset{B'}{\otimes} \overline{\mcal K} \ni \xi \uset{B}{\otimes}
\overline{\eta} \uset{B'}{\otimes} \bar{\hat{1}}
\overset{\theta}{\mapsto} \xi \uset{B}{\otimes} \overline{\eta} \in
        {\mcal H} \uset{B}{\otimes} \overline{\mcal H}$ is a $B'$-$B'$
        unitary intertwiner. This implies
        $\varphi\left(\tilde{P}_{RI_{+(2l+1)}}(x)
        \right)=\tilde{Q}_{RI_{+(2l+1)}}(\varphi(x))$. Equivariance of
        $\varphi$ with respect to the action of $RI_{+2l}$ is even more easier.
\vspace{2mm}

\noindent{\em Action of $E_{\pm 1}$:} From condition $(d)$ of part $(i)$ and 
definition of perturbation, we obtain
\begin{align*}
\varphi(\tilde{P}_{E_{+1}}) (\hat{x} \uset{B'}{\otimes} \bar{\hat{y}}) &
= \frac{1}{\delta_+} \sum_{j} v\tilde{P}_{E_{+1}} (xy^* \eta_j
\uset{B}{\otimes} \overline{\eta}_j) =
\frac{1}{\delta_+}\sqrt{\frac{\delta_-}{\delta_+}} \sum_{j, j'} v
\left(\,_A\langle xy^* \eta_j, \eta_j\rangle \eta_{j'} \otimes
\overline{\eta_{j'}} \right)\\ &=
\frac{1}{\delta_+}\sqrt{\frac{\delta_-}{\delta_+}} \sum_{j, j'}
E_A\left(\, _{B'}\langle xy^* \eta_j, \eta_j\rangle\right) \eta_{j'}
\otimes \overline{\eta_{j'}} =
\frac{1}{\delta_+}\sqrt{\frac{\delta_-}{\delta_+}} E_A(xy^*)
(dim\mcal{H}_B)^2 \, \hat{1}\uset{B'}{\otimes}\bar{\hat 1} \\ &=
\tilde{Q}_{E_{+1}}(\hat{x}\uset{B'}{\otimes}\hat{y}) 
\end{align*}
for all $x , y \in B'$. This gives $\varphi(\tilde{P}_{E_{+1}})
=\tilde{Q}_{E_{+1}}$. On the other hand, since
$\varphi(\tilde{P}_{E_{-1}})$ and $\tilde{Q}_{E_{-1}}$ are both
$B'$-$B'$ linear, in order to show their equality, it suffices to show
that $\varphi(\tilde{P}_{E_{-1}}) (\bar{\hat{1}} \uset{A}{\otimes}
\hat{1} ) =\tilde{Q}_{E_{-1}} (\bar{\hat{1}} \uset{A}{\otimes} \hat{1}
)$. Again, from condition $(d)$ of part $(i)$ and definition of
perturbation, we have
\begin{align*}
\varphi(\tilde{P}_{E_{-1}}) (\bar{\hat{1}} \uset{A}{\otimes} \hat{1}
)& = \delta_+^{-2} \sum_{j, j'} \left( (w \uset{A}{\otimes} u) \circ (
id_{\mcal{H}} \uset{B}{\otimes} \tilde{P}_{E_{-1}} \uset{A}{\otimes}
id_{\oline{\mcal{H}}} ) \right) ( \eta_{j} \uset{B}{\otimes}
\overline{\eta}_{j} \uset{A}{\otimes} \eta_{j'} \uset{B}{\otimes}
\overline{\eta}_{j'})\\ & = (\delta \delta_+)^{-1} \sum_{j,j', i} (w
\uset{A}{\otimes} u) ( \eta_{j} \uset{B}{\otimes} \langle \eta_j,
\eta_{j'} \rangle_B \overline{\xi}_{i} \uset{A}{\otimes} \xi_{i}
\uset{B}{\otimes} \overline{\eta}_{j'})\\ & = (\delta \delta_+)^{-1}
\sum_{j', i} \oline{ _{B'}\langle \eta_{j'} ,\, \xi_i
  \rangle^{\widehat{}}} \uset{A}{\otimes} {_{B'}\langle} \eta_{j'},
\xi_i \rangle^{\widehat{}} = (\delta \delta_+)^{-1} \sum_{k ,j', i}
\oline{ E_A( {_{B'}\langle} \eta_{j'}, \xi_i \rangle c_k^*) \hat{c}_k}
\uset{A}{\otimes} {_{B'}\langle} \eta_{j'}, \xi_i \rangle^{\widehat{}}
\\ & = (\delta \delta_+)^{-1} \sum_{k ,j', i} \bar{\hat{c}}_k
\uset{A}{\otimes} {_{B'}\langle} \eta_{j'} , { _{A}\langle} \xi_i ,
c_k \eta_{j'} \rangle \xi_i \rangle^{\widehat{}} = \delta ^{-1}
\sum_{k } \bar{\hat{c}}_k \uset{A}{\otimes}\, \hat{c}_k \\ &=
\tilde{Q}_{E_{-1}}( \bar{\hat{1}} \uset{A}{\otimes} \hat{1}),
\end{align*}
where $\{ c_k\}$ is a left basis for the subfactor $A \subset
B'$.  Thus, in view or Remark \ref{planar-morphism},
$\varphi$ is an isomorphism of $\ast$-planar algebras.
\end{pf}
\vspace*{2mm}

We now proceed to associate a unimodular bimodule planar algebar to a
finite index type $II_1$ subfactor $N \subset M$. For this, the
obvious thing to consider, is the unimodular bimodule planar algebra
associated to ${_N} L^2 (M)_M$. However, we would like to find out the
actions of tangles on the relative commutants (instead of intertwiner
spaces) and a set of conditions which uniquely determines the action
exactly the way \cite[Theorem 4.2.1]{Jon} states for extremal
subfactors. Before doing so, we first set up some notations and recall
certain standard subfactor theory facts - see \cite{Bis97, Jon83,
  JS97, PP88, PP86}.

Let $N \subset M$ be a subfactor with $\delta^2 : = [M:N] < \infty$
$(\delta > 0)$ and $\{M_k\}_{k \geq 1}$ be a tower of basic
constructions with $\{e_k \in \mscr{P} (M_k) \}_{k \geq 1}$ being a
set of Jones projections. We will have instances to apply the
following useful fact.

\begin{lem}\label{uniqueness-lemma}\cite{PP86}
For each $x_1 \in M_1$, there is a unique $x \in M$ satisfying $x_1
e_1 = x e_1$; this unique element is given by $x = [M:N] E_M(x_1
e_1)$.
\end{lem}

 For each $k \geq 1$, set $e_{[-1,k]} = \delta^{k (k+1)}(e_{k+1}e_k
 \cdots e_1) (e_{k+2} e_{k+1}\cdots e_2)\cdots (e_{2k+1} e_{2k}$
 $\cdots e_{k+1}) \in N^\prime \cap M_{2k+1}$, $ e_{[0, k]} =\delta^{k
   (k-1)} (e_{k+1}e_k \cdots e_2) (e_{k+2} e_{k+1}\cdots e_2)\cdots
 (e_{2k} e_{2k-1} \cdots e_{k+1}) \in M^\prime \cap M_{2k}$ and $v_k =
 \delta^k e_k e_{k-1} \cdots e_1 \in N^\prime \cap M_k$.  Then, the
 tower of $II_1$ factors $N \subset M_k\subset M_{2k+1}$ (resp., $M
 \subset M_k \subset M_{2k}$) is an instance of basic construction
 with ${e_{[-1,k]}}$ (resp., ${e_{[0,k]}}$) as Jones projection, that
 is, there exists an isomorphism $ {\varphi_{-1, k}} : M_{2k+1} {\lra}
 \mcal{L}_N(L^2(M_k))$ (resp., ${\varphi_{0, k}} : M_{2k} {\lra}
 \mcal{L}_M(L^2(M_k))$) given by
\begin{eqnarray*}
   \varphi_{-1, k}(x_{2k+1}) \hat{x}_{k} & = & \delta^{2(k+1)}
   E_{M_k}(x_{2k+1}x_{k}e_{[-1,k]})^{\widehat{}}\\ \text{(resp., }
   \varphi_{0,k} (x_{2k}) \hat{x}_{k} &= &\delta^{2k}
   E_{M_k}(x_{2k}x_{k}e_{[0,k]})^{\widehat{}} \text{ )}
\end{eqnarray*}
for all $x_i \in M_i$, $i = k,\, 2k,\, 2k+1$, which is identity
restricted to $M_k$ and sends $e_{[-1,k]}$ (resp., $e_{[0,k]}$) to the
projection with range $L^2(N)$ (resp., $L^2(M)$).  \comments{where
  $e_{[0,0]}:= 1_M $, give $\ast$-isomorphisms of towers $ ( N \subset
  M_k \subset^{e_{[-1,k]}} M_{2k+1}) \oset{\varphi_{-1,k}}{\cong} (N
  \subset M_k \subset^{\varphi_{-1,k}(e_{[-1,k]})}
  \mcal{L}_N(L^2(M_k)) )$ and $ (M \subset M_k \subset^{e_{[0,k]}}
  M_{2k}) \oset{\varphi_{0,k}}{\cong} (N \subset M_k
  \subset^{\varphi_{0,k}(e_{[0,k]})} \mcal{L}_M(L^2(M_k))) $.  It is
  easily seen that ${\varphi_{-1, k}}_{|_{M_{2k}}} = {\varphi_{0,
      k}}$.}  Also, $\varphi_{-1, k} (M_i'\cap M_{2k+1}) =\,
_{M_i}\mcal{L}_N(L^2(M_k))$ (resp., $\varphi_{0, k}(M_i'\cap M_{2k}) =
\, _{M_i}\mcal{L}_M(L^2(M_k))$) and $\vphi_{0,k} = \left. \vphi_{-1,k}
\right|_{M_{2k}}$ for all $k \geq 0 $, $-1 \leq i \leq k$. Further,
for each $k \geq 0$, we have an $M_k$-$M$ linear unitary given (on 
bounded vectors) by
\[
 L^2(M_k) \uset{N}{\otimes} L^2(M) \ni \hat{y} \uset{N}{\otimes}
 \hat{z} \oset{u_k }{\longmapsto} (y v_{k+1}z)^{\widehat{}} \in
 L^2(M_{k+1})
\]
for all $y \in M_k$ and $z \in M$. This unitary also satisfies the
equation $ \varphi_{0,k+1}(x) = u_{k}(\varphi_{-1,k}(x)
\uset{N}{\otimes} id_{L^2(M)}) u_k^*$  for all $ x \in M_{2k+1}$.

Apart from these, if ${_N}\mcal{H}_M$ denotes the bimodule
${_N}L^2(M)_M$, then for each $ k \geq 1$, we have a bunch of unitary
intertwiners (determined by the following actions on the bounded
vectors)
\[
\begin{array}{rcl}
\mcal{H}_{+
2k} \ni \hat{x}_1 \uset{M}{\otimes} \bar{\hat{x}}_2 \uset{N}{\otimes}
 \cdots \uset{M}{\otimes} \bar{\hat{x}}_{2k}
& \oset{w_{+2k}}{\uset{N\text{-}N}{\longmapsto}} & (x_1
x_2^*v_1x_3x_4^*v_2x_5 \cdots
x_{2k-2}^*v_{k-1}x_{2k-1}x_{2k}^*)^{\widehat{}} \in L^2(M_{k-1}),\\
&&\\
\mcal{H}_{-
2k} \ni \bar{\hat{x}}_1 \uset{N}{\otimes}
{\hat{x}}_2 \uset{M}{\otimes} \cdots \uset{N}{\otimes} {\hat{x}}_{2k}
& \oset{w_{-2k}}{\uset{M\text{-}M}{\longmapsto}}
& ( x_1^*v_1x_2x_3^*v_2x_4 \cdots
x_{2k-3}^*v_{k-1}x_{2k-2}x_{2k-1}^*v_kx_{2k})^{\widehat{}} \in
L^2(M_{k}),\\
&&\\
\mcal{H}_{+
(2k-1)} \ni \hat{x}_1 \uset{M}{\otimes} \bar{\hat{x}}_2 
\uset{N}{\otimes} \cdots
\uset{N}{\otimes}  \hat{x}_{2k-1}  
& \oset{w_{+(2k-1)}}{\uset{N\text{-}M}{\longmapsto}}& (x_1
x_2^*v_1x_3x_4^*v_2x_5 \cdots
x_{2k-2}^*v_{k-1}x_{2k-1})^{\widehat{}} \in L^2(M_{k-1}),\\
&&\\
\mcal{H}_{-(2k-1)} \ni \bar{\hat{x}}_1 \uset{N}{\otimes}
{\hat{x}}_2 \uset{M}{\otimes} \cdots
\uset{M}{\otimes}  \bar{\hat{x}}_{2k-1} 
& \oset{w_{-(2k-1)}}{\uset{M\text{-}N}{\longmapsto}} & (
x_1^*v_1x_2x_3^*v_2x_4 \cdots
x_{2k-3}^*v_{k-1}x_{2k-2}x_{2k-1}^*)^{\widehat{}} \in L^2(M_{k-1})
\end{array}
\]
 and a very useful formula (see \cite{Jon}) $x_1 v_1 x_2 v_2 \cdots
 v_k x_{k+1} = x_1 v^\ast_k x_2 v^\ast_{k-1} \cdots v^\ast_1 x_{k+1}$
 for all $x_i \in M$, $1 \leq i \leq 2k$.

We are now ready to present the extension of Jones' Theorem
\cite[Theorem 4.2.1]{Jon}, which associates a unimodular bimodule
planar algebra to a finite index subfactor and gives a natural
correspondence between extremality and sphericality.
\begin{thm}\label{jones-theorem}
Let $M_{-1} := N \subset M =: M_0$ be a subfactor with $\delta^2 : = [M:N] < \infty$
($\delta >0$) and $\{M_k\}_{k \geq 1}$ be a tower of basic
constructions with $\{e_k \in \mscr{P} (M_k) \}_{k \geq 1}$ being a
set of Jones projections. Then, $P$ defined by $P_{\vlon k} = N' \cap
M_{k-1} $ or $M' \cap M_{k}$ according as $\vlon = +$ or $-$, has a
unique unimodular bimodule planar algebra structure with the
$\ast$-structure given by the usual $\ast$ of the relative commutants
such that for each $k \in \N_0$,
\begin{enumerate}
\item the action of multiplication tangles is given by the usual
multiplication in the relative commutants,
\item the action of the left inclusion tangle $LI_{-k}$ is given by
  the usual inclusion $M' \cap M_k \subset N'\cap M_k$,
\item the action of the right inclusion tangle $RI_{+k}$ is given by
  the usual inclusion $ M_{k-1} \subset M_{k}$,
\item $P_{E_{+(k+1)}} = \delta e_{k+1}$,
\comments{\item $P_{RE_{+(k+1)}}  = \delta E^ {N'\cap M_k}_{ N'\cap M_{k-1}}$,}
\item $P_{LE_{+(k+1)}} = \delta^{-1} \us{i}{\sum} b_i^* x b_i $ for all $ x
  \in P_{+(k+1)}$,
\comments{ and
\item $P_{TR^r_{+(k+1)}}(x) = \delta^{(k+1)} tr_{M_{k}}(x)$ for all $x
  \in P_{+(k+1)}$,}
\end{enumerate}
where $\{ b_i\}_i$ is a left basis for the subfactor $N \subset M$.
In particular, $P$ is spherical if and only if the subfactor $N
\subset M$ is extremal. ($P$ will be referred as the planar algebra
associated to the tower $\{M_k\}_{k \geq -1}$ with Jones projections
$\{e_k\}_{k \geq 1}$.)
\end{thm}
\begin{pf}
Let $Q$ denote the bimodule planar algebra associated to the bifinite
bimodule ${_N}\mcal{H}_M := {_N}L^2(M)_M$ (as in Theorem
\ref{bimod-pa-theorem}) and $\tilde{Q} = Q^{(\delta, 1)}$ be its
normalization. With notations as above, we have the $\ast$-algebra
isomorphisms
\[
\begin{array}{rccl}
P_{+2k} = N^\prime \cap M_{2k-1} \ni x &
\oset{\chi_{+2k}}{\longmapsto} & Ad_{w^\ast_{+ 2k}} \left(
\vphi_{-1,k-1} (x) \right) & \in {_N}{\mcal L}_N ({\mcal H}_{+2k}) =
Q_{+2k} = \tilde{Q}_{+ 2k},\\ P_{-2k} = M^\prime \cap M_{2k} \ni x &
\oset{\chi_{-2k}}{\longmapsto} & Ad_{w^\ast_{- 2k}} \left( \vphi_{0,k}
(x) \right) & \in {_M}{\mcal L}_M ({\mcal H}_{-2k}) = Q_{-2k} =
\tilde{Q}_{- 2k},\\ P_{+(2k+1)} = N^\prime \cap M_{2k} \ni x &
\oset{\chi_{+(2k+1)}}{\longmapsto} & Ad_{w^\ast_{+ (2k+1)}} \left(
\vphi_{0,k} (x) \right) & \in {_N}{\mcal L}_M ({\mcal H}_{+(2k+1)}) =
Q_{+(2k+1)} = \tilde{Q}_{+ (2k+1)},\\ P_{-(2k+1)} = M^\prime \cap
M_{2k+1} \ni x & \oset{\chi_{-(2k+1)}}{\longmapsto} &
Ad_{w^\ast_{-(2k+1)}} \left( \vphi_{-1,k} (x) \right) & \in {_M}{\mcal
  L}_N ({\mcal H}_{-(2k+1)}) = Q_{-(2k+1)} = \tilde{Q}_{-(2k+1)}
\end{array}
\]
for $k \in \N_0$.  We provide $P$ with a
unimodular bimodule planar algebra structure from that of $\tilde{Q}$
as follows:\\ For each tangle $T: (\vlon_1 k_1\times \cdots \times
\vlon_b k_b) \ra \vlon_0 k_0$ (resp., $T : \emptyset \ra \vlon_0
k_0$), we define $P_T = \chi_{\vlon_0 k_0}^{-1} \circ \tilde{Q}_T
\circ ( \times_{i = 1}^b \chi_{\vlon_i, k_i})$ (resp., $P_T =
\chi_{\vlon_0 k_0}^{-1} ( \tilde{Q}_T)$). Thus, $P$ inherits a
unimodular bimodule planar algebra structure with modulus $( \delta,
\delta )$. We now show that $P$ satisfies all the conditions in the
statement, which forces the planar algebra structure on $P$ to be
unique by Remark \ref{planar-morphism} (2). To begin with, note that,
$(1)$ needs no further verification. We will establish the relations
in $(2)-(5)$ only for even $k$'s because the proofs for odd
$k$'s, are completely analogous to those in the even case.

(2) Suppose $k= 2l$ and $y \in P_{-2l} = M'\cap M_{2l}$. Unravelling
the definitions, we just need to show that $ \varphi_{0,l}(y) =
Ad_{w_{+(2l+1)}}(id_{\mcal{H}} \uset{M}{\otimes} \chi_{-2l}(y))$. For
$\xi = \hat{x}_1\uset{M}{\otimes} \bar{\hat{x}}_2\uset{N}{\otimes}
\cdots \uset{N}{\otimes} \hat{x}_{2l-1} \uset{M}{\otimes}
\bar{\hat{x}}_{2l} \uset{N}{\otimes} \hat{x}_{2l+1} \in
\mcal{H}_{+(2l+1)}^o$ where $x_i \in M$ for all $1 \leq i \leq 2l+1$,
we have
\begin{align*}
& w_{+(2l+1)}\circ (id_{L^2(M)} \uset{M}{\otimes} \chi_{-2l}(y))(\xi)\\
= & w_{+(2l+1)}\left(\hat{x}_1 \uset{M}{\otimes} (w_{-2l}^* \varphi_{0,l}(y) w_{-2l})\left(\bar{\hat{x}}_2\uset{N}{\otimes} \cdots \uset{N}{\otimes} \hat{x}_{2l-1} \uset{M}{\otimes} \bar{\hat{x}}_{2l} \uset{N}{\otimes} \hat{x}_{2l+1} \right) \right)
= \sum_i \left(x_1 z_{1, i}^* v_1 z_{2, i} z_{3, i}^*v_2 \cdots v_l z_{2l, i}\right)^{\widehat{}} \\
= & x_1 \varphi_{0,l}(y) w_{-2l}\left(\bar{\hat{x}}_2\uset{N}{\otimes} \cdots \uset{N}{\otimes} \hat{x}_{2l-1} \uset{M}{\otimes} \bar{\hat{x}}_{2l} \uset{N}{\otimes} \hat{x}_{2l+1} \right)
= \varphi_{0,l}(y) w_{-2l}\left(({{x}}_2x_1^*)^{\bar{\widehat{} }} \uset{N}{\otimes} \cdots \uset{N}{\otimes} \hat{x}_{2l-1} \uset{M}{\otimes} \bar{\hat{x}}_{2l} \uset{N}{\otimes} \hat{x}_{2l+1} \right) \\
= & \varphi_{0,l}(y) w_{+(2l+1)} (\xi),
 \end{align*}
 where $ (w_{-2l}^* \varphi_{0,l}(y)
 w_{-2l})\left(\bar{\hat{x}}_2\uset{N}{\otimes} \cdots
 \uset{N}{\otimes} \hat{x}_{2l-1} \uset{M}{\otimes} \bar{\hat{x}}_{2l}
 \uset{N}{\otimes} \hat{x}_{2l+1} \right) = \sum_i \bar{\hat{z}}_{1,i}
 \uset{N}{\otimes} \hat{z}_{2,i} \uset{M}{\otimes} \cdots
 \uset{N}{\otimes} \hat{z}_{2l,i}\in {\mcal H}^o_{-2l} $ for some
 $z_{j,i} \in M$ and $\{i\}$ finite, $1 \leq j \leq 2l$; as was
 desired.

(3) First consider $k =2l$, $y \in P_{+2l} = N^\prime \cap M_{2l-1}$
 and $\xi = \hat{x}_1\uset{M}{\otimes}
 \bar{\hat{x}}_2\uset{N}{\otimes} \cdots \uset{M}{\otimes}
 \bar{\hat{x}}_{2l} \uset{N}{\otimes} \hat{x}_{2l+1} \in {\mcal
   H}_{+(2l+1)}^o$ with $x_i$'s as above, we have
\begin{align*}
 \chi_{+(2l+1)}(y)(\xi) & = \delta^{2l} {w_{+(2l+1)}^*} \left(E_{M_l}
 (y x_1x_2^* v_1 \cdots x_{2l}^*v_l x_{2l+1} e_{[0,l]})^{\widehat{}}
 \right) \\ & = {w_{+(2l+1)}^*} \circ \varphi_{0,l} (y)\circ u_{l-1}
 \left( w_{+2l} (\hat{x}_1\uset{M}{\otimes}
 \bar{\hat{x}}_2\uset{N}{\otimes} \cdots \uset{N}{\otimes} {x}_{2l-1}
 \uset{M}{\otimes} \bar{\hat{x}}_{2l} )
 \uset{N}{\otimes}\hat{x}_{2l+1} \right) \\ & = {w_{+(2l+1)}^*} \circ
 u_{l-1,0} \circ (\varphi_{-1,l-1}(y) \uset{N}{\ot}id_{L^2(M)}) \left(
 w_{+2l} (\hat{x}_1\uset{M}{\otimes} \bar{\hat{x}}_2\uset{N}{\otimes}
 \cdots \uset{N}{\otimes} {x}_{2l-1} \uset{M}{\otimes}
 \bar{\hat{x}}_{2l} ) \uset{N}{\otimes}\hat{x}_{2l+1} \right) \\ & =
     {w_{+(2l+1)}^*} \left( \varphi_{-1, l-1}(y) \circ w_{+2l}
     \left(\hat{x}_1\uset{M}{\otimes} \bar{\hat{x}}_2\uset{N}{\otimes}
     \cdots \uset{N}{\otimes} {x}_{2l-1} \uset{M}{\otimes}
     \bar{\hat{x}}_{2l} \right) v_l {x}_{2l+1} \right) \\
& = \left(Ad_{w_{+2l}^*} (\varphi_{-1, l-1}(y)) \left(\hat{x}_1\uset{M}{\otimes} \bar{\hat{x}}_2\uset{N}{\otimes} \cdots \uset{M}{\otimes} \bar{\hat{x}}_{2l} \right)\right) \uset{N}{\otimes} {\hat{x}}_{2l+1}
= \left(\chi_{+2l}(y) \uset{N}{\otimes} id_{{L^2(M)}}\right)(\xi)\\
& = \tilde{Q}_{RI_{+ 2l}} (\chi_{+2l}(y))(\xi).
\end{align*}

On the other hand, if $k$ is odd, say $2l+1$, $y \in P_{+(2l+1)} =
N^\prime \cap M_{2l}$ and $\xi = \hat{x}_1\uset{M}{\otimes}
\bar{\hat{x}}_2\uset{N}{\otimes} \cdots \uset{N}{\otimes}
\hat{x}_{2l+1} \uset{M}{\otimes} \bar{\hat{x}}_{2l+2} \in {\mcal
  H}_{+(2l+2)}^o$, then
\begin{align*}
& \tilde{Q}_{RI_{+ (2l+1)}} (\chi_{+(2l+1)}(y))(\xi)\\
= & \left(\chi_{+(2l+1)}(y) \uset{M}{\otimes} id_{\oline{\mcal H}}\right)(\xi)
= w_{+(2l+1)}^* \left( \varphi_{0, l}(y) \left( x_1 x^\ast_2 v_1 \cdots v_{l-1} x_{2l-1} x^\ast_{2l} v_l x_{2l+1} \right)^{\hat{}} \right) \uset{M}{\otimes} \bar{{\hat{x}}}_{2l+2}\\
= & w_{+(2l+2)}^* \left( \varphi_{0, l}(y) \left( \left( x_1 x^\ast_2 v_1 \cdots v_l x_{2l+1} \right)^{\hat{}} \right) x^\ast_{2l+2} \right)
= w_{+(2l+2)}^* \left( \varphi_{0, l}(y) \left( x_1 x^\ast_2 v_1 \cdots v_l x_{2l+1} x^\ast_{2l+2} \right)^{\hat{}} \, \right)\\
= & w_{+(2l+2)}^* \left( \varphi_{-1, l}(y) \left( x_1 x^\ast_2 v_1 \cdots v_l x_{2l+1} x^\ast_{2l+2} \right)^{\hat{}} \, \right)
= \chi_{+(2l+2)}(y)(\xi).
\end{align*}

(4) The initial case $k = 1 $ is trivial. Let $k = 2l$.  For $\xi
= \hat{x}_1\uset{M}{\otimes} \bar{\hat{x}}_2\uset{N}{\otimes}
\cdots \uset{M}{\otimes} \bar{\hat{x}}_{2l} \uset{N}{\otimes} 
{\hat{x}}_{2l+1} \in {\mcal H}_{+(2l+1)}^o$ with $x_i$'s in $M$, we have
\begin{align*}
\chi_{2l+1}(\delta e_{2l})(\xi ) & =  \delta^{2l+1} w_{2l+1}^* 
\left(E_{M_l}( e_{2l} x_1 x_2^*v_1 \cdots x_{2l}^* v_l x_{2l+1}e_{[0,l]}
)^{\widehat{}} \right)\\
&= \delta^{2l+1} w_{2l+1}^* \left(x_1x_2^*v_1 \cdots v_lx_{2l+1} E_{M_l} ( e_{2l}e_{[0,l]})^{\widehat{}} \right)
= \delta w_{2l+1}^* \left( (x_1x_2^*v_1 \cdots v_lx_{2l+1} e_2)^{\widehat{}} \right)\\
& = \delta w_{2l+1}^* \left( (x_1x_2^*v_l^*  \cdots v_1^*x_{2l+1} e_2)^{\widehat{}} \right)
=  w_{2l+1}^* \left( (x_1x_2^*v_l^*  \cdots x_{2l-2}^*v_2^*x_{2l-1}x^*_{2l}x_{2l+1} )^{\widehat{}} \right)
\end{align*}
where, in the third equality, we have used the fact that
$E_{M_l}(e_{2l}e_{[0,l]}) = \delta^{-2l}e_2$, which holds by the
uniqueness condition in Lemma \ref{uniqueness-lemma} applied to the
equality $e_{2l}e_{[0,l]} = e_2e_{[0,l]}$ (a routine verification
involving Temperley-Lieb relations satisfied by the Jones' projections
$\{e_l\}$).  On the other hand, by the definition of $\tilde{Q}$ and
Theorem \ref{bimod-pa-theorem} (i)(d), we have
\begin{align*}
& \tilde{Q}_{E_{+2l}}(\xi) = \delta^{-1} \sum_i \hat{x}_1\uset{M}{\otimes} \bar{\hat{x}}_2 \uset{N}{\otimes} \cdots \uset{N}{\otimes} \hat{x}_{2l-1} \uset{M}{\otimes} (b_i x_{2l+1}^* x_{2l})^{\bar{\widehat{}}} \uset{N}{\otimes} \hat{b}_i \\
= & \; \delta^{-1} w_{+(2l+1)}^*\left( (x_1x_2^* v_1 \cdots x_{2l-2}^* v_{l-1} x_{2l-1}x_{2l}^* x_{2l+1} b_i^* v_l b_i )^{\widehat{}} \right)
= w_{+(2l+1)}^*\left( (x_1x_2^* v_l^* \cdots x_{2l-2}^* v_{2}^* x_{2l-1}x_{2l}^* x_{2l+1})^{\widehat{}} \right).
\end{align*}

(5) Let $k=2l$, $y \in P_{+(2l+1)} = N^\prime \cap M_{2l}$ and $\xi
= \bar{\hat{x}}_1\uset{N}{\otimes} {\hat{x}}_2\uset{M}{\otimes} \cdots
\uset{N}{\otimes} {\hat{x}}_{2l} \in {\mcal H}_{-2l}^o$ where $x_i \in
M$ for all $ 1 \leq i \leq 2l$. Then, using the precription of the
action of tangles in the planar algebra $Q$ (described in Figure \ref{bicat-pa}, we get
\begin{align*}
& \tilde{Q}_{LE_{+(2l+1)}}( \chi_{+(2l+1)}(y) ) (\xi)
= \delta^{-1} \sum_i (e_{\mcal H} \uset{M}{\otimes} id_{{\mcal H}_{-2l}})\left( \bar{\hat{b}}_i \uset{N}{\otimes} \chi_{+(2l+1)}(y) \left( \hat{b}_i \uset{M}{\otimes} \bar{\hat{x}}_1 \uset{N}{\otimes} \cdots \uset{N}{\otimes} \hat{x}_{2l}\right) \right)\\
= & \delta^{2l-1} \sum_i (e_{\mcal H} \uset{M}{\otimes} id_{{\mcal H}_{-2l}})\left( \bar{\hat{b}}_i \uset{N}{\otimes} w_{+(2l+1)}^* \left( E_{M_l} ( b_i x_1^*v_l^* x_2 \cdots x_{2l-1}^*v_1^*x_{2l} e_{[0,l]})^{\hat{}} \, \right) \right)\\
= & \delta^{2l-1} \sum_i (e_{\mcal H} \uset{M}{\otimes} id_{{\mcal H}_{-2l}})\left( w_{-(2l+2)}^* \left( (b_i^*v_{l+1}^* z_l)^{\hat{}} \, \right) \right)\\
= & \delta^{2l-1} \sum_{i, i'} (e_{\mcal H} \uset{M}{\otimes} id_{{\mcal H}_{-2l}}) \left( w_{-(2l+2)}^* \left( ( b_i^* v_{l+1}^*b_{i'}^*v_l^* E_{M_{l-1}}(v_l b_{i'} z_l) )^{\hat{}} \, \right) \right)\\
= & \delta^{2l-1} \sum_{i, i'} (e_{\mcal H} \uset{M}{\otimes} id_{{\mcal H}_{-2l}}) \left( \bar{\hat{b}}_i \uset{N}{\otimes} ({b^*_{i'}})^{\hat{}} \uset{M}{\otimes} \bar{\hat{1}} \uset{N}{\otimes} w_{+(2l-1)}^*\left( (E_{M_{l-1}}(v_l b_{i'} z_l) )^{\hat{}} \, \right) \right)\\
= & \delta^{2l-1} \sum_{i, i'} b_i^*b_{i'}^* \bar{\hat{1}} \uset{N}{\otimes} w_{+(2l-1)}^*\left( (E_{M_{l-1}}(v_l b_{i'} z_l))^{\hat{}} \, \right)
= \delta^{2l-1} \sum_{i, i'} w_{-2l}^* \left( (b_i^*b_{i'}^* v_l^* E_{M_{l-1}}(v_l b_{i'} z_l) )^{\hat{}} \right)\\
= & \delta^{2l-1} \sum_{i} w_{-2l}^* \left( (b_i^* z_l)^{\hat{}} \, \right)
= \delta^{2l-1} \sum_{i} w_{-2l}^* \left( E_{M_l} (b_i^*xb_i x_1^*v_l^* x_2 \cdots x_{2l-1}^*v_1^*x_{2l} e_{[0,l]})^{\hat{}} \, \right)
= \delta^{-1}\sum_i \chi_{-2l}\left( b_i^* x b_i \right)(\xi),
\end{align*}
where $z_l:= E_{M_l} ( xb_i x_1^*v_l^* x_2 \cdots
x_{2l-1}^*v_1^*x_{2l} e_{[0,l]})$ and in the fourth and eighth
equalities we have used the fact that $\{ v_lb_i\}$ is a left Pimsner-Popa
basis for $M_{l-1} \subset M_l$ - see \cite{JS97}.
\end{pf}
\begin{rem}
Starting with an extremal subfactor, the associated planar algebra in
Theorem \ref{jones-theorem} is indeed the same as that in
\cite[Theorem 4.2.1]{Jon} since the map ${_N\mcal{L}}(L^2(M_k) ) =: N'
\ni x \longmapsto \delta^{-2} \sum_i b_i^* x b_i \in M' :=
    {_M\mcal{L}} (L^2(M_k))$ is the unique $tr_{N'}$ preserving
    conditional expectation from $N'$ onto $M'$. Such extension of
    \cite[Theorem 4.2.1]{Jon} was first established by Michael Burns
    in his thesis. However, the techniques used in \cite{Bur03} are
    different from our proof which is built up using graphical
    calculus of morphisms in a pivotal bicategory. Further, Jones and
    Penneys, in \cite{JP}, also obtain an extension of \cite[Theorem
      4.2.1]{Jon} in a slightly general set up.
\end{rem}
\begin{rem}\label{jones-theorem-rem}
Apart from the action of the tangles given in conditions (1) - (5), we
also mention below the action of few other useful tangles.

(a) $P_{RE_{+k}} = \delta \left. E^{M_{k-1}}_{M_{k-2}}
\right|_{P_{+k}}$, following from conditions (3) and (4).

(b) $P_{TR^r_{+k}} = \delta^k \left. tr_{M_{k-1}} \right|_{P_{+k}}$,
following from the action of right conditional expectation tangle in
(a).

(c) $\delta^{-k} P_{TR^l_{+{2l}}}$ (resp., $\delta^{-k}
P_{TR^l_{+(2l-1)}}$) is given by the trace on $P_{+2l} = N' \cap
M_{2l-1}$ (resp., $P_{+(2l-1)} = N' \cap M_{2l-2}$) induced by the
canonical trace on ${_N} {\mcal L} ( L^2 (M_{l-1}) )$ via the map
$\vphi_{-1,l-1}$ (resp., $\vphi_{0,l-1}$); this could be derived from
the precription of the action. It also turns out that this trace on
$P_{+2l}$ (resp., $P_{+(2l-1)}$) matches with the one induced by the
canonical trace on ${_N} {\mcal L} ( L^2 (M_{2l-1}) )$ (resp., ${_N}
{\mcal L} ( L^2 (M_{2l-2}) )$) via the usual inclusion.
\comments{(d) $\text{Range} (P_{LI^n_{(-)^n k}}) = M'_{n-1} \cap M_{n+k-1}$
where $LI^n_{(-)^n k} : (-)^n k \rightarrow +(n+k)$ is the tangle
obtained from $LI_{-k}$ by replacing the single string not connected
to the internal disc, by $n$ many parallel strings and changing the
sign of the internal disc from $-$ to $(-)^n$.
}
\end{rem}
\begin{cor}\label{dualpa}
If $P$ is the planar algebra  associated to the tower $\{M_k\}_{k \geq -1}$ with Jones projections $\{e_k\}_{k \geq 1}$ (as in Theorem \ref{jones-theorem}), then

(a) $P_{E'_{-k}} (y) = \delta \us{i}{\sum} b^*_i e_1 y e_1 b_i$ for all $y \in P_{-k} = M' \cap M_k$, where $E'_{-k} = LI_{+(k-1)} \circ LE_{-k}$ and $\{ b_i \}_i$ is a left basis for $N \subset M$,

(b) $\lambda_n (P) =$ the planar algebra associated to the tower $\{ M_{k+n} \}_{k \geq -1}$ with Jones projections $\{ e_{k+n} \}_{k \geq 1}$.
\end{cor}
\begin{pf}
(a) We follow the same notations as in Theorem \ref{jones-theorem} and prove this only for the case $k = 2l \in 2\N$ because the case for odd $k$'s can be proved using similar arguments.
Let $a_{\mcal H}$ (resp., $e_{\oline{\mcal H}}$) denote the isomorphism ${\mcal H} \ni \alpha \os{a_{\mcal H}}{\longmapsto} \oline{ \oline{\alpha} } \in \oline{ \oline{\mcal H} }$ (resp., the usual evaluation map from $\oline{ \oline{\mcal H} } \us{M}{\otimes} \oline{\mcal H}$ to $L^2 (N)$), and $\xi = \bar{\hat{x}}_1\uset{N}{\otimes} {\hat{x}}_2\uset{M}{\otimes} \cdots \uset{N}{\otimes} {\hat{x}}_{2l} \in {\mcal H}_{-2l}^o$ where $x_i \in M$ for all $ 1 \leq i \leq 2l$.
Then, using the prescription of the actions of tangles, we get
\begin{align*}
& \tilde{Q}_{E'_{-2l}} \left( \chi_{- 2l}(y) \right) (\xi)\\
= & \delta \left( id_{\oline{\mcal H} } \us{N}{\otimes} \left[ e_{\oline{\mcal H} } \circ (a_{\mcal H} \us{M}{\otimes} id_{\oline{\mcal H} }) \right] \us{N}{\otimes} id_{{\mcal H}_{+(2l-1)}} \right) \left(\bar{\hat{x}}_1 \us{N}{\otimes} \hat{1}_M \us{M}{\otimes} \left[ w^*_{-2l} \circ \vphi_{0,l} (y) \circ w_{-2l} \left(\oline{\hat{1}}_M \uset{N}{\otimes} {\hat{x}}_2\uset{M}{\otimes} \cdots \uset{N}{\otimes} {\hat{x}}_{2l} \right) \right] \right)\\
= & \delta \us{i}{\sum} \left( id_{\oline{\mcal H} } \us{N}{\otimes} \left[ e_{\oline{\mcal H} } \circ (a_{\mcal H} \us{M}{\otimes} id_{\oline{\mcal H} }) \right] \us{N}{\otimes} id_{{\mcal H}_{+(2l-1)}} \right) \left( \oline{\hat{x}}_1 \us{N}{\otimes} \hat{1}_M \us{M}{\otimes} w^*_{-2l}  (b^*_i v^*_l E_{M_{l-1} } (v_l b_i z) )^{\widehat{}} \right)\\
&  \text{(where } \hat{z} := \vphi_{0,l} (y) \circ w_{-2l} \left( \bar{\hat{1}}_M \uset{N}{\otimes} {\hat{x}}_2\uset{M}{\otimes} \cdots \uset{N}{\otimes} {\hat{x}}_{2l} \right) \text{ and since } \{ v_l b_i \}_i \text{ forms a left basis for } M_{l-1} \subset M_l \text{)}\\
= & \delta \us{i}{\sum} \left( id_{\oline{\mcal H} } \us{N}{\otimes} \left[ e_{\oline{\mcal H} } \circ (a_{\mcal H} \us{M}{\otimes} id_{\oline{\mcal H} }) \right] \us{N}{\otimes} id_{{\mcal H}_{+(2l-1)}} \right) \left( \oline{\hat{x}}_1 \us{N}{\otimes} \hat{1}_M \us{M}{\otimes} \oline{\hat{b}}_i \us{N}{\otimes} w^*_{+(2l-1)}  ( E_{M_{l-1} } (v_l b_i z) )^{\widehat{}} \right)\\
= & \delta \us{i}{\sum} \left( \oline{\hat{x}}_1 \us{N}{\otimes} E_N (b^*_i) w^*_{+(2l-1)}  ( E_{M_{l-1} } (v_l b_i z) )^{\widehat{}} \right) = \delta \oline{\hat{x}}_1 \us{N}{\otimes}w^*_{+(2l-1)}  ( E_{M_{l-1} } (v_l z) )^{\widehat{}}
=  \delta w^*_{-2l} ( x^*_1 v^*_l  E_{M_{l-1} } (v_l z) )^{\widehat{}}\\
= & \delta^3 w^*_{-2l} ( x^*_1 v^*_{l-1} e_l  E_{M_{l-1} } (e_l v_{l-1} z) )^{\widehat{}}
= \delta w^*_{-2l} ( x^*_1 v^*_{l-1} e_l v_{l-1} z )^{\widehat{}}
= w^*_{-2l} \left( x^*_1 v_1 \vphi_{0,l} (y) (v_1 x_2 x^*_3 v_2 \cdots v_l x_{2l}) \right)^{\widehat{}}
\end{align*}
On the other hand,
\begin{align*}
\delta \us{i}{\sum} \chi_{-2l} ( b^*_i e_1 y e_1 b_i ) (\xi)
& = \delta^{2l+1} \us{i}{\sum} w^*_{-2l} \left( E_{M_l} ( b^*_i e_1 y e_1 b_i x^*_1 v_1 x_2 x^*_3 v_2 \cdots v_l x_{2l} e_{[0,l]}) \right)^{\widehat{}}\\
& = \delta^{2l+1} \us{i}{\sum} w^*_{-2l} \left( b^*_i E_N (b_i x^*_1) e_1 E_{M_l} ( y v_1 x_2 x^*_3 v_2 \cdots v_l x_{2l} e_{[0,l]}) \right)^{\widehat{}}\\
& = w^*_{-2l} \left( x^*_1 v_1 \vphi_{0,l} (y) (v_1 x_2 x^*_3 v_2 \cdots v_l x_{2l}) \right)^{\widehat{}}.
\end{align*}

(b) Since $\lambda_n$ is additive with respect to $n$, it is enough to prove the statement for $n=1$. Note that $\lambda_1 (I_{\vlon k}) = E'_{\vlon (k+1)}$ for all colors $\vlon k$.
By part (a) (resp., conditions (2) and (5) of Theorem \ref{jones-theorem}), we get $\lambda_1 (P)_{-k} = M'_1 \cap M_{k+1}$ (resp., $\lambda_1 (P)_{+k} = M' \cap M_{k}$). So, by the uniqueness part of Theorem \ref{jones-theorem}, it remains to check whether the conditions (1) to (5) therein hold in this case; these verifications are completely straight forward and we will skip them.
\end{pf}

\subsection{Reconstruction of bimodule}\label{bimodrecon} 
Starting with a bimodule planar algebra, we will construct a bifinite
bimodule whose associated bimodule planar algebra is isomorphic to the
one which we started with. We extensively use the techniques of
constructing an extremal subfactor from a spherical unimodular
bimodule planar algebra in \cite{JSW08} and \cite{KS08}; in fact,
we first show that their construction with necessary modifications,
works without the assumption of sphericality.

Let $P$ be a unimodular bimodule planar algebra with modulus
$(\delta,\delta)$ such that $\delta >1$. We will work only with the
positive part of $P$; so, for the time being, we will write $P_k$ in
place of $P_{+k}$. Henceforth, we will exactly follow the set-up of
\cite{JSW08} and \cite{KS08} and not mention this fact at every step;
whenever some modifications become necessary, we will explicitly
mention them.  Set $P^k = \uset{l\in \N_0}{\oplus} P^k_l$ where $P^k_l
:= P_{k+l}$ for $k, l \in \N_0$. Define the following structures on $P^k$:
\begin{enumerate}
\item {\em Multiplication:}\\ \psfrag{lhs}{$P^k \times P^k \supset
  P^k_l \times P^k_m \ni (x,y) \overset{\cdot}{\mapsto} x \cdot y :=
  \oset{2(l \wedge m)}{\uset{i=0}{\sum}} P$} \psfrag{rhs}{$\in
  \overset{2(l \wedge m)}{\uset{i=0}{\oplus}} P^k_{l+m-i} \subset
  P^k$} \psfrag{k}{$k$} \psfrag{+}{$+$} \psfrag{i}{$i$}
  \psfrag{x}{$x$} \psfrag{y}{$y$} \psfrag{2l-i}{$2l-i$}
  \psfrag{2m-i}{$2m-i$}
  \includegraphics[scale=0.2]{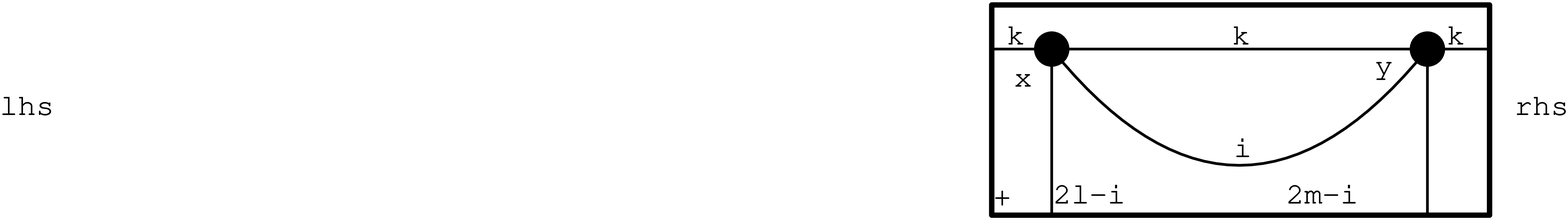}\\ We will
  also denote the element in the sum when $i=0$, by $x \odot y$; the
  $i$-th element will be denoted by $x \uset{i}{\odot} y$.
\item {\em ${\ast}$-structure:}\\
\psfrag{lhs}{$P^k \supset P^k_l \ni x \overset{\dagger}{\mapsto}  
x^\dagger := P$}
\psfrag{rhs}{$\in P^k_l \subset P^k$}
\psfrag{k}{$k$}
\psfrag{+}{$+$}
\psfrag{2l}{$2l$}
\psfrag{x*}{$x^\ast$}
\includegraphics[scale=0.15]{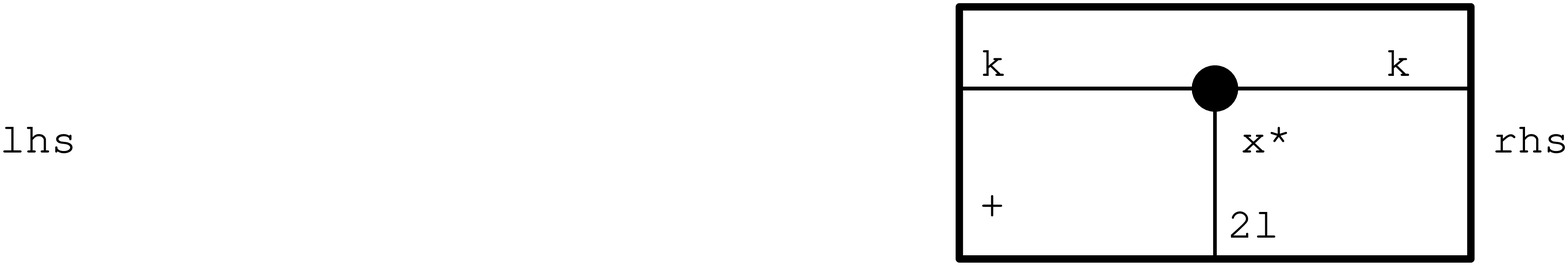}
\item {\em Trace:}
\[ 
P^k \supset P^k_l \ni x \overset{t_k}{\mapsto} t_k(x) := \delta_{l=0}
\delta^{-k} P_{TR^l_k} (x) \in \C
\]
\end{enumerate}
$P^k$ becomes an associative $*$-algebra with respect to (1) and (2).
\begin{rem}
The trace-functional on $P^k$ which was defined in \cite{JSW08,KS08},
does not satisfy tracial property in the absence of sphericality. So,
we had to take this specific spherical isotopy of the trace in
\cite{JSW08,KS08} so that the functional defined in (3), is indeed a
trace on $P^k$.
\end{rem}
To see the tracial property, note that for $x \in P^k_l$, $y \in
P^k_m$,
\[
\psfrag{lhs}{$
t_k(x \cdot y)=\delta_{l = m} \delta^{-k} P$}
\psfrag{rhs}{=$t_k(y \cdot x)$.}
\psfrag{k}{$k$}
\psfrag{x}{$x$}
\psfrag{y}{$y$}
\psfrag{2l}{$2l$}
\includegraphics[scale=0.1]{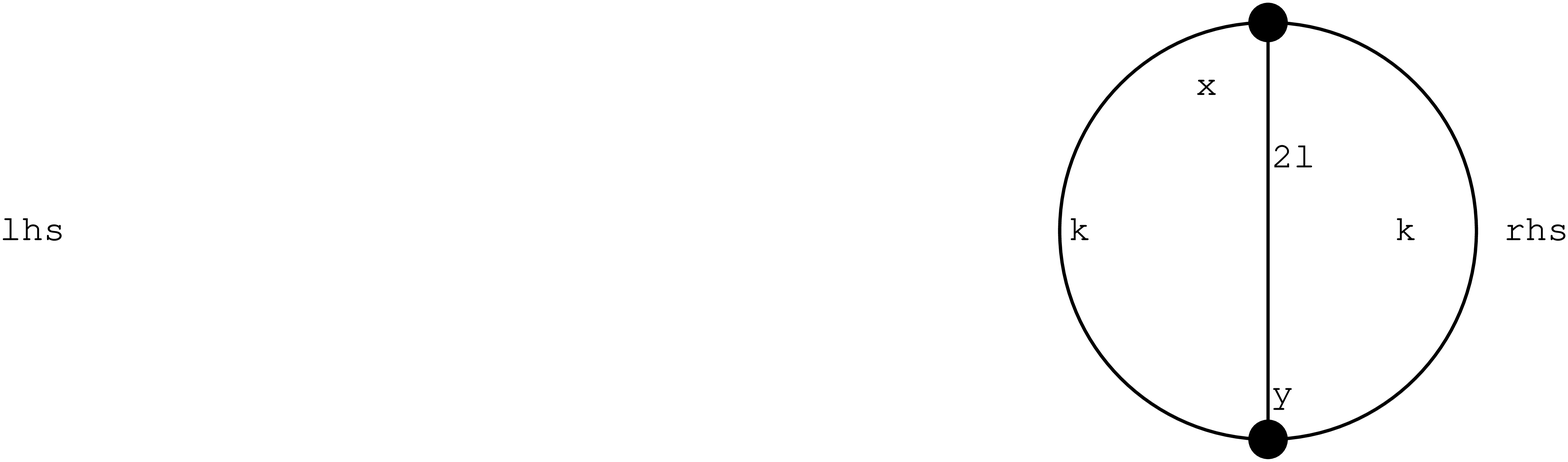}
\]
We define a sesquilinear form $\langle \cdot, \cdot \rangle$ on $P^k$
by $\langle x, y \rangle = t_k (x^\dagger \cdot y)$ for $x,y\in
P^k$. It is straight forward to check that $(i)$ $t_k$ is positive
definite and hence the sesquilinear form gives an inner product,
$(ii)$ $\{P^k_l\}_{l\in \N_0}$ are mutually orthogonal subspaces of
$P^k$ with respect to this inner product. Let $L: P^k \rightarrow
End_{\C}(P^k)$ and $R: P^k \rightarrow End_{\C}(P^k)$ denote the left
and right multiplication operators respectively. The proof of the
boundedness of these operators in \cite{JSW08,KS08} needs a little
modification in the non-spherical case because the inner product is
slightly different; however, the main idea of the proof will remain
the same.
 \begin{lem}
 $L_a$ and $R_a$ are bounded for $a \in P^k$.
 \end{lem}
 \begin{pf}
 Without loss of generality, let $a \in P^k_l$ for $l\in \N_0$. Now,
 $L_a = \oset{2l}{\uset{i=0}{\sum}} L^i_a$ where $L^i_a : P^k
 \rightarrow P^k$ is defined by
 \[
P^k \supset P^k_m \ni x \oset{L^i_a}{\longmapsto} \delta_{i\leq 2m}
\left( a \uset{i}{\odot} x\right) \in P^k_{l+m-i} \subset P^k.
 \]
 Since $\{P^k_l\}_{l\in \N_0}$ are mutually orthogonal, therefore it
 is enough to show that there exists an $M>0$ such that $L^i_a (x)
 \leq M \|x\|$ for all $x \in P^k_m$, $m \in \N_0$. Without loss of
 generality, let $m > 2l$. Note that\\
\psfrag{lhs}{$\| L^i_a (x) \|^2 \! = \! {\delta^{-k}} \! P$}
\psfrag{rhs}{ where $y \in \! P_{(-)^i(k+i)}$ is the positive square root of}
\psfrag{k}{$k$}
\psfrag{k+i}{$k+i$}
\psfrag{a}{$a$}
\psfrag{a*}{$a^\ast$}
\psfrag{x}{$x$}
\psfrag{y}{$y$}
\psfrag{k+i}{$k+i$}
\psfrag{2l-i}{$2l-i$}
\psfrag{2m-i}{$2m-i$}
\includegraphics[scale=0.2]{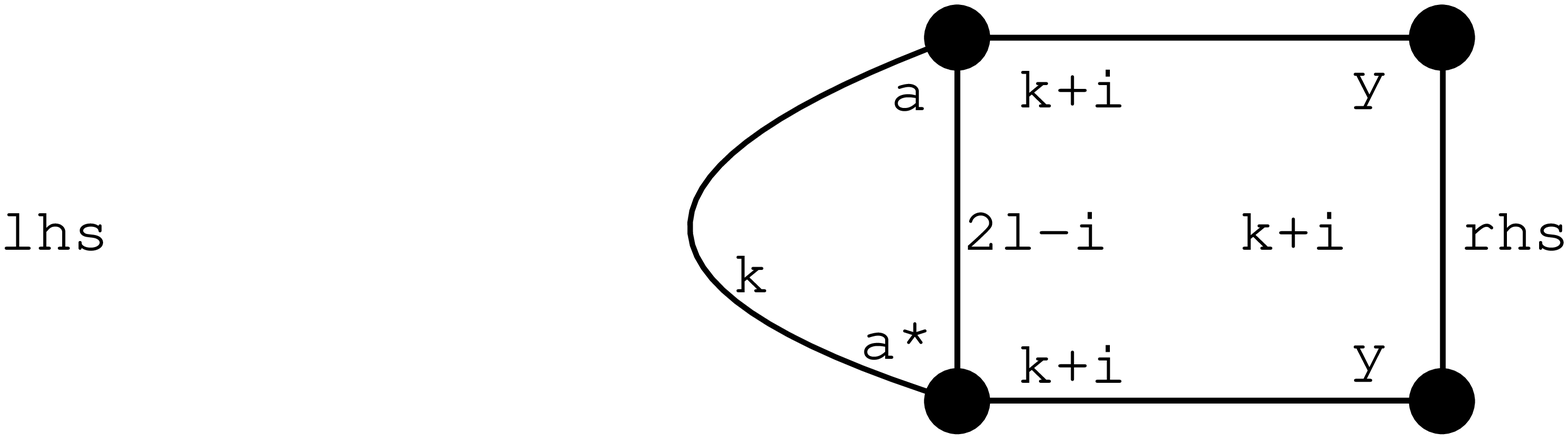}
\psfrag{space}{} \psfrag{lhs}{$P$} \psfrag{rhs}{.}
\psfrag{k+2m-i}{$k+2m-i$} \psfrag{k}{$k$} \psfrag{i}{$i$}
\psfrag{x*}{$x^\ast$} \psfrag{x}{$x$} \psfrag{*}{$(-)^i$}
\includegraphics[scale=0.2]{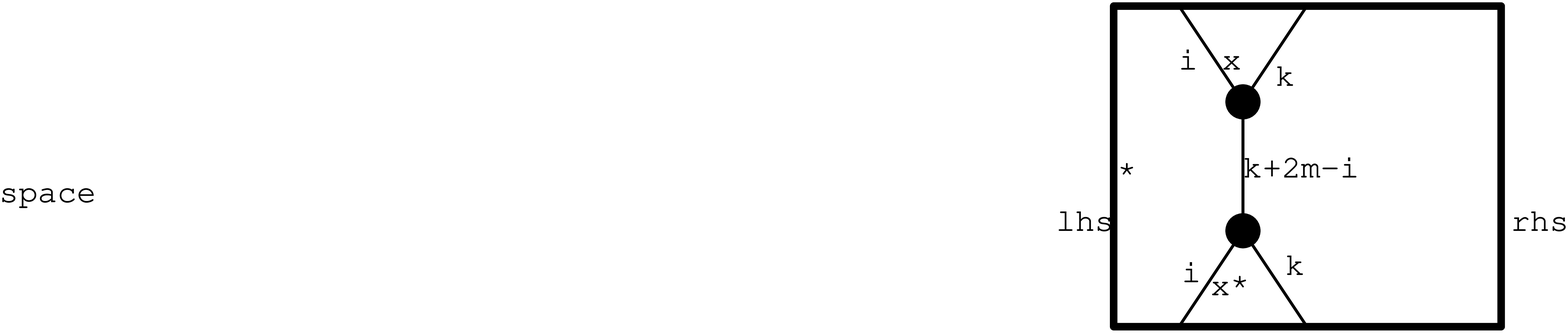}\\
 To see
 that this is indeed a positive element, one needs to take an
 appropriate rotation of $x$ and use positivity of the action of the
 right conditional expectation tangle from $(-)^i (k+m)$ to $(-)^i
 (k+i)$. Thus, $\| L^i_a (x) \|^2$ can also be expressed as $\|
 P_{T_a} (y) \|^2$ where $T_a$ is the semi-labelled tangle
 \psfrag{k+2m-i}{$k+2m-i$}
 \psfrag{k}{$k$}
 \psfrag{k+i}{$k+i$}
 \psfrag{a}{$a$}
 \psfrag{2l-i}{$2l-i$}
 \psfrag{*}{$(-)^i$}
 \psfrag{+}{$+$}
 \includegraphics[scale=0.2]{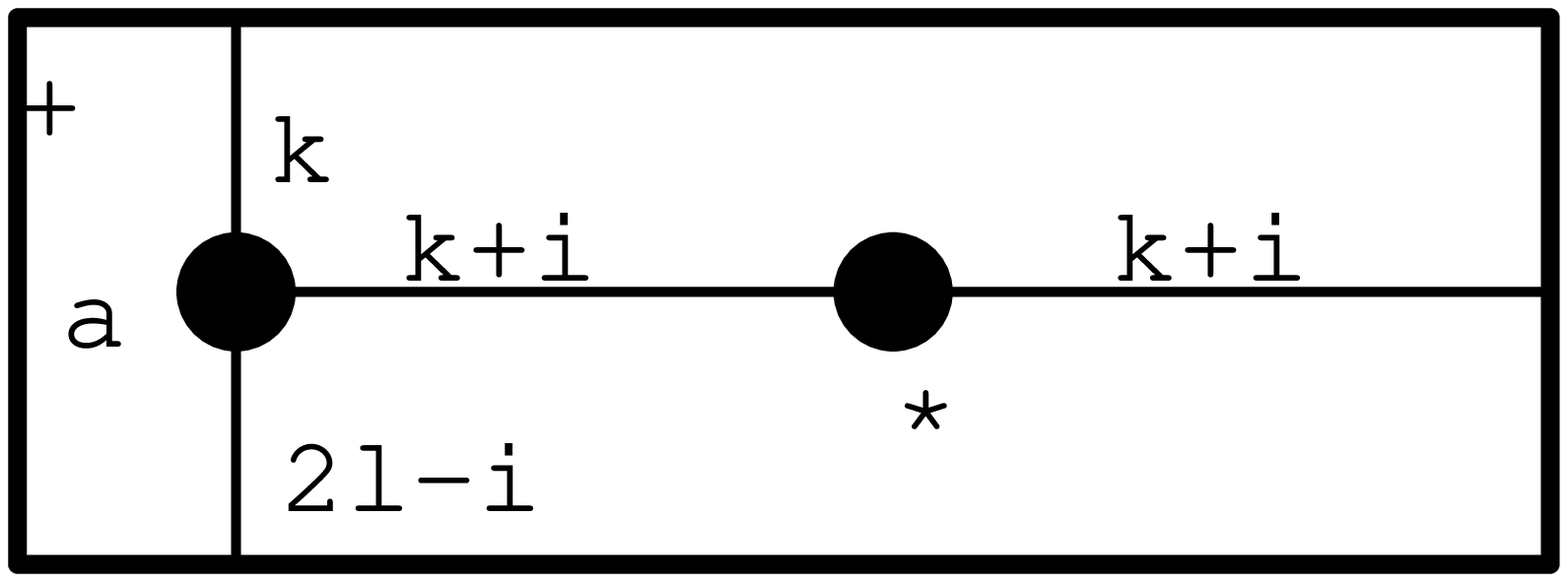}.
 We equip
 the domain $P_{(-)^i (k+i)}$ of $P_{T_a}$, with the norm $\| \cdot
 \|'$ coming from $\delta^{-k} P_{TR^l_{(-)^i (k+i)}}$ so that $\| y
 \|' = \|x\|$. Hence, the desired inequality is satisfied by setting
 $M :=$ the operator norm of $P_{T_a}$ which is independent of $m >
 2m$ and $x \in P^k_m$. Boundedness of $R_a$ will follow using the
 same kind of arguments.
 \end{pf}
 
 Let ${\mcal H}_k$ be the completion of $P^k$ and $M_k := (LP^k)''
 \subset {\mcal L}({\mcal H}_k)$. Note that $\hat{1}_{P^k} \in {\mcal
   H}_k$ is a cyclic, separating, trace vector for $M_k$. The unital
 $\ast$-algebra inclusion $\! \vcenter{ \psfrag{+}{$+$}
   \psfrag{k}{$k$} \psfrag{x}{$x$} \psfrag{2l}{$2l$} \psfrag{lhs}{$P^k
     \supset P^k_l \ni x \mapsto P$} \psfrag{rhs}{$\in P^{k+1}_l
     \subset P^{k+1} \text{ induces an}$}
   \includegraphics[scale=0.2]{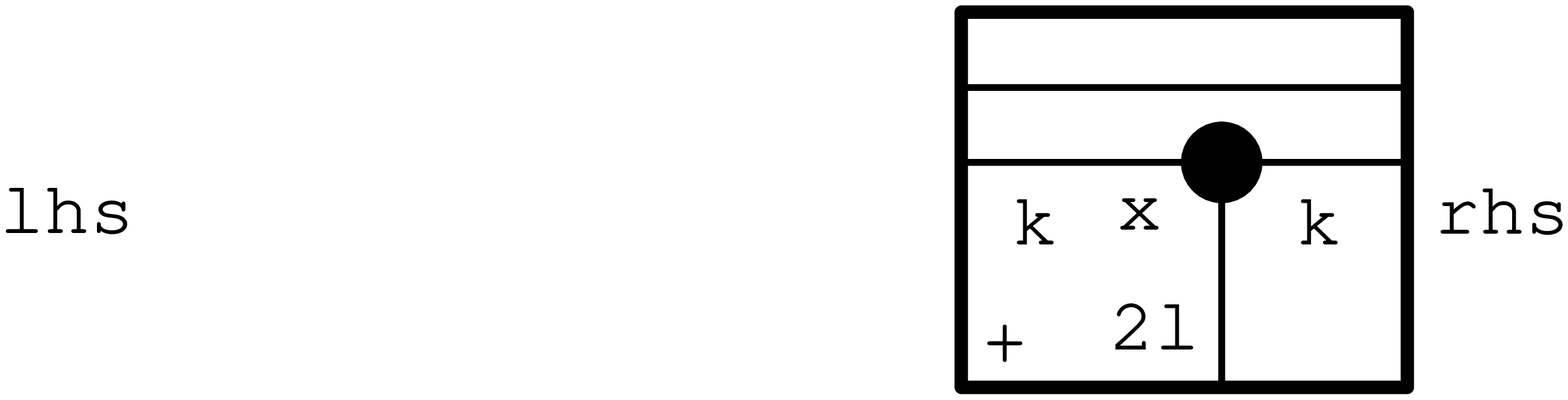} }$ an
 inclusion $M_k \hookrightarrow M_{k+1}$. Note that $LP^k_0 \subset
 M'_0 \cap M_k$. One can also prove that $P_k \ni x \mapsto L_x \in
 M'_0 \cap M_k$ is a $C^\ast$-algebra isomorphism and $M_k$'s are
 $II_1$-factors; for proofs, see \cite{JSW08,KS08} where $\delta > 1$
 is crucially used.
\begin{rem}\label{trace problem}
The above mentioned isomorphism takes the normalized left picture
trace on $P_k$ to the $M_k$-trace on $M'_0 \cap M_k$. So, even if we
assume that $\{M_k\}_{k \geq 2}$ is a tower of basic constructions of
$M_0 \subset M_1$, there is no hope of the unimodular bimodule planar
algebra associated to $M_0 \subset M_1$ being isomorphic to $P$
(unless $P$ is spherical) because of the condition in Remark
\ref{jones-theorem-rem}(c). We will fix this issue in the following
theorem using a certain perturbation.
\end{rem}

 If $Q$ is a unimodular bimodule planar algebra, then there exists a
 unique invertible positive central element $z \in Q_{+1}$ such that
 $Q_{TR^l_{+1}} (\cdot) = Q_{TR^r_{+1}} (\cdot \; z)$. It is easy to
 check that $z$ is a weight of $Q$ and $Q_{TR^l_{\vlon k}} (\cdot) =
 Q_{TR^r_{\vlon k}} (\cdot \; z_{\vlon k})$. We will refer this as the
 {\em trace intertwiner weight of $Q$}.

 \begin{thm}\label{sfrecon}
 If $Q$ is a unimodular bimodule planar algebra with modulus
 $(\delta,\delta)$, then there exists a finite
 index subfactor $M_0 \subset M_1$ of type $II_1$ whose associated
 planar algebra is isomorphic to $Q$.
\end{thm}
\begin{pf}
Without loss of generality, we may assume $\delta > 1$ because index
of a bimodule planar algebra is at least $1$ (which follows from
positivity of the action of the trace tangles) and the index $1$ case
is a triviality.  Let $P$ be the perturbation of $Q$ by the
decomposition $z = z^{1/2}\cdot z^{1/2}$ of the trace intertwiner
weight $z$ of $Q$. Clearly, $P$ is a unimodular bimodule planar
algebra. Consider the tower of $II_1$-factors $\{M_k\}_{k\geq 0}$
constructed from $P$ right before Remark \ref{trace problem} and the
isomorphism $Q_k = P_k \ni x \oset{\phi}{\mapsto} L_x \in M'_0 \cap
M_k$. Note that the map $Q_{RI_k} = P_{RI_k}$ under $\phi$, is given
by the inclusion $M_k \subset M_{k+1}$. So, by Theorem
\ref{jones-theorem}, it remains to show that $(i)$ $\{M_k\}_{k \geq
  2}$ is a tower of basic construction of $M_0 \subset M_1$ with Jones
projections $\{e_k := \delta^{-1} \phi(Q_{E_{+k}})\}_{k \geq 1}$ and
$(ii)$ $\phi \circ Q_{E'_{+k}} = \delta^{-1} \uset{b \in B}{\sum} b
\phi (\cdot) b^\ast$ for all $k \geq 1$, where $B$ is a right
Pimsner-Popa basis for $M_0 \subset M_1$.
\vspace*{1mm}

\noindent {\em Proof of (i)}: First note that the unique trace
preserving conditional expectation from $M_k$ to $M_{k-1}$ is induced
(via $L$) by the map $P^k \supset P^k_l \ni x
\oset{E^k_{k-1}}{\longmapsto} \left\{ \vcenter{
  \psfrag{lhs}{$\delta^{-1} P$} \psfrag{rhs1}{$\in P^{k-1}_l \subset
    P^{k-1}$, if $k$ is even,} \psfrag{rhs2}{$\in P^{k-1}_l \subset
    P^{k-1}$, if $k$ is odd.}  \psfrag{k-1}{$k-1$} \psfrag{+}{$+$}
  \psfrag{2l}{$2l$} \psfrag{x}{$x$} \psfrag{z}{$z$}
  \psfrag{zinv}{$z^{-1}$}
  \includegraphics[scale=0.2]{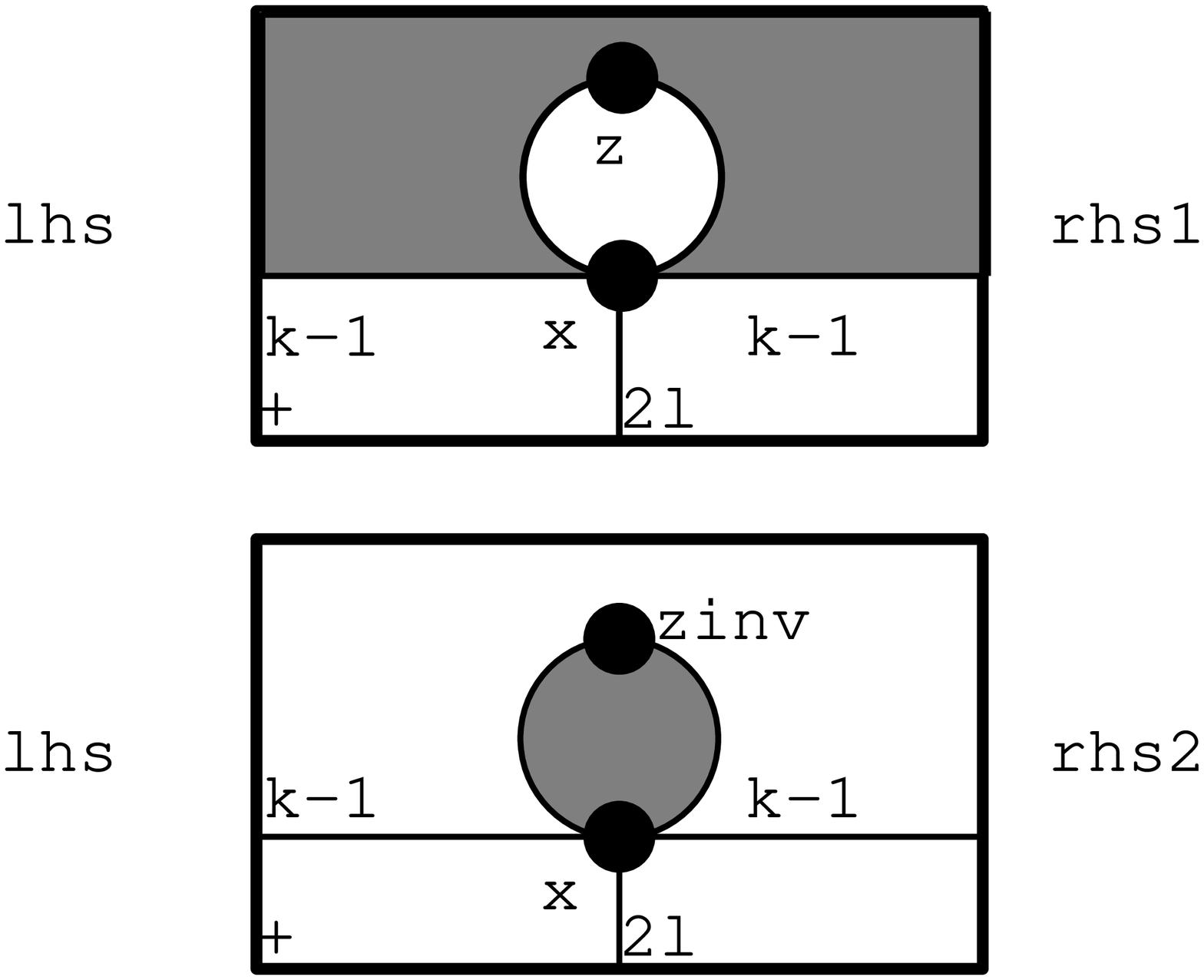}} \right.$
Using the definition of $P=Q^{(z^{1/2},z^{1/2})}$, it immediately
follows that $Q_{E_{+k}} x Q_{E_{+k}} = \delta E^k_{k-1} (x)
Q_{E_{+k}}$ for all $x \in P^k$. Thus, for all $y \in M_k$, $e_k y e_k
= E^{M_k}_{M_{k-1}} (y) e_k$. Moreover, it is straight forward to
check that for all $x \in P^{k+1}$, $Q_{E_{+k}} x = \delta Q_{E_{+k}}
E^k_{k-1}(Q_{E_{+k}} x)$ which implies $e_k M_{k+1} = e_k M_k$. Thus,
by Proposition \ref{bc-fact}, $M_{k+1}$ is a basic construction of
$M_{k-1} \subset M_k$ with Jones projection $e_k$ and $[M_k : M_{k-1}]
= \delta^2$.
\vspace*{1mm}

\noindent {\em Proof of (ii)}: Since $L^2 (M_1) \cong {\mcal H}_1
\cong \uset{k \in \N_0}{\oplus} P^1_k$, therefore for each $b \in B$
and $k \in \N_0$, there exists $b_k \in P^1_k$ such that $\hat{b} =
\uset{k \in \N_0}{\sum} \hat{L}_{b_k}$. Also, $\hat{b^*} = \uset{k \in
  \N_0}{\sum} \hat{L}_{b^{\dagger}_k}$. Since $B$ is a basis for $M_0
\subset M_1$, we have
\[
\hat{L}_{1_{P^2}} = \hat{1}_{M_2} = \uset{b\in B}{\sum} (be_1
b^\ast)\hat{} = \delta^{-1}\uset{b \in B, k,l \in \N_0}{\sum}
 (L_{b_k \cdot Q_{E_{+1}} \cdot
  b^{\dagger}_l })\hat{} = \delta^{-1}\uset{m\in \N_0}{\sum}\ 
\uset{l+k=m+i,\; i \leq 2(k\wedge l) }{\uset{ b \in B,\; i,k,l \in
    \N_0 \text{ s.t.}  }{\sum}} (L_{\left(b_k \odot Q_{E_{+1}}\right) \uset{i}{\odot}
  b^{\dagger}_l })\hat{}.
\]
\comments{({\bf Pending: Right boundedness used in 2nd equality})}
Using the
relation $Q_{E_{+1}}=P_{E^{(z^{-1/2},z^{-1/2})}_{+1}}$, the above
gives the following formula in terms of pictures:
\[
\psfrag{left}{$\uset{l+k=m+i,\; i \leq 2(k\wedge l)
}{\uset{
b \in B,\; i,k,l \in \N_0 \text{ s.t.}
}{\sum}} P$}
\psfrag{middle}{$=\delta \; \delta_{m = 0} P$}
\psfrag{right}{}
\psfrag{+}{$+$}
\psfrag{i}{$i$}
\psfrag{z}{$z$}
\psfrag{bk}{$b_k$}
\psfrag{bl*}{$b^\ast_l$}
\psfrag{2k-i}{$2k-i$}
\psfrag{2l-i}{$2l-i$}
\includegraphics[scale=0.2]{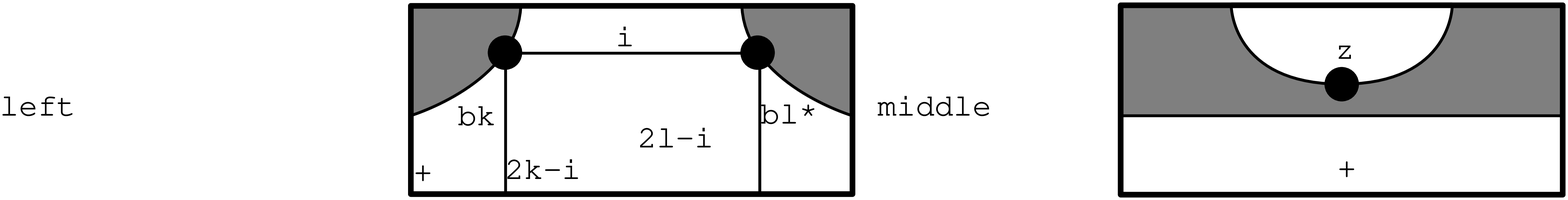}
\]
for all $m \in \N_0$ where the sum in the left hand side is with
respect to the Hilbert space norm in $P^2$. Now, let $L_x \in M'_0
\cap M_k = LP^k_0$ for $x \in P^k_0 = Q^k_0$. So, $\uset{b \in
  B}{\sum} (b L_x b^\ast)\hat{} = \uset{m\in \N_0}{\sum}\,
\uset{l+k=m+i,\; i \leq 2(k\wedge l) }{\uset{ b \in B,\; i,k,l \in
    \N_0 \text{ s.t.}  }{\sum}} (L_{\left( b_k \odot x \right)
  \uset{i}{\odot} b^{\dagger}_l })\hat{}$. Note that
\[
\psfrag{left}{$\left(b_k \odot x\right) \uset{i}{\odot} b^{\dagger}_l = P$}
\psfrag{middle}{.}
\psfrag{+}{$+$}
\psfrag{i}{$i$}
\psfrag{x}{$x$}
\psfrag{z}{$z$}
\psfrag{bk}{$b_k$}
\psfrag{bl*}{$b^\ast_l$}
\psfrag{k-1}{$k-1$}
\psfrag{2k-i}{$2k-i$}
\psfrag{2l-i}{$2l-i$}
\includegraphics[scale=0.2]{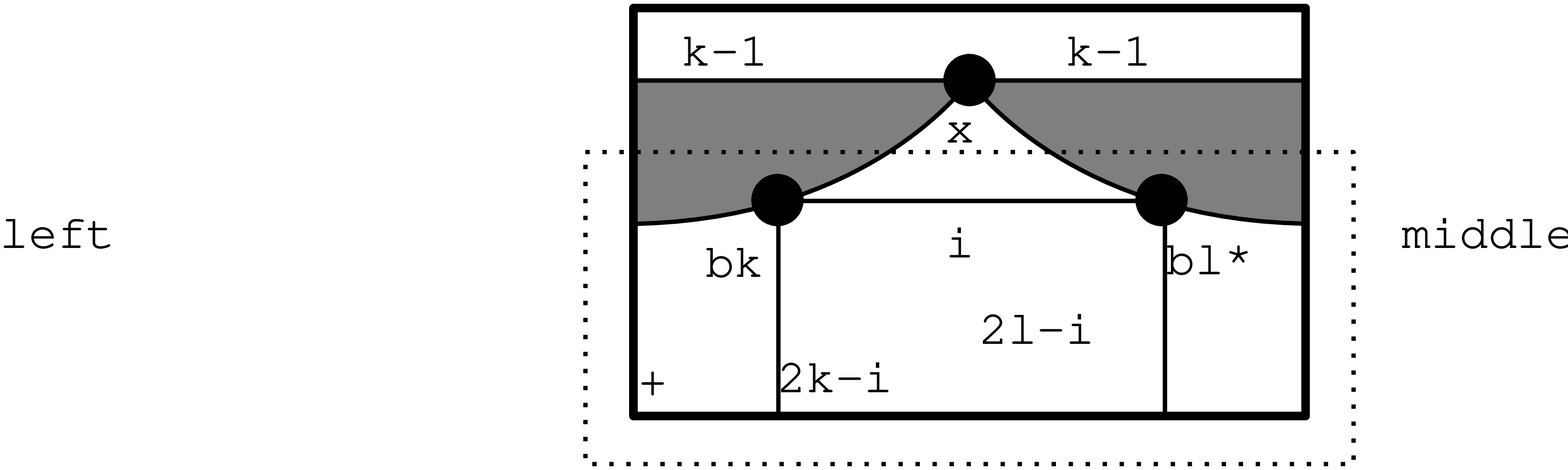}
\]
Using the isomrphism ${\mcal H}_k \supset (P^k)\hat{} \ni \hat{x}
\mapsto \hat{L}_x \in (M_k)\hat{} \subset L^2(M_k)$ and replacing the
sum (in the expansion of $\uset{b \in B}{\sum} (b L_x b^\ast)\hat{}$)
of the stuff inside the dotted box by the right side of the above
formula, we obtain
\[
\psfrag{left}{$\uset{b \in B}{\sum} b \phi(x) b^\ast = \uset{b \in
    B}{\sum} b L_x b^\ast = \delta \; L_P$
}
\psfrag{middle}{$= \delta \; L_{Q_{E'_{+k}} (x)} =\delta \; \phi (Q_{E'_{+k}} (x))$.}
\psfrag{+}{$+$}
\psfrag{i}{$i$}
\psfrag{x}{$x$}
\psfrag{z}{$z$}
\psfrag{bk}{$b_k$}
\psfrag{bl*}{$b^\ast_l$}
\psfrag{k-1}{$k-1$}
\psfrag{2k-i}{$2k-i$}
\psfrag{2l-i}{$2l-i$}
\includegraphics[scale=0.2]{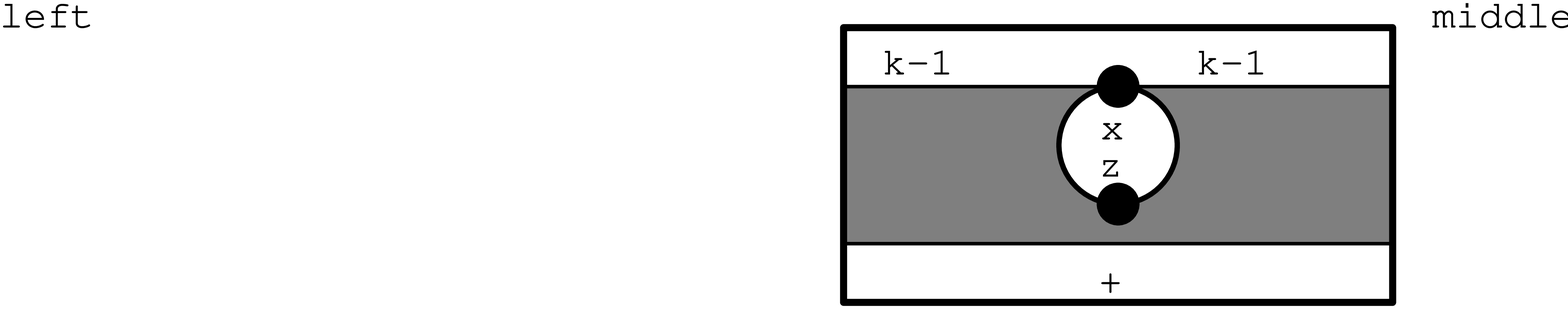}
\]\end{pf}
\begin{rem}
We must mention here that Popa (in \cite{Pop02}) started with his set
up of `generalized $\lambda$-lattice' (which, we believe, is analogous
to a unimodular bimodule planar algebra in our case) and gave a
general method of constructing a (non-extremal) subfactor whose standard
invariant corresponds to the generalized $\lambda$-lattice.
\end{rem}
Using Theorem \ref{sfrecon}, we now construct a bimodule from a
bimodule planar algebra in the next theorem.
\begin{thm}\label{bmrecon}
Let $P$ be a bimodule planar algebra with modulus $(\delta_- ,
\delta_+)$ such that $(\delta_+ \delta_-) > 1$. Then there exists a
bifinite bimodule whose associated bimodule planar algebra is
isomorphic to $P$.
\end{thm}
\begin{pf}
Let $Q$ be the normalization of $P$ and $\delta := \sqrt{\delta_+
  \delta_-}$.  By Theorem \ref{sfrecon}, there exist a subfactor $N
\subset M$ of type $II_1$, whose associated planar algebra is
isomorphic to $Q$.  Without loss of generality, let $\delta_+ >1$.
Choose $p \in {\mscr P}(M)$ such that $tr_M(p)=\delta^{-1}_+$. Set
$p':= JpJ \in M' := {_M}{\mcal L} (L^2 (M))$, ${\mcal H} := L^2 (M) p
= Range (p')$, $A:= N p'$ and $B:= M_{p}$ where $J$ is the canonical
anti-unitary involution on $L^2 (M)$.  Consider the bimodule $_A
{\mcal H}_B$. Note that (i) $dim({_A}{\mcal H}) = tr_{M'}(p')
dim({_N}L^2(M))= \delta^{-1}_+ \delta^2 = \delta_-$, (ii) $dim({\mcal
  H}_B) = (tr_M (p))^{-1} dim (L^2 (M)_M) = \delta_+$ and (iii) $B' :=
{\mcal L}_B ({\mcal H}) = M p'$. Thus, the subfactors $A \subset B'$
and $N \subset M$ are isomorphic, and so are their associated planar
algebras. The proof of Theorem \ref{sfrecon} implies that $Q$ ($=$
normalized $P$) is isomorphic to the normalized bimodule planar
algebra associated to the bimodule $_A L^2(B')_{B'}$ which, by Theorem
\ref{bimod-pa-theorem} (iii), is isomorphic to that of $_A {\mcal
  H}_B$. So, $P$ is isomorphic to the bimodule planar algebra
associated to $_A {\mcal H}_B$.
\end{pf}
\begin{rem}
Since we have an extension of Jones' theorem \cite[Theorem 4.2.1]{Jon}
as well as its converse for any finite index subfactor, we could have
very well referred a unimodular bimodule planar algebra as a
`subfactor-planar algebra', but we abstain from doing so as the term
subfactor planar algebra has already been in use for a spherical
unimodular bimodule planar algebra.
\end{rem}

\section{Examples of perturbations}\label{eg}
In this section, we will provide two examples of perturbations of
bimodule planar algebras by non-scalar weights. The first holds for
any non-spherical planar algebras whereas the second example is a more
concrete one which involves {\em diagonal subfactors}.
\subsection{Non-spherical to spherical}\label{ns2s} 
In Section \ref{pa}, we saw that the perturbation class of every
bimodule planar algebra contains a unimodular one, namely, its
normalization. In view of this, we analyze the natural question
whether there is more than one unimodular bimodule planar algebra
(upto isomorphism) in the perturbation class of a bimodule planar
algebra. The answer is negative for irreducible ones because the
weights have to be scalars.

Let $P$ be a unimodular bimodule planar algebra (not necessarily
spherical) with modulus $(\delta, \delta)$. Suppose $z \in P_{+1}$
denotes the trace intertwiner weight of $P$ (discussed in Section
\ref{bimodrecon}). Note that $z^{1/2}$ is also a weight of
$P$. Consider the planar algebra $Q:= P^{(z^{1/4},z^{1/4})}$. It is
easy to check that $Q$ is also a bimodule planar algebra. Now, for $x
\in Q_{+1} = P_{+1}$, we have
\[
Q_{TR^l_{+1}} (x)
= P_{{TR^l_{+1}}^{(z^{1/4},z^{1/4})}} (x) 
= P_{TR^l_{+1}} (z^{-1/2}x) = P_{TR^r_{+1}} (z^{1/2}x)
= P_{{TR^r_{+1}}^{(z^{1/4},z^{1/4})}} (x) 
= Q_{TR^r_{+1}} (x),
\]
where the third equality follows from the trace intertwining property
of $z$ and the rest follow directly from the definition of
perturbation. Setting $x = 1_{P_{+1}}$ in the above equation, we
conclude that $Q$ is unimodular with modulus $P_{TR^l_{+1}} (z^{-1/2})
= P_{TR^r_{+1}} (z^{1/2})$ and thereby, is the same as its
normalization; thus, $Q$ is spherical. We include this observation in
the following proposition.
\begin{prop}\label{pertclass}
The perturbation class of every bimodule planar algebra contains a
unique spherical unimodular bimodule planar algebra which can also be
characterized by the unimodular one having the minimal index value. In
other words, any bimodule planar algebra assuming the minimal index in
its perturbation class must be spherical.
\end{prop}
\begin{pf} To show minimality of the index, consider a spherical
unimodular bimodule planar algebra $P$ with modulus $\delta$ and a
positive weight $z$ of $P$. Set $Q:=P^{(z^{1/2},z^{1/2})}$. Let
$\{p_i\}^n_{i=1}$ be the set of minimal central projections of $P_{+1}
= Q_{+1}$. Note that $\delta = \us{i}{\sum} \omega_i$ where $\omega_i
= P_{TR^r_{+1}} (p_i) = P_{TR^l_{+1}} (p_i)$ for all $i$. Also, for
all $i$, there exists $\lambda_i > 0$ such that $z = \us{i}{\sum}
\lambda_i p_i$. If $(\delta_- , \delta _+)$ denotes the modulus of
$Q$, then $\delta_\pm = \us{i}{\sum} \lambda^{\pm1}_i \omega_i$. Thus,
\[\begin{array}{cll}
& index(Q) & \\
= &
\ous{n}{\sum}{i,j=1} \frac{\lambda_i}{\lambda_j} \omega_i \omega_j = \left(\ous{n}{\sum}{i=1} \omega^2_i \right) + \us{1\leq j < k \leq n}{\sum} \left( \frac{\lambda_j}{\lambda_k} + \frac{\lambda_k}{\lambda_j} \right) \omega_i \omega_j \geq \left(\ous{n}{\sum}{i} \omega^2_i \right) + \us{1\leq j < k \leq n}{\sum} 2 \omega_i \omega_j &
= \left( \us{i}{\sum} \omega_i \right)^2\\
& & = index(P)
\end{array}\]
where the equality occurs if and only if $\lambda_i = \lambda_j$ for
$1 \leq i,j \leq n$, that is, $z$ is a scalar weight. Now, if
$\tilde{P}$ is any unimodular bimodule planar algebra assuming minimal
index in its perturbation class, then there exists a weight $w$ such
that $\tilde{Q} := \tilde{P}^{(w^{1/2},w^{1/2})}$ is unimodular
spherical. Again, by the above argument, $index(\tilde{Q})$ being
minimal, is same as $index(\tilde{P})$ and $w$ is a scalar
weight. Since both $\tilde{P}$ and $\tilde{Q}$ are unimodular,
therefore $w=1_{P_{+1}}$ and hence, $\tilde{P} = \tilde{Q}$ is
spherical. This also shows the uniqueness of a spherical unimodular
bimodule planar algebra in the perturbation class of a bimodule planar
algebra.\end{pf}
\begin{rem}
Instance of minimizing index of a conditional expectation
onto a subfactor already appeared in the work of Hiai (see \cite{Hia88}) and
then Popa (see \cite{Pop94}). Proposition \ref{pertclass} gives a
nice way of minimizing index using perturbation of planar algebra.
\end{rem}
\begin{rem}
From the proof of Proposition \ref{pertclass}, it is clear that if the
perturbation of a spherical bimodule planar algebra by a positive
weight, is spherical, then the weight must be a scalar one.
\end{rem}
\begin{rem}
It is easy to check that a finite depth bimodule planar algebra is
always spherical because by Perron-Frobenius theorem, the index must
be equal to the norm-square of the pricipal graph and perturbation
does not change the principal graphs.
\end{rem}
\subsection{Spherical to non-spherical}\label{s2ns}
Here we try to find whether we can perturb one of the known spherical
unimodular planar algebras and get a non-spherical one. For this, we
study the case of diagonal subfactors.

In Section \ref{ns2s}, we found that the perturbation of a bimodule
planar algebra $P$ by a weight $z$, has modulus $(P_{TR^l_{+1}}
(z^{-1/2}) , P_{TR^r_{+1}} (z^{1/2}))$; so, for the perturbation to be
unimodular, we need to find $z$ satisfying $P_{TR^l_{+1}} (z^{-1/2}) =
P_{TR^r_{+1}} (z^{1/2})$. Consider the diagonal planar algebra $P$
constructed in \cite{BDG08} with respect to the free group ${\mathbb
  F}_2$ generated by two free generators $a_\eta$ for $\eta \in
I:=\{-,+\}$ and the trivial cocycle. We briefly recall (from
\cite{BDG08}) few aspects of this planar algebra $P$ which will be
needed for further analysis.
\vspace*{2 mm}\\ {\em Vector spaces:} $P_{\eta 0}:= \C$ and $P_{\eta k}
= \C \left\{ \uline{\vlon} \in I^{2k} : alt_\eta (\uline{\vlon}) = e
\right\}$ for all $k \geq 1$, where $e$ denotes the identity of
${\mathbb F}_2$ and $alt_\eta$ is given by $I^n \ni (\vlon_1, \cdots,
\vlon_n) \oset{alt_\eta}{\longmapsto} \left\{
\begin{array}{ll}
a^{-1}_{\vlon_1} a_{\vlon_2} a^{-1}_{\vlon_3} \cdots
a^{(-1)^n}_{\vlon_n} \in {\mathbb F}_2 & \text{if }
\eta=+,\\ a_{\vlon_1} a^{-1}_{\vlon_2} a_{\vlon_3} \cdots
a^{(-1)^{n-1}}_{\vlon_n} \in {\mathbb F}_2& \text{if } \eta=-.
\end{array}\right.$
\vspace{3 mm}\\ 
{\em Action of tangles:} For a tangle $T : (\eta_1
k_1, \cdots, \eta_b k_b) \rightarrow \eta_0 k_0$, the action $P_T :
P_{\eta_1 k_1}\times \cdots \times P_{\eta_b k_b} \rightarrow
P_{\eta_0 k_0}$ is given by
\[
\left\langle P_T(\uline{\vlon}^1, \cdots, \uline{\vlon}^b),
\uline{\vlon}^0 \right\rangle := \# \left\{ f:\{\text{strings in }T\}
\rightarrow I \left| \left. f \right| _{\text{marked points of $i$-th
    disc in $T$}} = \uline{\vlon}^i \text{ for } 0\leq i \leq b
\right.  \right\}
\]
where $\uline{\vlon}^i$ belongs to the distinguished basis of
$P_{\eta_i k_i}$ for $0\leq i \leq b$.

Fix a positive scalar $\lambda_{-} = \lambda^{-1}_+$. Consider $z :=
\left[\lambda_- (-,-) + \lambda_+ (+,+)\right] \in P_{+1}$ and
$\lambda_{\eta \uline{\vlon}} := \oset{k}{\uset{i=1}{\prod}}
\lambda^{\eta (-1)^{i-1}}_{\vlon_i}$ for $\uline{\vlon} = (\vlon_1,
\cdots, \vlon_k) \in I^{k}$.  Using the relation $z^{-1} =
\left[\lambda^{-1}_- (-,-) + \lambda^{-1}_+ (+,+)\right]$ $ \in
P_{+1}$ and the action of tangles in Definition \ref{weight}, one can
derive that $z_{\eta k} = \uset{\uline{\vlon} \in I^{k}}{\sum}
\lambda_{\eta \uline{\vlon}} \; (\uline{\vlon}, \tilde{\uline{\vlon}})
\in P_{\eta k}$, where $\tilde{\uline{\vlon}}$ is the sequence
obtained by reversing the order in $\uline{\vlon}$.  From the action
of the multiplication tangle $M_{\eta k}$, we have $z_{\eta k} \;
\left(\uline{\vlon} , \uline{\tilde{\nu}}\right) \; z^{-1}_{\eta k} =
\lambda_{\eta \uline{\vlon}} \; \lambda^{-1}_{\eta \uline{{\nu}}} \;
\left(\uline{\vlon} , \uline{\tilde{\nu}} \right) = \lambda_{\eta
  \left(\uline{\vlon} , \uline{\tilde{\nu}}\right)} \;
\left(\uline{\vlon} , \uline{\tilde{\nu}} \right)$ for $\uline{\vlon}
, \uline{\nu} \in I^k$ such that $alt_\eta \left(\uline{\vlon} ,
\uline{\tilde{\nu}}\right) = e$.  Freeness of $a_-$ and $a_+$ implies
that there exists a configuration of non-crossing pairings of matching
signs (abbreviated as `NC-pairing') in the sequence
$\left(\uline{\vlon} , \uline{\tilde{\nu}}\right)$, which implies that
$\lambda_{\eta \left(\uline{\vlon} , \uline{\tilde{\nu}}\right)} = 1$.
Thus, $z_{\eta k}$ is central and hence, $z$ is a weight of $P$.
Further, positivity of $z$ implies that $Q := P^{(z^{1/2},z^{1/2})}$
is a bimodule planar algebra. To check unimodularity of $Q$, note that
$P_{TR^l_{+1}} (z^{-1}) = \left(\lambda_-\right)^{-1} +
\left(\lambda_+\right)^{-1} = \lambda_- + \lambda_+ = P_{TR^r_{+1}}
(z)$. Now, $Q$ is spherical if and only if, for all $s , t \in \C$,
\[
s \lambda_-^{-1} + t \lambda_+^{-1} = Q_{TR^l_{+1}} \left(s (-,-) + t
(+,+) \right) = Q_{TR^r_{+1}} (s (-,-) + t (+,+)) = s \lambda_- + t
\lambda_+
\]
if and only if $\lambda_- = 1 = \lambda_+$, that is, $z =
1_{P_{+1}}$. Clearly, the range of the index values of these perturbed
planar algebras is $[4,\infty)$.  We gather the above observation in
  the following proposition.
\begin{prop}\label{pertdiagpa}
Let $P$ denote the diagonal planar algebra constructed in \cite[\S
  3]{BDG08} with respect to the free group ${\mathbb F}_2$ generated
by two free generators and the trivial cocycle. Then one can perturb
$P$ to obtain a continuum of unimodular bimodule planar algebras
$\{Q^\gamma\}_{\gamma \geq 4}$ such that $Q^\gamma$ has (i) index
$\gamma$, (ii) is equal to $P$ if $\gamma =4$ and (iii) is
non-spherical if and only if $\gamma \neq 4$.
\end{prop}
\begin{rem}\label{jones-remark}
These perturbed planar algebras turn out to be the ones associated to
the subfactors with index greater than $4$ constructed in \cite{Jon83}
(which are all non-extremal), as was pointed out by Vaughan Jones.
\end{rem}
We give an explicit proof of the above remark in the following
proposition.
\begin{prop}\label{jones-prop}
Let $N\subset M$ be a subfactor of type $II_{1}$ with a partition of
unity $\left\{ \left. p_{\vlon}\in\mscr{P}(N^{\prime}\cap M) \,
\right| \vlon\in I \right\} $ satisfying (i)
$Np_{\vlon}=M_{p_{\vlon}}$ and (ii) $c_{-}:=tr_{M}(p_{-})\neq
tr_{M}(p_{+})=:c_{+}$. (Such subfactors were constructed by Jones in
\cite{Jon83} with $M$ having full fundamental group.) Then, the planar
algebra associated to $N\subset M$ is isomorphic to $Q^{(c_-
  c_+)^{-1}}$ (as in Proposition \ref{pertdiagpa}).
\end{prop}
\begin{pf}
Condtions $(i)$ and $(ii)$ imply
$p_{\vlon}\in\mscr{P}_{minl}(N^{\prime}\cap M)$, $N^{\prime}\cap
M\cong\mathbb{C}p_{-}\oplus\mathbb{C}p_{+}$ and
$c_{\vlon}tr_{N'}(p_{\vlon})[M:N]=[M_{p_{\vlon}}:Np_{\vlon}]=1$.  So,
$[M:N]=c_{-}^{-1}+c_{+}^{-1}=(c_{-}c_{+})^{-1}$ and
$tr_{N'}(p_{\vlon})=c_{-\vlon}$ and hence, $N\subset M$ is
non-extremal.

Suppose $\{M_{k}\}_{k \geq 1}$ is a tower of basic constructions of $N
\subset M$ with Jones projections $\{e_{k}\}_{k \geq 1}$ and
$\tilde{P}$ be the planar algebra associated to this tower.  Let
$p_{\eta}^{(n)}$ denote the element $P_{\psfrag{n-1}{$n-1$}
  \psfrag{+}{$+$} \psfrag{peta}{$p_\eta$}
  \includegraphics[scale=0.2]{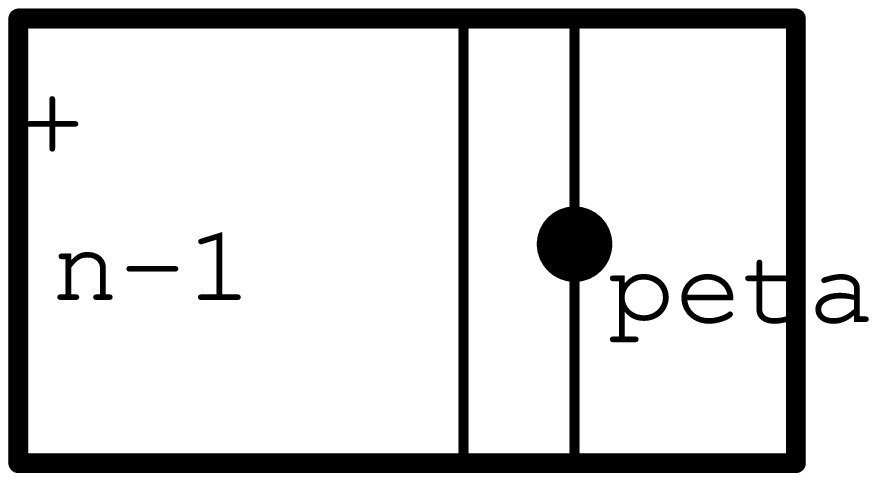}}$ or $P_{
  \psfrag{n-1}{$n-1$} \psfrag{+}{$+$} \psfrag{peta}{$p_\eta$}
  \includegraphics[scale=0.2]{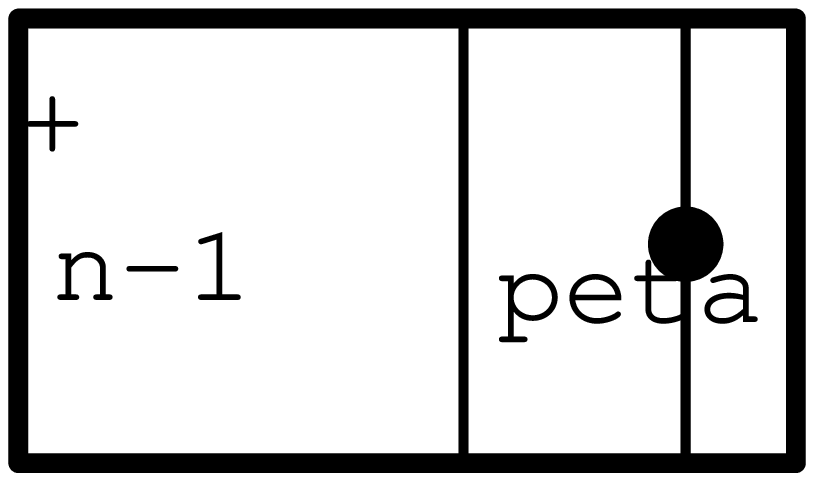}}$ (lying in
$\mscr{P}(M_{n-2}^{\prime}\cap M_{n-1})$ by Corollary \ref{dualpa}(b)) according as $n$ is even or odd,
$c_{\eta}^{(n)}:=tr_{M_{n-1}}(p_{\eta}^{(n)}) = c_{(-)^{n-1}\eta}$,
$p_{\underline{\vlon}}:=p_{\vlon_{1}}^{(1)}\cdots
p_{\vlon_{n}}^{(n)}\in\mscr{P}(N^{\prime}\cap M_{n-1})$ and
$c_{\ul{\vlon}}:=\ous
{n}{\prod}{i=1}c_{\vlon_i}^{(i)}=tr_{M_{n-1}}(p_{\ul{\vlon}})$ for
$\eta \in I$, $\underline{\vlon}=(\vlon_{1},\cdots,\vlon_{n})\in
I^{n}$.  Clearly, $\{p_{\underline{\vlon}} : \underline{\vlon}\in
I^{n}\}$ forms a partition of unity. Further, the formula involving
indices of cut-down subfactors implies $Np_{\underline{\vlon}} =
\left(M_{n-1}\right)_{p_{\underline{\vlon}}}$ and hence
$p_{\ul{\vlon}}$ is minimal in $\tilde{P}_{n}$. The set of
$tr_{M_{n-1}}$-values of these minimal projections is given by
$\left\{ \left.c(n,k):=c_{-}^{k}c_{+}^{n-k}\right|0\leq k\leq
n\right\} $.

In order to establish an isomorphism between $\tilde{P}$ and
$Q^{(c_-c_+)^{-1}}$, we first prove that for $\ul{\vlon},\ul{\eta}\in
I^{n}$, $c_{\ul{\vlon}}=c_{\ul{\eta}}$ if and only if there exists an
NC-pairing for the sequence
$\left(\ul{\vlon},\tilde{\ul{\eta}}\right)$.

\noindent {\em `if' part}: The pairings give the equation
$c_{\left(\ul{\vlon},\tilde{\ul{\eta}}\right)}=\left(c_{-}c_{+}\right)^{n}$.
On the other hand, $c_{\left(\ul{\vlon},\tilde{\ul{\eta}}\right)} =
c_{\ul{\vlon}}c_{-\ul{\eta}}=c_{\ul{\vlon}}
\left[\frac{\left(c_{-}c_{+}\right)^{n}}{c_{\ul{\eta}}}\right]$ where
$-(\eta_{1},\ldots,\eta_{n}) := (-\eta_{1},\ldots,-\eta_{n})$.

\noindent {\em `only if' part}: We use induction on the length $n.$
The intial case of $n=1$ is a triviality since $c_{-}\neq c_{+}$.  In
the induction step, let $n>1$ and $0\leq k \leq n$ such that
$c_{\ul{\vlon}}=c(n,k)=c_{\ul{\eta}}$. We may also assume $0<k<n$
without loss of generality because there is exactly one sequence in
each of the cases $k=0$ and $k=n$.  Note that if both
$\left(\ul{\vlon},\tilde{\ul{\eta}}\right)$ and
$\left(\ul{\eta},\tilde{\ul{\nu}}\right)$ have NC-pairings, then so
does $\left(\ul{\vlon},\tilde{\ul{\nu}}\right)$ (resp.,
$\left(\ul{\eta},\tilde{\ul{\vlon}}\right)$) which can easily be
obtained by taking the product (resp., $\ast$) of the Temperley-Lieb
diagrams associated to the NC-pairings.  Consider
$\underline{{\vlon}}(n,k) \in I^n$ given by
$\vlon_{i}(n,k)=\left\{ \begin{array}{ll} (-)^{i} & \text{if }1\leq
  i\leq k, \\ (-)^{i+1} & \text{if }k+1\leq i\leq
  n\end{array}\right.$.  Clearly, $c_{\ul{\vlon}(n,k)}=c(n,k)$; so, it
  is enough to show that
  $\left(\underline{\vlon},\underline{\tilde{\vlon}}(n,k)\right)$ has
  an NC-pairing. We begin by showing that there exists a $\ul{\nu}\in
  I^{n}$ such that $\nu_{1}=-=\vlon_{1}(n,k)$ and
  $\left(\ul{\vlon},\tilde{\ul{\nu}}\right)$ has an NC-pairing (and
  hence $c_{\ul{\vlon}}=c_{\ul{\nu}}$ by the `if' part). Suppose first
  that $\vlon_1 = +$. Since $0<k<n$, $\underline{\vlon}$ must have a
  pair of consecutive matching signs. Choose the left most consecutive
  matching pair in $\ul{\vlon}$ and change both the signs in the pair
  to obtain $\ul{\gamma}\in I^{n}$.  Observe that there exists an
  NC-pairing for $\left(\ul{\vlon},\tilde{\ul{\gamma}}\right)$.  Also,
  the first consecutive matching pair in $\ul{\gamma}$ gets closer to
  the left side than the one in $\ul{\vlon}$. Apply this method
  succesively to obtain the desired $\ul{\nu}$. And, if $\vlon_1 = -$,
  then simply take $\uline{\nu} = \uline{\vlon}$. Let
  $\ul{\nu}^\prime$ (resp., $\ul{\vlon}^\prime (n,k)$) be the sequence
  obtained from $\ul{\nu}$ (resp., $\ul{\vlon} (n,k)$) after removing
  the minus sign in the first entry.  Since $c_{\uline{\nu'}} =
  c_{\uline{\vlon'}(n,k)}$, by induction assumption, there exists an
  NC-pairing for $\left(\underline{\vlon}^\prime (n,k),
  \tilde{\underline{\nu}^\prime} \right)$ and hence for
  $\left(\underline{\vlon}(n,k),\tilde{\underline{\nu}}\right)$.
\vspace*{1mm}

It is also true that for $\ul{\vlon},\ul{\eta}\in I^{n}$,
$c_{\ul{\vlon}}=c_{\ul{\eta}}$ if and only if $p_{\ul{\vlon}}$ and
$p_{\ul{\eta}}$ are equivalent in $\tilde{P}_{n}$. Suppose
$c_{\ul{\vlon}}=c_{\ul{\eta}}$. So, there exists an NC-pairing for
$\left(\ul{\vlon},\tilde{\ul{\eta}}\right)$.  Let $T$ be the
Temperley-Lieb diagram associated to this NC-pairing.  Since
$p_{\ul{\vlon}}$ and $p_{\ul{\eta}}$ are minimal in $\tilde{P}_{n}$,
therefore it is enough to check
$p_{\ul{\vlon}}\left(\tilde{P}_{T}\right)p_{\ul{\eta}}\neq0$ which
easily follows by considering the scalar
$tr_{M_{n-1}}\left(p_{\ul{\vlon}}
\left(\tilde{P}_{T}\right)p_{\ul{\eta}}\left(\tilde{P}_{T^{*}}\right)
p_{\ul{\vlon}}\right)$ and then using the action of tangles to show it
is non-zero (in fact, positive). In particular,
$p(\uline{\vlon},T,\tilde{\uline{\eta}}):=
\sqrt{\frac{c_{\ul{\vlon}}}{tr_{M_{n-1}}\left(p_{\ul{\vlon}}
    \left(\tilde{P}_{T}\right)p_{\ul{\eta}}\left(\tilde{P}_{T^{*}}
    \right)\right)}}p_{\ul{\vlon}}\left(\tilde{P}_{T}\right)p_{\ul{\eta}}$
is a partial isometry with $p_{\ul{\vlon}}$ (resp., $p_{\ul{\eta}}$)
as its final (resp., initial) projection. To show that
$p(\uline{\vlon},T,\tilde{\uline{\eta}})$ is independent of $T$,
consider two Temperley-Lieb elements $T_{1}$ and $T_{2}$ for two
NC-pairings of $(\uline{\vlon}, \uline{\tilde{\eta}})$. It immediately
follows from the action of tangles that $p(\uline{\vlon},
T_1,\tilde{\uline{\eta}})\left(p(\uline{\vlon},T_2,\tilde{\uline{\eta}})
\right)^\ast \in \R_+ \left[ p_{\ul{\vlon}}
  \left(\tilde{P}_{T_{1}}\right)p_{\ul{\eta}}\left(\tilde{P}_{T_{2}^{*}}
  \right)p_{\ul{\vlon}} \right]= \R_+
\left[p_{\ul{\vlon}}\left(\tilde{P}_{T_{1}\cdot
    T_{2}^{*}}\right)p_{\ul{\vlon}}\right] = \R_+
p_{\ul{\vlon}}$. Hence, $p(\uline{\vlon},T_1,\tilde{\uline{\eta}}) =
p(\uline{\vlon},T_2,\tilde{\uline{\eta}})$. So, we will denote the
partial isometry simply by $p(\uline{\vlon},\tilde{\uline{\eta}})$.

Consider the diagonal planar algebra $P$ in Proposition
\ref{pertdiagpa} and its perturbation $Q=Q^{(c_-c_+)^{-1}} =
P^{(z^{1/2}, z^{1/2})}$ with respect to the weight $z :=
\sqrt{\frac{c_-}{c_+}}(-,-)+\sqrt{\frac{c_+}{c_-}}(+,+)$. Define
$\vphi:Q_+ \rightarrow \tilde{P}_{+}$ by
\[
Q_{+k} = P_{+k} =
(\uline{\vlon},\tilde{\uline{\eta}})\os{\vphi}{\mapsto}
p(\uline{\vlon},\tilde{\uline{\eta}}) \in \tilde{P}_{+k}.
\]
Clearly, $\vphi$ is a $\ast$-algebra isomorphism. By Proposition
\ref{pa-positive-isom}, it is enough to show that $\vphi$ is
equivarent with respect to all tangles having positive colors on their
discs and in fact, for the following tangles:

\noindent {\em Right-inclusion tangles:} For $\uline{\vlon} ,
\uline{\eta} \in I^n$ such that $alt_+(\uline{\vlon} ,
\tilde{\uline{\eta}}) = e$ (equivalently, $(\uline{\vlon} ,
\tilde{\uline{\eta}})$ has an NC-pairing), we have $\vphi \circ
Q_{RI_{+ n}} (\uline{\vlon} , \tilde{\uline{\eta}}) = p(\uline{\vlon}
, -, -, \tilde{\uline{\eta}}) + p (\uline{\vlon} , +, +,
\tilde{\uline{\eta}}) \in \tilde{P}_{+(n+1)}$. Thus, it is enough to
show that $p(\uline{\vlon} , \tilde{\uline{\eta}}) = p(\uline{\vlon} ,
-, -, \tilde{\uline{\eta}}) + p (\uline{\vlon} , +, +,
\tilde{\uline{\eta}})$.  Let $T$ be the Temperley-Lieb diagram
corresponding to an NC-pairing for $(\uline{\vlon} ,
\tilde{\uline{\eta}})$; then, $RI_{+n} \circ T$ gives an NC-pairing
for $(\uline{\vlon} ,\pm,\pm, \tilde{\uline{\eta}})$. Further, we have
$c_{(\uline{\vlon},\pm)} = c_{\uline{\vlon}} \, c^{(n+1)}_\pm$ and
\[tr_{M_n}\left( p_{(\ul{\vlon},\pm)} \left( \tilde{P}_{RI_{+n} \circ
  T} \right) p_{(\ul{\eta},\pm)} \left( \tilde{P}_{(RI_{+n} \circ
  T)^{*}} \right) \right) =
tr_{M_{n-1}}\left(p_{\ul{\vlon}}\left(\tilde{P}_{T}\right)p_{\ul{\eta}}
\left(\tilde{P}_{T^{*}}\right)\right) c^{(n+1)}_\pm\] which imply
$p(\uline{\vlon} ,\pm,\pm, \tilde{\uline{\eta}}) = p(\uline{\vlon} ,
\tilde{\uline{\eta}}) p^{(n+1)}_\pm$. Hence, the required equality
holds.

\noindent{\em Jones projection tangles:} From the action of tangles on
$\tilde{P}$, we have
\begin{align*}
\tilde{P}_{E_{+n}} =& \us{\uline{\vlon}^\prime,\uline{\vlon}^{\prime
    \prime} \in I^{n+1}}{\sum} p_{\uline{\vlon}^\prime}
\tilde{P}_{E_{+n}} p_{\uline{\vlon}^{\prime \prime}} =
\us{\uline{\vlon} \in I^{n-1}, \eta,\nu \in I}{\sum}
p_{(\uline{\vlon},\eta,\eta)} \tilde{P}_{E_{+n}}
p_{(\uline{\vlon},\nu,\nu)}\\ = & \us{\uline{\vlon} \in I^{n-1},
  \eta,\nu \in I}{\sum}
\sqrt{\frac{tr_{M_{n}}\left(p_{(\ul{\vlon},\eta,\eta)}\left(\tilde{P}_{E_{+n}}
    \right)p_{(\ul{\vlon},\nu,\nu)}\left(\tilde{P}_{E_{+n}}
    \right)\right)}{c_{(\ul{\vlon},\eta,\eta)}}} \; p (\uline{\vlon},
\eta,\eta,\nu,\nu,\tilde{\uline{\vlon}}).
\end{align*}
Now, for $\uline{\vlon} \in I^{n-1}$ and $\eta,\nu \in I$, note that
$\left(\tilde{P}_{E_{+n}}\right)p_{(\ul{\vlon},\nu,\nu)}\left(\tilde{P}_{E_{+n}}\right)
= \frac{c^{(n)}_{\nu}}{\sqrt{c_- c_+}} \; p_{\ul{\vlon}}
\left(\tilde{P}_{E_{+n}}\right)$ and
\begin{align*}
tr_{M_n}\left(p_{(\ul{\vlon},\eta,\eta)} p_{\ul{\vlon}}
\left(\tilde{P}_{E_{+n}} \right) \right) & = \sqrt{c_- c_+}\,
tr_{M_n}\left( \left(\tilde{P}_{E_{+n}}\right)
p_{(\ul{\vlon},\eta,\eta)} \left(\tilde{P}_{E_{+n}}\right)\right)\\ &
= c^{(n)}_{\eta} tr_{M_n}\left( p_{\ul{\vlon}}
\left(\tilde{P}_{E_{+n}}\right) \right) = c^{(n)}_{\eta} \sqrt{c_-
  c_+} \, tr_{M_{n-2}}( p_{\ul{\vlon}}) = c^{(n)}_{\eta}
c_{\ul{\vlon}} \sqrt{c_- c_+}
\end{align*}
Since $c_{(\ul{\vlon}, \eta, \eta)} =
c_- c_+ c_{\ul{\vlon}}$, therefore $\tilde{P}_{E_{+n}} =
\us{\uline{\vlon} \in I^{n-1}, \eta, \nu \in I}{\sum}
\sqrt{\frac{c^{(n)}_{\eta} c^{(n)}_{\nu}}{c_- c_+}} p (\uline{\vlon},
\eta,\eta,\nu,\nu,\tilde{\uline{\vlon}})$.  On the other hand, using
the weight $z$ of $P$, $Q_{E_{+n}}$ is equal to
\[
\left[\us{\ul{\vlon} \in I^{n-1}}{\sum}
\sqrt{\frac{c_{(-)^{n-1}-}}{c_{(-)^{n-1}+}}}
(\uline{\vlon},-,-,-,-, \tilde{\uline{\vlon}}) +
(\uline{\vlon},-,-,+,+, \tilde{\uline{\vlon}}) +
(\uline{\vlon},+,+,-,-, \tilde{\uline{\vlon}}) +
\sqrt{\frac{c_{(-)^{n-1}+}}{c_{(-)^{n-1}-}}}
(\uline{\vlon},+,+,+,+, \tilde{\uline{\vlon}})
\right].\] Thus, $\vphi$ preserves the action of $E_{+n}$.
\vspace*{1mm}

\noindent{\em Left conditional expectation tangle:} Fix
$\nu,\nu^\prime \in I$ and $\ul{\vlon}, \ul{\eta} \in I^{n-1}$ such
that $alt_+ (\nu, \ul{\vlon}, \tilde{\ul{\eta}}, \nu^\prime) =
e$. Note that $Q_{E'_{+n}} (\nu, \ul{\vlon}, \tilde{\ul{\eta}},
\nu^\prime) = \delta_{\nu = \nu^\prime} Q_{E'_{+n}} (\nu, \ul{\vlon},
\tilde{\ul{\eta}}, \nu) = \delta_{\nu = \nu^\prime}
\sqrt{\frac{c_{-\nu}}{c_\nu}} \us{\gamma \in I}{\sum} (\gamma,
\ul{\vlon}, \tilde{\ul{\eta}}, \gamma) \os{\vphi}{\mapsto} \delta_{\nu
  = \nu^\prime} \sqrt{\frac{c_{-\nu}}{c_\nu}} \us{\gamma \in I}{\sum}
p (\gamma, \ul{\vlon}, \tilde{\ul{\eta}}, \gamma)$.  On the other
hand, $\tilde{P}_{E'_{+n}} (p (\nu, \ul{\vlon}, \tilde{\ul{\eta}},
\nu^\prime)) = \tilde{P}_{E'_{+n}} (p_\nu \, p (\nu, \ul{\vlon},
\tilde{\ul{\eta}}, \nu^\prime) \, p_{\nu^\prime}) = \delta_{\nu =
  \nu^\prime} \tilde{P}_{E'_{+n}} ( p (\nu, \ul{\vlon},
\tilde{\ul{\eta}}, \nu))$.  Now, we may assume that $alt_+ (\nu,
\ul{\vlon}, \tilde{\ul{\eta}}, \nu) = e$ which implies $alt_-
(\ul{\vlon}, \tilde{\ul{\eta}}) = e$ equivalently, $(\ul{\vlon},
\tilde{\ul{\eta}})$ has an NC-pairing given by a Temperley-Lieb diagram
$T$ (say).  Then, $T_1 := LI_{-(n-1)}\circ T$ induces an NC-pairing for
$(\nu, \ul{\vlon}, \tilde{\ul{\eta}}, \nu)$. So, we have
\begin{align*}
\tilde{P}_{E'_{+n}} ( p (\nu, \ul{\vlon}, \tilde{\ul{\eta}}, \nu)) &=
\sqrt{\frac{c_{(\nu,\ul{\vlon})}}{tr_{M_{n-1}}\left(p_{(\nu,
      \ul{\vlon})} \left(\tilde{P}_{T_1}\right)p_{(\nu, \ul{\eta})}
    \left(\tilde{P}_{T^{*}_1}\right)\right)}} \tilde{P}_{E'_{+n}}
\left( p_{(\nu, \ul{\vlon})} \left( \tilde{P}_{T_1} \right) p_{(\nu ,
  {\ul{\eta}})} \right)\\ &=
\sqrt{\frac{c_{(\nu,\ul{\vlon})}}{tr_{M_{n-1}}\left(p_{(\nu,
      \ul{\vlon})} \left(\tilde{P}_{T_1}\right)p_{(\nu, \ul{\eta})}
    \left(\tilde{P}_{T^{*}_1}\right)\right)}}
\sqrt{\frac{c_{-\nu}}{c_\nu}} \us{\gamma \in I}{\sum} \; p_{(\gamma,
  \ul{\vlon})} \left( \tilde{P}_{T_1} \right) p_{(\gamma ,
  {\ul{\eta}})}.
\end{align*}
The last equality is obtained by using the action of tangles on
$\tilde{P}$.  To show $\vphi \circ Q_{E'_{+n}} = \tilde{P}_{E'_{+n}}
\circ \vphi$, it remains to be shown that
$\frac{c_{(\nu,\ul{\vlon})}}{tr_{M_{n-1}}\left(p_{(\nu, \ul{\vlon})}
  \left(\tilde{P}_{T_1}\right)p_{(\nu, \ul{\eta})}
  \left(\tilde{P}_{T^{*}_1}\right)\right)}$ is independent of $\nu \in
I$. Observe that the denominator
\begin{align*}
& tr_{M_{n-1}}\left(p_{(\nu, \ul{\vlon})}
  \left(\tilde{P}_{T_1}\right)p_{(\nu, \ul{\eta})}
  \left(\tilde{P}_{T^{*}_1}\right)\right) = (c_- c_+)^{\frac{n}{2}}
  \tilde{P}_{TR^r_{+n}} \left(p_{(\nu, \ul{\vlon})}
  \left(\tilde{P}_{T_1}\right)p_{(\nu, \ul{\eta})}
  \left(\tilde{P}_{T^{*}_1}\right)\right)\\ =& (c_- c_+)^{\frac{n}{2}}
  \sqrt{\frac{c_\nu}{c_{-\nu}}} \sqrt{c_- c_+} \tilde{P}_{TR^r_{+n}}
  \left(p^\prime_{\ul{\vlon}}
  \left(\tilde{P}_{T_1}\right)p^\prime_{\ul{\eta}}
  \left(\tilde{P}_{T^{*}_1}\right)\right) = c_\nu (c_-
  c_+)^{\frac{n}{2}} \tilde{P}_{TR^r_{+n}} \left(p^\prime_{\ul{\vlon}}
  \left(\tilde{P}_{T_1}\right)p^\prime_{\ul{\eta}}
  \left(\tilde{P}_{T^{*}_1}\right)\right),
\end{align*}
where $p^\prime_{(\omega_1,\cdots,\omega_{n-1})} :=
p^{(2)}_{\omega_1} p^{(3)}_{\omega_2} \cdots p^{(n)}_{\omega_{n-1}}
\cdots, \omega_{n-1}) \in I^{n-1}$ and in the second equality, we
replace the loop with $p_\nu$ in it by a loop with nothing on it, and
the numerator $c_{(\nu,\ul{\vlon})} = c_\nu c_{-\ul{\vlon}}$.  Hence,
the required fraction becomes independent of $\nu \in I$.
\end{pf}


\section{Questions}
\begin{defn}
A bimodule planar algebra is said to have trivial perturbation class
if all its perturbation by positive weights, are spherical.
\end{defn}
It will be interesting to find a set of necessary and sufficient
conditions for a bimodule planar algebra having trivial perturbation
class in terms of its principal graph(s); one can also consider this
question in the more specific case of the Bisch-Haagerup planar
algebras (see \cite{BH96}, \cite{BDG09} and \cite{BDG10}). Note that
all finite depth or irreducible bimodule planar algebras have trivial
perturbation class; so, for the Bisch Haagerup planar algebras, this
question is relevant ony when the two subgroups have nontrivial
intersection and the group generated by them is infinite.
\vspace{2mm}

Another interesting problem is to obtain a method of perturbing a
bifinite bimodule to a new one whose associated planar algebra is the
perturbation of the one associated to the bimodule which we start
with. Note that, by Theorem \ref{bmrecon}, we may find a bifinite
bimodule corresponding to a perturbation of the planar algebra
associated to a given bifinite bimodule but we don't have any direct
relation between these bimodules in terms of the weight.
We will give answers to these questions in a forthcoming article.

\subsection*{Acknowledgements} The authors would like to thank Stefaan
Vaes for his helpful comments, questions and views which helped us to
be more precise at many instances. Thanks are also due to Steven
Deprez, Vijay Kodiyalam and V. S. Sunder for several useful
discussions. The authors would like to thank Vaughan Jones for
suggesting Remark \ref{jones-remark} to us during the NCGOA 2010
conference at Vanderbilt University.

\bibliographystyle{alpha}

\vspace*{8mm}
\end{document}